\documentclass[12pt]{amsart}
\usepackage{epsf,amscd,amssymb}
 
\newtheorem{thm}{Theorem}[section]
\newtheorem{lem}[thm]{Lemma}
\newtheorem{cor}[thm]{Corollary}
\newtheorem{prop}[thm]{Proposition}

\newtheorem{thmA}{Theorem A}

\newtheorem{thmB}{Theorem B}

\newtheorem{thmC}{Theorem C}

\newtheorem{corD}{Corollary D}

\theoremstyle{remark} 
\newtheorem{rem}[thm]{Remark}
\newtheorem{hyp}[thm]{Hypothesis}

\newtheorem{app}[thm]{State}

\theoremstyle{definition}
\newtheorem{defn}[thm]{Definition}

\newcommand\bR{{\mathbb{R}}}
\newcommand\bC{{\mathbb C}}
\newcommand\bE{{\mathbb E}}

\newcommand\bH{{\mathbb{H}}}

\newcommand\PO{{\rm PO}}
\newcommand\PSL{{\rm PSL}}

\newcommand\dev{{\bf dev}}
\newcommand\SI{{\bf S}}

\newcommand\clo{{\rm Cl}}

\newcommand\ra{\rightarrow}

\newcommand\emp{\emptyset}
\newcommand\eps{\epsilon}

\newcommand\vth{\vartheta}

\newcommand\ovl{\overline}

\newcommand\calR{{\mathcal R}}
\newcommand\calS{{\mathcal S}}
\newcommand\calT{{\mathcal T}}
\newcommand\calC{{\mathcal C}}
\newcommand\calI{{\mathcal I}}
\newcommand\calB{{\mathcal B}}

\newcommand\hyps{\mathbb{H}^3} 
\newcommand\hypp{\mathbb{H}^2}
\newcommand\tri{\triangle}

\newcommand\dist{\mathbb{d}}
\newcommand\Isom{\rm{Isom}}
\newcommand\Area{\rm{Area}}
\newcommand\convh{\rm{convh}}

\newcommand\sing{\rm sing}
\newcommand\core{\mathcal C}
\newcommand\cl{\bf c}

\newcommand\cat{{\rm{CAT}}}

\newcommand\rmc{\hbox{\rm :}}

\begin{document}

\title[The PL-methods for hyperbolic $3$-manifolds]
{The PL-methods for hyperbolic $3$-manifolds to prove tameness}
\author{Suhyoung Choi}
\address{Department of Mathematics \\
Korea Advanced Institute of Science and Technology\\  
305--701 Daejeon, South Korea}
\email{schoi@math.kaist.ac.kr}
\date{January 10, 2006 (version 3)} 
\thanks{The author gratefully acknowledges  
support from the Korea Reseach Foundation Grant (KRF-2002-070-C00010)}

\begin{abstract}
Using PL-methods, we prove the Marden's conjecture that a hyperbolic $3$-manifold $M$ 
with finitely generated fundamental group and with no parabolics 
are topologically tame. Our approach is to form 
an exhaustion $M_i$ of $M$ and modify the 
boundary to make them $2$-convex. 
We use the induced path-metric, which  
makes the submanifold $M_i$ negatively curved 
and with Margulis constant independent of $i$.
By taking the convex hull in the cover of $M_i$ 
corresponding to the core, we show that 
there exists an exiting sequence of surfaces $\Sigma_i$.
Some of the ideas follow those of Agol. 
We drill out the covers of $M_i$ by a core $\core$ again 
to make it negatively curved.
Then the boundary of the convex hull of $\Sigma_i$ is shown 
to meet the core. By the compactness argument of Souto, 
we show that infinitely many of $\Sigma_i$ are homotopic in $M - \core^o$. 
Our method should generalize to a more wider class of piecewise hyperbolic manifolds. 
\end{abstract}  

\maketitle

\tableofcontents

%This paper is a modification of the paper titled 
%``Drilling cores of hyperbolic $3$-manifolds and tameness"
%which is expanded into three papers. 
%In the first paper \cite{2convhull}, 
%titled ``$2$-convex hulls of  hyperbolic manifolds"
%we showed how to isotopy the boundary of a submanifold 
%of codimension zero so that the submanifold becomes 
%$2$-convex. In the second paper \cite{genhypmfld} titled
%``General hyperbolic $3$-manifolds and convex hulls of their cores"
%we showed that $2$-convex hyperbolic manifolds are Gromov hyperbolic, 
%we discussed the Gauss-Bonnet type estimates for the area of hyperbolic 
%surfaces, and the boundary components of 
%convex hulls of the cores of general hyperbolic
%manifolds. The aim of this paper is to prove the Marden's conjecture 
%using these methods. 

%** motivation
Recently, Agol initiated a really interesting 
approach to proving Marden's conjecture
by drilling out closed geodesics. 
In this paper, we use truncation of the hyperbolic 
manifolds and drilling out compact cores of them
to prove the tameness conjecture. 
We do use the Agol's idea of using covering spaces and 
taking convex hulls of the cores. 
Our approach is somewhat different
in that we do not use end-reductions and pinched Riemannian 
hyperbolic metrics; however, 
we use the incomplete hyperbolic metric itself.
The hard geometric analysis and geometric convergence techniques can be avoid 
using the techniques of this paper. 
By a negatively curved space, we mean a metric space whose 
universal cover is \cat$(-1)$.
Except for developing a somewhat complicated theory of 
deforming boundary to make the submanifolds of codimension $0$ 
negatively curved, we do not need any other highly developed 
techniques.
Also, we might be able to generalize the techniques to
negatively curved polyhedral $3$-manifolds and complexes obtained from groups.

Note also that there is a recent paper by Calegari and Gabai
\cite{CalGab} using modified least area surfaces and closed geodesics.
The work here is independently developed from their line of 
ideas. Also, we note that there were earlier attempts by 
Freedman \cite{Freedman1}, Freedman-McMullen \cite{FM}, 
which were very influential for the later success by Agol 
and Calegari-Gabai, 
and another earlier unsuccessful attempt by Ohshika, using 
the least area surfaces. 
%(The author did formulated a strategy using 
%the boundary surfaces of convex hulls similar to what is in this paper 
%to prove the tameness conjecture
%as early as May of 2003 in a topology conference in Arkansas and 
%a note of which is to be posted in my site).

%** main result
In this paper, we let $M$ be a hyperbolic $3$-manifold 
with a Scott's core homeomorphic to a compression body. 
We suppose that $M$ has a finitely generated fundamental group
and the holonomy is purely loxodromic
and has ends $E, E_1, \dots, E_n$. 
Let $F_1, \dots, F_n$ be the incompressible surfaces in 
neighborhoods of the ends $E_1, \dots, E_n$ respectively.  
Let $N(E)$ be a neighborhood of an end $E$ with no incompressible surface 
associated. 

The Marden's conjecture states that a hyperbolic $3$-manifold 
with a finitely generated fundamental group is homeomorphic to the interior of 
a compact $3$-manifold. It will be sufficient to prove for 
the above $M$ to prove the complete conjecture. 

The cases when the group contains parabolic elements are left out,
which we will work out on a later occasion. 

\begin{thmA}
Let $M$ be as above with ends $E, E_1, \dots, E_n$,  
and $\core$ be a compact core of $M$. Then 
$E$ has an exiting sequence of surfaces of genus equal to 
that of the boundary component of $\partial \core$ corresponding to $E$.  
\end{thmA} 

\begin{thmB}
$M$ is tame{\rm ;} that is, 
$M$ is homeomorphic to the interior of a compact manifold. 
\end{thmB}

This paper has three parts: 
In Part 1, 
let $M$ be a codimension $0$ submanifold of a hyperbolic 
$3$-manifold $N$ of infinite volume with certain nice 
boundary conditions. $M$ is locally finitely 
triangulated. Suppose that $M$ 
is {\em $2$-convex} in $N$ in the sense that any tetrahedron $T$ in $N$ with 
three of its side in $M$ must be inside $M$. 
Now let $L$ be a finitely triangulated, compact codimension $0$-submanifold $M$ 
so that $\partial L$ is incompressible in $M$ with a number 
of closed geodesics $\cl_1, \dots, \cl_n$ removed.
Given $\eps > 0$, 
we show that $\partial L$ can be isotopied to a hyperbolically
triangulated surface so that it bounds in $M$
a $2$-convex submanifold whose $\eps$-neighborhood contains 
$\cl_1, \dots, \cl_n$. 
The isotopy techniques will be based on PL-type arguments and deforming by 
crescents. An important point to be used in the proof is 
that the crescents avoid closed geodesics and geodesic laminations. 
Thus, the isotopy does not pass through the closed geodesics
by as small amount as we wish.  (These will be the contents of 
Theorem C and Corollary D.) 

Part 2 is as follows: 
A {\em general hyperbolic manifold} is a manifold with boundary modeled 
on subdomains in the hyperbolic space. 
A general hyperbolic manifold 
is {\em $2$-convex} if every isometry from a tetrahedron with 
an interior of one of its side removed extends to the tetrahedron itself. 
We show that a $2$-convex general hyperbolic manifold is 
negatively curved. The proof is based on 
the analysis of the geometry of the vertices of the boundary 
required by the $2$-convexity.
We will also define a hyperbolic surface as 
a triangulated surface where each triangle gets mapped to 
a geodesic triangle and the sum of the induced angles at each 
vertex is always greater than or equal to $2\pi$. We 
show the area bound of such surfaces. Finally, we 
show that the boundary of the convex hull of a core in 
a general hyperbolic manifold with finitely generated fundamental
group can be deformed to a nearby hyperbolic-surface, which
follows from the local analysis of geometry.

In Part 3, we will give the proof of Theorems A and B
using the results of Part 1 and 2. The outline is 
given in the abstract and in the beginning of Part 3. The proof 
itself is rather short spanning 9-10 pages only.

We like to thank the valuable comments and encouragements from 
Professors Francis Bonahon, Richard Canary, Bill Goldman, Misha Kapovich, Sadayoshi Kojima,
Kenichi Ohshika, Curt McMullen, Yair Minsky and Bill Thurston.

\part{$2$-convex hulls of submanifolds of 
hyperbolic manifolds}

\section{Introduction to Part 1}

%** 0.2 purpose

The purpose of Part I is to deform a submanifold of 
codimension $0$ of general hyperbolic manifold into 
a negatively curved one, i.e., Corollary D. 

%** 0.3 definitions

%Loosely define
%special hyperbolic manifolds

%general hyperbolic manifolds 

We will be working in a more general setting.
A {\em general hyperbolic manifold} is a Riemannian manifold $M$ with 
corner and a geodesic metric that admits
a geodesic triangulation so that each
$3$-simplex is isometric with a compact hyperbolic 
simplex. Let $\tilde M$ be a universal cover of $M$ and 
$\pi_1(M)$ the group of deck transformations. 
$M$ admits a local isometry, so-called developing map, 
$\dev: \tilde M \ra \hyps$ for the hyperbolic space 
$\hyps$ equivariant with respect to
a homomorphism $h: \pi_1(M) \ra \PSL(2, \bC)$. 
The pair $(\dev, h)$ is only determined up to action 
\[(\dev, h) \mapsto (g\circ \dev, g\circ h(\dot) \circ g^{-1})
\hbox{ for }  
g \in\PSL(2, \bC). \]

We remark that in Thurston's notes \cite{Thnote} 
a locally convex general hyperbolic manifold 
is shown to be covered by a convex domain in $\hyps$ or, 
equivalently, it can be 
extended to a complete hyperbolic manifold. 
However, a general hyperbolic manifold can be much worse
although here we would be looking at mostly general hyperbolic 
manifolds that are covered by 
coverings of some domains in the hyperbolic 
spaces.

%%% October 5, 2004

%Furthermore, $M$ is special if a regular cover of $M$ 
%is a domain in $\hyps$. In this case, $M$ is isometric 
%with the domain quotient out by a Kleinian group. 

A general hyperbolic manifold $M$ is {\em $2$-convex} if 
given a hyperbolic $3$-simplex $T$, 
a local-isometry $f: T - F^o \ra M^o$ for a face $F$ of $T$
extends to an isometry $T \ra M^o$. 
(See \cite{psconv} for more details. Actually 
projective version applies here by the Klein 
model of the hyperbolic $3$-space.)

\begin{itemize}
\item By a {\em totally geodesic hypersurface}, we mean the union of 
components of the inverse image under a developing
map of a totally geodesic plane in $\hyps$. 
\item A {\em local totally geodesic hypersurface} is 
an open neighborhood of a point in the hypersurface. 
\item For a point, a {\em local half-space} is the closure 
in the ball of the component of an open ball around it with 
a totally geodesic hypersurface passing through it. 
\item The local totally geodesic hypersurface intersected with
the local-half space is said to be the {\em side} of the local half-space. 
\item A local half-space with its side removed is said to be an {\em open} local
half-space.
\end{itemize}

%saddle vertices hyperbolic-vertices

A surface $f: S \ra M$ is said to be {\em triangulated} 
if $S$ is triangulated and each triangle is mapped 
to a totally geodesic triangle in $M$. 
(We will generalize this notion a bit.)
\begin{itemize}
\item An interior vertex of $f$ is a {\em saddle-vertex} if 
every open local half-space associated with 
the vertex does not contain the local image 
of $f$ with the vertex removed. 
\item A {\em strict saddle-vertex} is a saddle-vertex where 
every associated closed local half-space does not 
contain the local image.
\item An interior vertex of $f$ is a {\em convex vertex} if 
a local open half-space associated with
the vertex contains the local image of 
$f$ removed with the vertex. 
\end{itemize}
(In these definitions, we take a slightly larger ambient open manifold of $M$ 
to make sense of open local half-space. Also, 
we can see these definitions 
better by looking at the unit 
tangent bundle $U$ at the vertex:
Then the image of $f$ corresponds to a path 
in $U$.)
If $P$ is a local totally geodesic 
plane passing through the vertex whose 
one-sided neighborhood contains the local image of $f$, 
then we say $P$ is a {\em supporting plane}. 

An interior vertex of $f$ is an {\em hyperbolic-vertex} if 
the sum of angles of the triangles around a vertex 
is $\geq 2\pi$. 
A saddle-vertex is a hyperbolic-vertex by Lemma \ref{lem:lengthhemisphere}.
%lem:lengthhemisphere.
%\cite{genhypmfld}.

% map case
A map $f:S \ra M$ is a saddle-map if each vertex is a saddle-vertex,
and $f$ is a hyperbolic-map if each vertex is a hyperbolic-one.
A saddle-map is a hyperbolic-map but not conversely in general.
For an imbedding $f$ and an orientation, a convex vertex is 
said to be a {\em concave vertex} if the local half-space is 
in the exterior direction. Otherwise, the convex vertex is 
a {\em convex vertex}.
We also know that if the boundary of 
a general hyperbolic submanifold $N$ of $M$ is saddle-imbedded,
then $N$ is $2$-convex (see Proposition \ref{prop:s-bd2-conv}).

%our tools: crescents

%** 0.4 main theorem
\begin{thmC} 
Let $M$ be an orientable $2$-convex general hyperbolic 3-manifold, 
and $\Sigma$ be a closed surface in $M$.
Suppose that each component of 
$\Sigma$ is incompressible in $M$ if 
we remove a finite number of the image $\cl_1, \dots, \cl_n$ of 
the closed geodesics in the interior $M^o$ of $M$. 
Then for arbitrarily given small $\eps > 0$,
we can isotopy $\Sigma$ so that $\Sigma$ becomes a saddle-imbedded surface.
Finally, during the isotopy $\Sigma$ might pass through some of $\cl_i$ 
but with as small an amount as possible.
\end{thmC}

%% Main Corollary
\begin{corD} 
Let $N$ be an orientable $2$-convex general hyperbolic 3-manifold, 
and $M$ be a compact codimenion-zero submanifold of $N$ with boundary
$\partial M$. Suppose that each component of 
$\partial M$ is incompressible in $M$ if 
we remove a finite number of the image $\cl_1, \dots, \cl_n$ of 
the closed geodesics in the interior $M^o$ of $M$. 
Then for arbitrarily given small $\eps > 0$,
we can isotopy $M$ to a homeomorphic general hyperbolic 
$3$-manifold $M'$ in $N$ so that $M'$ is 2-convex
and an $\eps$-neighborhood of 
$M'$ contains the collection of closed geodesics $\cl_1, \dots, \cl_n$.
\end{corD}

We say that $M'$ obtained from $M$ by the above 
process is a {\em $2$-convex hull} of $M$. 
Although $M'$ is not necessarily a subset of $M$,
the curves $\cl_1, \dots, \cl_n$ is a subset of an $\eps$-neighborhood 
of $M'$ and hence we have certain amount of control. 
(Here the geodesics are allowed to self-intersect.)

%%** 0.5 outline

%main strategy 

%mention crescents.

%outline here
In section \ref{sec:pre}, we review some hyperbolic manifold theory
and discuss saddle vertices and relationship with 
$2$-convexity.

In section \ref{sec:cres}, we introduce so-called crescents:
Let $M$ be a $2$-convex general hyperbolic manifold and 
$\Sigma$ a closed subsurface, possibly with many components. 
We take the inverse image $\tilde \Sigma$ in $\tilde M$ of 
$\Sigma$ in $M$, incompressible 
in the ambient $2$-convex general hyperbolic 
manifold with a number of geodesics 
$\cl_1, \dots, \cl_n$ removed.
A crescent is a connected domain bounded by a totally geodesic hypersurface 
and an open surface in $\tilde \Sigma$. 
The portion of boundary in the totally geodesic hypersurface 
is said to be the $I$-part and the portion in 
$\tilde \Sigma$ is said to be the $\alpha$-part. 
A crescent may contain another crescents and so on, and 
the folding number of 
a crescent is the maximum intersection number of 
the generic path from the outer part in the surface 
to the innermost component of the crescent with 
the surface removed. We show that for a given closed surface $\Sigma$,
the folding number is bounded above.

A {\em highest-level crescent} is an innermost one that is contained 
in a crescent with highest folding number which achieves the folding number. 
We show that a highest-level crescent is always contained 
in an innermost crescent; i.e., so called the secondary 
highest-level crescent. In a secondary highest level crescent, 
the closure of the $\alpha$-part and the $I$-part are isotopic.
We also show that the secondary highest-level crescents meet nicely
extending their $\alpha$-parts in $\tilde \Sigma$, following 
\cite{psconv}.

In section \ref{sec:crescent-iso}, we introduce the crescent-isotopy theory
to isotopy a surface in a general hyperbolic 
manifold so that all of its vertices become saddle-vertices:
We form the union of secondary highest-level crescents
and can isotopy the union of 
their $\alpha$-parts to the complement $I$ in 
the boundary of their union. This is essentially 
the crescent move. 
(In this paper, by isotopies, we mean the deck-transformation group equivariant isotopies unless 
we specify otherwise. At least, if the isotopy in each step is not equivariant, we will make it so 
after the final isotopy is completed.)

However, there might be some parts of $\tilde \Sigma$ 
meeting $I$ tangentially from above. We need 
to first push these parts upward first
using so-called convex truncations.

Also, after the move, there might be pleated parts 
which are not triangulated. We present 
a method to perturb these parts to triangulated parts
without increasing the levels or the set of crescents by much. 

Next, we use the crescent isotopy and perturbations  
to obtain saddle-imbedded surface isotopic to $\Sigma$:
\begin{itemize}
\item We take the highest folding number and take all 
outer secondary highest-level crescents, 
\begin{itemize}
\item do some convex truncations,  
\item do the crescent isotopies and 
\item convex perturbations.
\end{itemize}
\item Next, we take all inner highest-level crescents of the same level
as above, do some truncations,  
and do crescent moves and perturb as we did above. 
Now the highest folding number decreases by one.
\item We do the next step of the induction
until we have no crescents any more. 
\end{itemize}
In this case, all the vertices are saddle-vertices. 
This completes the proof of Theorem C.

Finally, we prove Corollary D by
applying our results to a codimension-zero submanifold $M$ with
incompressible boundary in the ambient manifold with some
geodesics removed.

% Dec. 15 3:45
%*1. preliminary notation settings 
\section{Preliminary} \label{sec:pre}

%** 1 intro and outline of the section

In this section, we review the hyperbolic space 
and the Kleinian groups briefly. We 
discuss the relationship between the $2$-convexity of 
general hyperbolic manifolds.

%** 1.1. hyperbolic space and Kleinian groups 
\subsection{Hyperbolic manifolds}
%%% NN \hyp^3, \Isom, U, U^2, x, B^2, B, \dist 

The hyperbolic $n$-space is a complete Riemannian metric 
space $(\bH^n, d)$ of constant curvature equal to $-1$. 
We will be concerned about hyperbolic plane and 
hyperbolic spaces, i.e., $n=2,3$, in this paper. 

The upper half space model for $\hypp$ is the pair 
\[(U^2, \PSL(2, \bR) \cup \overline{\PSL}(2, \bR))\]
where $U^2$ is the upper half space.

The Klein model of $\hypp$ is the pair $(B^2, \PO(1, 2))$
where $B^2$ is the unit disk and $\PO(1,2)$ is 
the group of projective transformations acting on $B^2$. 

A {\em Fuchsian} group is a discrete subgroup of 
the group of isometries of $\PO(1, 2)$ of the group 
$\Isom(\hypp)$ of isometries of $\hypp$.

There are many models of the hyperbolic $3$-space:
The upper half-space model consists of 
the upper half-space $U$ of $\bR^3$ and 
the group of isometries are identified as 
the group of similarities of $\bR^3$ preserving $U$, 
which is identified as the union of $\PSL(2, \bC)$ and 
its conjugate $\overline{\PSL}(2, \bC)$. 

We shall use the Klein model mostly:
The Klein model consists of the unit ball 
in $\bR^3$ and the group of isometries are 
identified with the group of projective transformations 
preserving the unit ball, which is identified 
as $\PO(1, 3)$. 

A {\em Kleinian} group is a discrete subgroup of the group of 
isometries $\PO(1, 3)$, i.e., the group $\Isom(\hyps)$
of isometries of $\hyps$. 

A {\em parabolic} element $\gamma$ of a Kleinian group is 
a nonidentity element such that $\dist(\gamma(x), x)$, $x \in \hyps$, has 
no lower bound other than $0$.
A {\em loxodromic} element of a Kleinian group is 
an isometry with a unique invariant axis. 
A {\em hyperbolic} element is a loxodromic 
one with invariant hyperplanes. 

For Fuchsian groups, a similar terminology holds. 

In this paper, we will restrict our Kleinian groups 
to be torsion-free and have no parabolic elements
and all elements are orientation-preserving.

\subsection{Saddle-vertices}
Let $M$ be a general hyperbolic manifold. 
We now classify the vertices of a triangulated map $f:S \ra M$
where we do not yet require the general position property of $f$ 
but identify the vertex with its image. 

By a {\em straight geodesic} in a general hyperbolic manifold, we 
mean a geodesic that maps to geodesics in $\hyps$
under the developing maps.

%%% October 7th 11:24
\begin{lem}\label{lem:vertex}
Let $f:S \ra M$ be a triangulated map. 
\begin{itemize}
\item An interior vertex of $S$ 
is either a convex-vertex, a concave vertex, or a saddle-vertex.
\item A saddle-vertex which is not strict one 
has to be one of the following: 
\begin{itemize}
\item[(i)] a vertex with a totally geodesic local image.
\item[(ii)] a vertex on an edge in the intersection of 
two totally geodesic planes 
where $f$ locally maps into one sides of each plane.
\item[(iii)] A vertex which is contained in at least three edges in the image of $f$ 
in a local totally geodesic plane $P$ and the edges are not contained in 
any closed half-plane of $P$. The local half-space bounded by $P$ 
is the unique one containing the image of $f$. 
\end{itemize}
\item If $f$ is a general position map, then 
a saddle-vertex is a strict saddle-vertex.
\end{itemize}
\end{lem}
\begin{proof} 
Suppose $v$ is a saddle-vertex and not a strict one 
and not of form (i) or (ii). 
Since $v$ is not strict, there is a supporting plane. 
If there are more than three supporting planes in general position,
then $v$ is a strict convex-vertex. If there are two supporting plane, 
then an edge in the image of $f$ is in the edge of 
intersection of the two planes in order that $v$ be a nonstrict 
saddle vertex, which is absurd. 

Hence the supporting plane is unique. If there are no three edges 
as described in (iii), we can easily find another supporting plane.
\end{proof}

%%March 30, 11:55

\begin{lem}\label{lem:deforming} 
A saddle-vertex of type {\rm (ii)} and {\rm (iii)} of Lemma \ref{lem:vertex}
can be deformed to a strict saddle-vertex by 
an arbitrarily small amount by pushing if necessary the vertex  
from the boundary of the closed local half-space
containing the local image of $f$ in the direction of 
the open half-space. 
\end{lem}
\begin{proof} 
If the saddle-vertex is a strict one, then we leave it alone.
If the saddle-vertex is not a strict one, a closed local half-space
contains the local image of $f$.
Let $U$ be the unit tangent bundle at the vertex.
A closed hemisphere $H$ contains the path corresponding to 
the local image of $f$. 

%In case (i) of Lemma \ref{lem:vertex}, we simply push our vertex $v$ in 
%the normal direction to the supporting plane. 
In case (ii) of Lemma \ref{lem:vertex}, there are actually two closed hemispheres 
$H_1$ and $H_2$ whose intersection contains the local
image of $f$ in $U$. Therefore, we choose a direction 
in the interior of the intersection of $H_1$ and $H_2$. 
Then by Lemma \ref{lem:cross}, the result of a sufficiently 
small deformation is a strict saddle vertex. 

In case (iii) of Lemma \ref{lem:vertex}, let $w_1, w_2, w_3$ be the points on 
the unit tangent bundle at the vertex corresponding the three edges. 
Then by moving vertex in the direction, 
the corresponding directions $w'_1, w'_2, w'_3$ of the perturbed 
edges form a strictly convex triangle in an open hemisphere in $U$
and the direction vector $v$ of the movement not in the hemisphere.  
$v, w'_1, w'_2, w'_3$ are vertices of a geodesic triangulation of $U$ into 
a $2$-skeleton of the topological tetrahedron and every triple of them form 
vertices of a strictly convex triangle in an open hemisphere.
Therefore, there is no closed hemisphere in $U$ containing all of them.
Hence, we obtain a strict saddle vertex.
\end{proof}

We will need the following much later:
\begin{lem}\label{lem:cross} 
Suppose that $f:\Sigma \ra M$ is a triangulated imbedding. 
Let $v$ be a vertex of $f$ and $f':U_v \ra U^M_v$ be 
the induced map from the link of $v$ in $\Sigma$
to that of $v$ in $M$. Suppose that there exists 
a segment $l$ of length $> \pi$ in $U^M_v$ with 
endpoints in the image of $f'$ separating 
two points in the image of $f'$ so that the minor 
arc $\ovl{xy}$ meets $l$ transversely.
Then $v$ is a strict saddle vertex.  
\end{lem}
\begin{proof} 
If a closed hemisphere contains $l$, then it must contain $l$ in its boundary. 
Therefore, $x$ or $y$ is not in the hemisphere, and there is no closed hemisphere 
containing $l$ and $x$ and $y$. 
\end{proof}

Given an oriented surface, a convex vertex is 
either a {\em convex vertex} or a {\em concave vertex} 
depending on whether the supporting local half-space 
is in the outer normal direction or in the inner normal 
direction.

\begin{prop}\label{prop:alts-vertex} 
A vertex of an oriented imbedded triangulated surface 
is either a saddle-vertex or a convex vertex or 
a concave vertex.
\end{prop}
\begin{proof} 
Straightforward.
\end{proof}

%*** 1.2.3. saddle-boundary iff 2-convexity
In this paper, we consider only metrically 
complete submanifolds, i.e., locally compact ones. 

\begin{prop}\label{prop:s-bd2-conv} 
A general hyperbolic manifold
$M$ is $2$-convex if and only if each 
vertex of $\partial M$ is a convex vertex or a saddle-vertex.
\end{prop}
\begin{proof}
Suppose $M$ is $2$-convex.
If a vertex $x$ of $\partial M$ is a concave, 
we can find a local half-open space in $M$ with its side 
passing through $x$. The side meets $\partial M$ 
only at $x$. From this, we can find a $3$-simplex 
inside with a face in the side. This contradicts 
$2$-convexity of $M$. 

Conversely, suppose that $\partial M$ has only 
convex vertices or saddle-vertices.
Let $f:T - F^o\ra M^o$ be a local-isometry from 
a $3$-simplex $T$ and a face $F$ of $T$. 
We may lift this map 
to $\tilde f: T - F^o \ra \tilde M^o$ where 
$\tilde M$ is the universal cover of $M$ where $\tilde f$ is an imbedding. 

Since $\tilde M$ is metrically complete, 
$\tilde f$ extends to $\tilde f': T \ra \tilde M$. 
Suppose that $f$ does not extend to $f':T \ra M^o$. 
This implies that $\tilde f'(F)$ meets $\partial \tilde M$ where
$\tilde f'(\partial F)$ does not meet $\partial \tilde M$. 
The subset $K = \partial \tilde M \cap \tilde f'(F)$ 
has a vertex $x$ of $\partial \tilde M$ 
which is an extreme point of the convex 
hull of $K$ in the image of $F$. 
We can tilt $\tilde f'(T)$ by a supporting line $l$ at $x$ a bit and 
the new $3$-simplex meets $\partial \tilde M$ at $x$ 
only. This implies that $x$ is not a saddle-vertex but 
a concave vertex, a contradiction.
\end{proof}

%* 2. 2-convex hulls of hyperbolic manifolds
%N D, \Gamma 

\section{Crescents}\label{sec:cres}

%% Outline 

Let $M$ be a metrically complete $2$-convex 
general hyperbolic manifold from now on and $\tilde M$ its universal cover. 
Let $\Gamma$ denote the deck transformation group of 
$\tilde M \ra M$. 

Let $\Sigma$ be a properly imbedded compact subsurface of 
an orientable general hyperbolic manifold $M$ with more than one components in general.
We denote by $\tilde \Sigma$ the inverse image of 
$\Sigma$ in the universal cover $\tilde M$ of $M$. 
($\tilde \Sigma$ is not connected in general and 
components may not be universal covers of $\Sigma$.)
We assume that the triangulated $\tilde M$ is in general position
and so is $\tilde \Sigma$ under the developing maps. 

For each component $\Sigma_0$ of $\Sigma$ and 
a component $\Sigma_0'$ of $\tilde \Sigma$
mapping to $\Sigma_0$, there exists a subgroup 
$\Gamma_{\Sigma_0'}$ acting on $\Sigma_0'$
so that the quotient space is isometric to $\Sigma_0$.

\begin{hyp}\label{hyp:incompressible}
We will now assume that $\Sigma$ is incompressible in 
$M$ with a number of straight closed geodesics $\cl_1, \dots, \cl_n$ removed. 
\end{hyp}

First, we introduce crescents 
for $\tilde \Sigma$ which is the inverse image 
of a surface $\Sigma$ in a $2$-convex general hyperbolic manifold.
We define the folding number of crescents and 
show that they are bounded above. 

We define the highest level crescents, i.e., the innermost crescents in the 
crescent with the highest folding number
incurring the highest folding number.
We show that closed geodesics avoid the interior of crescents. 
Given a highest level crescent, we show that
there is an innermost crescent that has a connected 
$I$-part to which the closure of the $\alpha$-part is isotopic
in the crescent by the incompressibility of $\Sigma$.
These are the {\em secondary highest-level crescents}.
We show that the secondary highest-level crescent is 
homeomorphic to its $I$-part times the unit interval.

Next, we show that if two highest-level crescents meet each other 
in their $I$-parts tangentially, then they both are included 
in a bigger secondary highest-level crescent.

Furthermore, if two secondary highest-level crescents meet
in their interiors, then they meet nicely extending their 
$\alpha$-parts. This is the so-called transversal intersection 
of two crescents. 

%** 2.1. Crescents and crescent-moves
%*** 2.1.1. Definition of crescents

\subsection{Definition of crescents}

%%% April 13 14:21
We list here key notions associated with 
crescents that are used very often in this paper. 
\begin{defn}\label{defn:triangulations}
Clearly, we require $\tilde \Sigma$ to be a properly and tamely 
imbedded subsurface disjoint from $\tilde M$. 
\begin{itemize}
\item If $\tilde \Sigma$ is not necessarily triangulated but 
each point of it has a convex $3$-ball neighborhood $B$
where the closure of one of the component of $B - \tilde \Sigma$ 
is the closure of a component of a convex $3$-ball with 
a closed triangulated disk with boundary in $\partial B$ removed.
In this case, $\tilde \Sigma$ is said to be {\em nicely} imbedded.
\item If $\tilde \Sigma$ is a union of triangulated compact triangles
but the vertices are not necessarily in general position, 
we say that $\tilde \Sigma$ is {\em triangulated}. 
\item If $\tilde \Sigma$ is triangulated by compact 
triangles whose vertices are in general position, 
we say that $\tilde \Sigma$ is {\em well-triangulated}. 
\end{itemize}
\end{defn}

\begin{defn}\label{defn:crescent} 
We assume that $\tilde \Sigma$ is nicely imbedded at least. 
A {\em crescent} $\mathcal R$ for $\tilde \Sigma$ is 
\begin{itemize}
\item a connected domain in $\tilde M$ which is a closure of 
a connected open domain in $\tilde M$, 
\item so that its boundary is a disjoint union of 
a (connected) open domain in $\tilde \Sigma$ and 
the closed subset that is 
the disjoint union of totally geodesic $2$-dimensional domains 
in $\tilde M$ that develops into 
a {\bf common} totally geodesic hypersurface in $\hyps$ under $\dev$. 
\end{itemize}
We denote by $\alpha_{\mathcal R}$ the domain in $\tilde \Sigma$ 
and $I_{\mathcal R}$ the union of totally geodesic domains. 
To make the definition canonical, we require 
$I_{\mathcal R}$ to be a maximal totally geodesic set in 
the boundary of $\mathcal R$.
We say that $I_\calR$ and $\alpha_\calR$ the {\em $I$-part} 
and the {\em $\alpha$-part} of $\calR$. 
\end{defn}

%%06/14 11:45

As usual $\Sigma$ is oriented so that there are outer and 
inner directions to normal vectors. 

\begin{defn}\label{defn:foldingnumber} 
The subset $\tilde \Sigma \cap \mathcal R$ may have more than 
one components. For each component of ${\mathcal R} - \tilde \Sigma$,
we can assign a {\em folding number} which is the minimal generic 
intersection number of that a path in the interior of $\calR$
from $\alpha_{\mathcal R}$ meeting $\tilde \Sigma$ to reach to the component. 
The folding number of $\calR$ is the maximum of the folding numbers 
for all of the components. 
\end{defn}

\begin{defn}
A {\em sub-crescent} $\calS$ of a crescent $\calR$ is 
the closure of a component of $\calR - A$ where $A$ is a union 
of components of $\tilde \Sigma \cap \calR$ where 
we define the $I$-part to be the union of maximal domains 
in the intersection of the $I_\calR$ with $\calS$ 
and $\alpha_\calS$ to be $\partial \calS - I_\calS$. 
In these case $\calR$ is a {\em super-crescent} of $\calS$.

Given a crescent $\calS$, an {\em ambient folding number} is 
the maximum of the folding number of super-crescents of $\calS$. 
\end{defn} 

A priori, a crescent may have an infinite folding number. 
However, we will soon show that the folding number is finite. 

Note that a proper sub-crescent has a strictly less folding number
than the original crescent and a strictly greater ambient folding number
than the original one. 

\begin{defn}
The {\em $I$-part hypersurface} is the inverse image in $\tilde M$ 
under the developing map of the totally geodesic plane $P$ 
containing the developing image of the $I$-part of the crescent.
\end{defn}

A closed subset $K$ of $\tilde M$ is a {\em geometric limit} or just {\em limit} of 
a sequence of closed subsets $K_i$ if for each compact ball $B$ in $\tilde M$, 
$K_i \cap B$ converges to $K \cap B$ in the Hausdorff metric sense. 

\begin{defn}
A noncompact domain will be called a {\em crescent} if 
it is bounded by a (connected) domain in $\tilde \Sigma$ and 
the union of totally geodesic domains 
developing into a common totally geodesic plane 
in $\hyps$ and is a geometric limit of 
compact crescents. Again the $I$-part is the maximal totally geodesic 
subset of the boundary of the crescent.
The $\alpha$-part is the complement in the boundary 
of the crescent and is a connected open subset of $\tilde \Sigma$.
Of course they need not be limits of the corresponding 
subsets of the compact crescents. 
(See Proposition \ref{prop:closedness} for 
a related idea.)
\end{defn}

\begin{defn}
A crescent is an {\em outer} one if its interior to 
the $\alpha$-part is in the outer normal direction 
of $\tilde \Sigma$. It is an {\em inner} one otherwise. 
\end{defn}

\begin{defn}
The boundary $\partial I_{\mathcal R}$ of $I_{\mathcal R}$
is the set of boundary points in the $I$-part hypersurface of 
$\calR$. Also, for any subset $A$ of $I_\calR$, 
we define $\partial A$ to be the set of boundary points 
in the $I$-part hypersurface. We define $I^o_\calR$ to be the interior, 
i.e., $I_\calR - \partial I_\calR$. 
\end{defn}

\begin{defn}
A {\em pinched simple closed curve} is a simple curve
pinched at most three points or pinched at a connected arc. 
The boundary of the $I$-part is 
a disjoint union of pinched simple curves. 
\end{defn}

\subsection{Properties of crescents}
%% Oct. 15, 2004

%**** 2.1.1.1. 
The following is a really important property
since this shows we can use crescents in general hyperbolic manifolds 
without worrying about whether the $I$-parts meet the boundary of
the ambient manifold. 

\begin{prop}\label{prop:disjcrescent} 
Let $\calR$ be a crescent in a $2$-convex ambient general hyperbolic 
manifold $M$. Then $\calR$ is disjoint from $\partial \tilde M$.
In fact, if $\calR$ is compact, then $\calR$ is uniformly bounded 
away from $\partial \tilde M$.
\end{prop}
\begin{proof}
Suppose that $\calR$ meets $\partial \tilde M$. 
Since the closure of $\alpha_\calR$ being a subset of $\tilde \Sigma$ 
is disjoint from $\partial \tilde M$, 
it follows that $I_{\calR}$ meets $\partial \tilde M$ 
in its interior points and away from the boundary points 
in the ambient totally geodesic subsurface $P$ in $\tilde M$. 

We find the extreme point of $I_{\calR} \cap \partial M$ 
and find the supporting line. This point has a local half-space 
in $\calR$. By tilting the $I$-part a bit by the supporting line, 
we find a local half-space in $M$ and in it a local totally geodesic 
hypersurface meeting $\partial M$ at a point. This contradicts 
the $2$-convexity of $\tilde M$.

The second part follows from the disjointness of $\calR$ to 
$\partial M$ and the compactness of $\calR$. 
\end{proof}

\begin{figure}[ht]
\centerline{\epsfxsize=2.5in \epsfbox{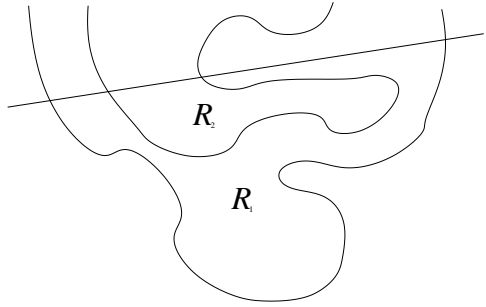}}
\caption{\label{fig:crescent} 
A $2$-dimensional section of crescents 
$\calR_1$ and $\calR_2$ where $\calR_2$ is $1$-nested 
and is innermost.}

\end{figure}
\typeout{<>}

%**** 2.1.2.2 A Limit of crescents is a subset of a crescent

The following shows the closedness of set of points of 
$\tilde M - \tilde \Sigma$ in crescents. 

\begin{prop}\label{prop:closedness} 
Let ${\calR}_i$ be a sequence of crescents. 
Suppose that $x$ is a point of $\tilde M - \tilde \Sigma$ 
which is a limit of a sequence of points in the union of 
${\calR}_i$. Then $x$ is contained in a crescent. 
\end{prop} 
\begin{proof} 
Let $x_i \in \bigcup_{j} {\calR}_j$ be a sequence 
converging to $x$. We may assume that $x$ is not an 
element of any ${\calR}_j$. 

Using geometric convergence, there exists a totally geodesic hypersurface
$P$ through $x$ and a geometric limit of a sequence of 
$I_{\calR_j}$ converging to a subset $D$ in $P$. 

Then $D$ is separating in $\tilde M - \tilde \Sigma$.
If not, there exists a simple closed curve $\gamma$ in $\tilde M$
meeting $D$ only once, which means $I_{\calR_j}$ for a sufficiently 
large $j$ meets $\gamma$ only once as well.

Now, $\tilde M - D$ may have more than one components. 
Since $x \in D$, we take an open domain $L$ bounded by $D$ and a subset of $\tilde \Sigma$
and whose closure contains $x$. 
The closure of $L$ is a crescent containing $x$ since $R_j$ is a geometric limit of a sequence of 
compact crescents $R_{j, i}$, $i=1,2,3, \dots$ and we can choose a subsequence from these 
converging to $L$. 
\end{proof}

\begin{defn}\label{defn:limitcr} 
Let $\Sigma_i$ be a sequence of nicely imbedded surfaces isotopic to $\Sigma$.
By the method of above Proposition \ref{prop:closedness}, for any 
sequence of crescents $R_i$ for surface $\Sigma_i$ with a geometric limit $K$ with 
nonempty interior, it follows that 
there is a maximal crescent $R$ contained in $K$ with $I_R$ in the geometric limit of 
a subsequence of totally geodesic hypersurfaces $P_i$ containing $I_{R_i}$.  
$R$ is said to be {\em generalized limit} of the sequence $R_i$.
\end{defn}

The definition is needed since we might have some parts of domain degenerating 
to lower-dimensional objects. 
The generalized limit might not be unique. In fact, there could be 
two with disjoint interiors. 

\begin{prop}\label{prop:alphalimit} 
Suppose that a sequence of crescents $R_i$ converges to a crescent $R$
in the generalized sense. 
Then $\alpha_R$ is a subset of a limit of any subsequence of $\clo(\alpha_{R_i})$. 
\end{prop}
\begin{proof} 
The union of a totally geodesic $2$-dimensional domain and $\clo(\alpha_{R_i})$ 
is the boundary of $R_i$. The limit of a subsequence of 
the boundary of $R_i$ converges to a subset containing the boundary of $R$.
\end{proof}

%**** 2.1.2.3. Upper bound to size and folding number

As usual, 
we assume that the holonomy group of $\Sigma$ 
does not consist strictly of parabolic or elliptic or identity 
elements. 

A {\em size} of a crescent is the supremum of 
the distances $d(x, \alpha_{{\mathcal R}})$ for
$x \in I_{{\mathcal R}}$. We show that this is globally 
bounded by a constant depending only on $\Sigma$. 

First, a complementary result is proved:
\begin{lem}\label{lem:upper1} 
Every $x \in \alpha_{\calR}$ satisfies 
\[ d(x, I_{\calR}) \leq N \] 
for a uniform constant $N$ depending on $M$ and $\Sigma$ only. 
\end{lem}
\begin{proof} 
If not, since there is a compact fundamental domain 
in $\tilde \Sigma$, 
using deck transformations acting on $\tilde \Sigma$,
we obtain a sequence of bigger and bigger 
compact crescents where the corresponding sequence of 
the $I$-parts leave any compact subset of $\tilde M$
and the corresponding sequence of $\alpha$-parts meets 
a fixed compact subset of $\tilde M$. 
Therefore, we form a subsequence of
the developing images of the $I$-parts converging to a point
of the sphere at infinity of $\hyps$. 

Let ${\mathcal R}_i$ be the corresponding crescents.
Then $\alpha_{{\mathcal R}_i}$ is a subsurface with boundary 
in the $I$-parts, and 
$\alpha_{{\mathcal R}_i}$ contains any compact 
subset of $\tilde \Sigma$ eventually.
 
Let $c$ be a closed curve in $\Sigma$ with nonidentity holonomy.
Let $\tilde c$ be a component of its inverse image in $\tilde \Sigma$.
Since $\tilde c$ must escape any compact subset of 
$\tilde \Sigma$, $\tilde c$ escape $\alpha_{{\mathcal R}_i}$. 
Thus, $\tilde c$ must meet all $I_{{\mathcal R}_i}$ for $i$ sufficiently large.
Since the developing image of $\tilde c$ has two well-defined endpoints, 
this means that the limit of the sequence of $I$-parts must 
contain at least two points, a contradiction. 
\end{proof}

\begin{prop}\label{prop:cresupper}
Let $M$ and $\Sigma$ be as above. 
Then $d(x, \alpha_{\calR})$ for $x \in I_\calR$ is uniformly bounded above by 
a constant depending only on $M$ and $\Sigma$ and, hence,
there is an upper bound to the size of a crescent. 
There is an upper bound to the folding number of crescents 
depending only on $\Sigma$ and $M$. 
\end{prop}
\begin{proof}
By above Lemma \ref{lem:upper1}, $\alpha_{\calR}$ is in 
the $N$-neighborhood $A$ of $I_{\calR}$. We draw perpendicular 
geodesics to $I_{\calR}$ foliating a subset of $\tilde M$. 
Each geodesic must meet $\alpha_{\calR}$ eventually 
in $A$ since otherwise $\alpha_{\calR} \cup I_\calR$ do not 
form a boundary of a domain. Therefore, 
the first statement is proved. 

The second statement follows from the fact that the perpendicular 
geodesic meets $\tilde \Sigma \cap A$ since $\tilde \Sigma$ is 
properly imbedded. 
\end{proof}

\begin{cor}\label{cor:cressv} 
Suppose that $\calR$ is an outer crescent.
Then there exists a convex vertex in $\alpha_{\calR}$. 
That is, the set of vertices of $\calR$ cannot consist only of 
concave vertices and saddle vertices.
\end{cor}
\begin{proof} 
We choose a function $f$ so that $I_{\calR}$ is contained in 
the zero set and other level sets are totally geodesic. 
$f$ is bounded on $\alpha_{\calR}$ by Proposition \ref{prop:cresupper}.
Hence, there is a maximum point. By tilting the totally geodesic plane 
by a little, we obtain a strictly convex vertex. 
\end{proof}

%%% Oct 05 10:55
%**** 2.1.3.0. highest level crescents

\subsection{Highest-level crescents}
Given $\Sigma$, there is an upper bound to 
the folding-number of all crescents associated with $\Sigma$
by Proposition \ref{prop:cresupper}. 
We call the maximum the {\em highest folding number} of $\Sigma$. 
We perturb $\Sigma$ to minimize the highest folding number 
which can change only by $\pm 1$ under perturbations.
After this, the folding number is constant under 
small perturbations of $\Sigma$. 
If there are no crescents, then 
the {\em folding number} of $\Sigma$ is defined to be $-1$.

Also, the union of all crescents for $\tilde \Sigma$ is 
in a uniformly bounded neighborhood of $\tilde \Sigma$
with the bound depending only on $\Sigma$. 

%%%% April 25 11:15

We say that a $0$-folded crescent $\calR$ is 
a {\em highest-level} crescent 
if it is an innermost crescent of 
an $n$-folded crescent ${\calR}'$ where
$n$ is the highest-folding number of $\tilde \Sigma$ 
and $\calR$ is the innermost one
that achieves the highest-level.

Suppose that $\calR$ is a compact highest-level crescent. 
Let $A_1, \dots, A_n$ be components of $I_{\calR}$ with pinched points 
removed. Recall that $I_\calR$ lies in a totally geodesic hypersurface.
The outermost pinched simple closed curve $\alpha_i$ 
in the boundary of $A_i$ has a trivial holonomy.
Since $\calR$ is of highest-level, 
$\alpha_i$ is an innermost curve itself 
or bounds some closed curves
in $\tilde \Sigma \cap \partial I_\calR$.
If each $\alpha_i$ is as in the former case,
then $\calR$ is said to be an {\em innermost ball-type crescent}, 
which is homeomorphic to a $3$-ball by the incompressibility condition for $\Sigma$.

We have the following important definition:
\begin{defn}
The {\em outer-folding number} of $\tilde \Sigma$ is the maximum of 
the ambient folding number of an outer highest-level crescents. 
The {\em inner-folding number} of $\tilde \Sigma$ is the maximum of 
the ambient folding number of an inner highest-level crescents.
The outer-folding number is $-1$ if there are 
no highest-level outer crescents. 
The inner-folding number is $-1$ if there are no
highest-level inner crescents.  
\end{defn}

By Proposition \ref{prop:cresupper}, the numbers are finite. 
The maximum of the both of the numbers are the folding number 
of $\tilde \Sigma$.

\subsection{Outer- and inner-contact points}
%***** 2.1.3.0.1. outer and inner points $\Sigma \cap I_{\mathcal R}$. 

We can classify the points of $\tilde \Sigma \cap I_{\mathcal R}$
when $\tilde \Sigma$ is well-triangulated: 
A point of it is an {\em outer-contact point} if the point 
is not in $\partial I_\calR$
and  has a neighborhood in $\tilde \Sigma$ outside
${\mathcal R}^o \cup \alpha_{\mathcal R}$; 
a point is an {\em inner-contact point} if the point
is not in $\partial I_\calR$ and 
has a neighborhood in $\tilde \Sigma$ contained in 
${\mathcal R}$. A point is either an outer-point
or an inner point or can be both. 

The following classifies the set of outer-contact points. 
(A similar result holds for the set of inner points
except for (d).)

\begin{prop}\label{prop:sigmaI}
Suppose that $\tilde \Sigma$ is well-triangulated.
For a highest-level crescent $\calR$, the intersection points of 
$I^o_{\calR}$ and $\tilde \Sigma$ are either outer-contact points 
or inner-contact points. 
The set of outer-contact points of $I_{\mathcal R}$ for 
a highest-level crescent $\mathcal R$ is one of the following: 
\begin{itemize}
\item[(a)] a union of at most three isolated points. 
\item[(b)] a union of at most one point and a segment 
or a segment with some endpoints removed. 
\item[(c)] a union of two segments with a common endpoint
with some of the other endpoints removed. 
\item[(d)] a triangle with some of the vertices or 
a boundary segments removed.  
\end{itemize}
The same statement are true for inner-contact points.
\end{prop}
\begin{proof} 
The set of outer-contact points is obviously a union 
of open cells of dimension $0$, $1,$ or $2$.
The vertices of the closure of 
each objects are the vertices of $\tilde \Sigma$. 

This follows from the general position of vertices of $\tilde \Sigma$. 
If the set of outer-contact points are union of $0$- and $1$-dimensional
objects, then (a), (b), or (c) follows.

If there is a $2$-dimensional object, then it contains an open triangle 
and there cannot be any other objects not in the closure of it. 

Either the interior of an edge is in $\partial I_\calR$ 
or it is disjoint from $\partial I_\calR$ 
since otherwise we might have four coplanar vertices. 
\end{proof}

\begin{defn}\label{defn:I} 
Given a crescent $\calR$, we define $I^O_{\calR}$ to be the $I_{\calR}$
with the pinched points, boundary points in the $I$-part hypersurface, and the segments 
and triangles in the outer-contact set as above removed.
(We don't remove the isolated points.)
\end{defn} 
Note that $I^O_\calR$ may not equal the topological 
interior $I^o_\calR$ of $I_\calR$ in the totally geodesic 
hypersurface.

\begin{rem}
We remark that in cases (b), (c), (d), 
the set of outer-contact points (inner-contact points) 
can separate $I_\calR$. The set is a {\em disconnecting set of 
outer-contact points}. 
\end{rem}

\subsection{Closed geodesics and crescents}
%***** 2.1.3.0.2 closed geodesics and crescents. 

\begin{prop}\label{prop:avoidg}
Suppose that $c$ is a straight closed geodesic in $M$ 
not meeting $\Sigma$. 
Let $\calR$ be a highest-level crescent.
Then 
\begin{itemize}
\item Each component $\tilde c$ of
the inverse image of $c$ in $\tilde M$ 
does not meet $\calR$ in its interior 
and the $\alpha$-parts. 
\item $\tilde c$ could meet $\calR$ in its $I$-part
tangentially and hence be contained in the $I$-part.
In this case, $\calR$ is not compact.
\item If $l$ is a geodesic in $\tilde M$ 
eventually leaving all compact subsets, 
then the above two statements hold as well.
In particular, this is true if $\tilde M$ is a special hyperbolic 
manifold and $l$ ends in the limit set of 
the holonomy group associated with $\tilde M$. 
\end{itemize}
\end{prop}
\begin{proof}
If $\tilde c$ meets the $\alpha$-part of $\calR$, then 
$\tilde c$ meets the interior of $\calR$.

If a portion of $\tilde c$ meets the interior of $\calR$, then 
$\tilde c \cap \calR$ is a connected arc, say $l$
since $\calR$ is a closure of a component cut out by 
a totally geodesic hyperplane in $\tilde M - \tilde \Sigma$ --(*).

Since $\tilde c$ is disjoint from $\tilde \Sigma$, 
both endpoints of $l$ must be in $I^O_{\calR}$ or 
in $\SI^2_\infty \cap \ovl{I_{\calR}}$ for 
the closure $\ovl{I_{\calR}}$ of $I_{\calR}$ in 
the compactified $\hyps \cup \SI^2_\infty$. 
If at most one point of $l$ is in $I_{\calR}$, then 
$l$ is transversal to $I_{\calR}$ and the other endpoints 
$l$ must lie in $\alpha_{\calR}$ by (*). This is absurd. 

If at least two points of $l$ are in $I_{\calR}$, then 
$l$ is a subset of $I_{\calR}$. Since $\tilde c$ is disjoint from
$\tilde \Sigma$, it follows that $\tilde c$ is a subset of $I_{\calR}$, 
and $\calR$ is not compact. 

The only remaining possibility for $l$ is the third one
that $l$ ends in $\SI^2_\infty$, $l$ is a subset of $\calR^o$, $\tilde c = l$,  
and $\calR$ is noncompact.
Suppose that $l$ is a subset of $\calR^o$.
A point of $\calR^o$ is a point of some compact crescent $R_i$ in $\calR$ since 
a noncompact crescent is a generalized limit of compact ones. 
Therefore, $l$ meets a compact crescent as above, which was shown to be 
not possible above. This proves the first two items.

The third item follows similarly.

\end{proof}

\subsection{A highest-level crescent is included in 
a secondary highest-level crescent} 

%***** 2.1.3.0.3 Fundamental highest-level crescent Proposition
\begin{prop}\label{prop:alphapt} 
Let $\tilde \Sigma$ be well-triangulated
and $\calR$ be a highest-level crescent. Then 
there exists an innermost crescent 
$\calR'$ containing $\calR$ so that 
\begin{itemize}
\item $I_{\calR'}$ is connected 
and has no pinched points or a disconnecting set of 
outer-contact points, i.e., $I^O_{\calR'}$ is connected.
\item The closure of $\alpha_{\calR'}$ is homeomorphic 
to $I_{\calR'}$ and is isotopic to $I_{\calR'}$ in $\calR'$. 
\item $\calR'$ is homeomorphic to 
$I_{\calR'} \times [0,1]$.
\item $\calR'$ is inner-most.
\end{itemize}
\end{prop}
\begin{proof}

We assume the hypothesis \ref{hyp:incompressible}.

The basic idea is to add some domains by disks in $\tilde \Sigma$ 
obtained by incompressibility of $\tilde \Sigma$. 

By Proposition \ref{prop:avoidg}, 
the interior of $\calR$ is disjoint from any 
lifts of $\cl_1, \dots, \cl_n$. 

{\em The first step is to cut off by the $I$-part hypersurface to 
simplify the starting crescent}\/: 
Suppose that $\calR$ is compact to begin with.
Then $\calR$ is disjoint from the lifts of $\cl_1, \dots, \cl_n$
since $\tilde \Sigma$ is disjoint from these. 
We let $\calS$ be a crescent obtained 
from $\calR$ by cutting through the $I$-part hypersurface $P$
and taking the closure of a component of $\calR - \tilde \Sigma - P$.
We say that $\calS$ is a {\em cut-off} crescent from $\calR$.
(The ambient folding number and the folding number may change by 
cutting off.)

Again $\calS$ is innermost since otherwise we will have 
points in $\calS$ in the other side of $\tilde \Sigma$ than those of $\calR$ 
but $\calS \subset \calR$. 

The ambient level of $\calS$ is equal to that of $\calR$: 
To reach a point in $\calR$ from any point of $\alpha_{\calR'}$ for 
a super-crescent $\calR'$, one needs to traverse on a generic path
at least the highest-level times of $\hat \Sigma$ in $\calR'$. 
$\calS$ is a sub-crescent of the cut-off crescent $\calS'$ of $\calR'$. 
Since a path in $\calS'$ is a path in $\calR'$ and 
the ambient level is a maximum value, 
the ambient level may increase. 
However, since we are already at the highest-level, the equality holds. 

We introduce a height function $h$ on $\calS$ defined by 
introducing a parameter of hyperbolic hypersurfaces perpendicular 
to a common geodesic passing through $I_{\calS}$ 
in the perpendicular manner. (It will not matter which 
parameter we choose).
We may assume that $h$ is Morse in the combinatorial sense. 
(See Freedman-McMullen \cite{FM}.)

{\em Now, we fatten up $\calS$ a bit so that we can work 
with surfaces instead of just topological objects}\/:
If $I_\calS$ does not meet any lifts of $\cl_1, \dots, \cl_n$, 
let $N_\eps(\calS)$ be the neighborhood of $\calS$ in
the closure of the component of $\tilde M - \tilde \Sigma$ 
containing the interior of $\calS$. 

We let $N_{\eps}(\alpha_\calS)$ to be the intersection of 
$\tilde \Sigma$ with $N_{\eps}(\calS)$.
$N_\eps(\calS)$ can be chosen so that $N_\eps(\alpha_\calS)$ 
in $\tilde \Sigma$ becomes an open surface compactifying 
to a surface. There exists a part $\calI$ in of the boundary 
which is a complement in the boundary of $N_\eps(\calS)$ 
of $N_\eps(\alpha_\calS)$ and lies on
a properly imbedded surface $P'$ perturbed away from 
the $I$-part hypersurface $P$ of $\calS$. 

Topologically, $N_\eps(\alpha_\calS)$ is homeomorphic 
to a surface possibly with $1$-handles attached
from $\alpha_{\calS}$ and $\calI$ is obtained from $I_\calS$ by 
removing $1$-handles corresponding to the pinched points 
or the disconnecting set of outer-contact points.

{\em Now we aim to show that $N_\eps(\calS)$ is a compression body 
with $N_\eps(\alpha_\calS)$ as the compressible surface in 
the boundary}\/:

We may extend $h$ to an $\eps$-neighborhood of 
$\calS$, which may introduce only saddle type singularity 
in $N_\eps(\alpha_\calS) - \alpha_\calR$ where there are only
one handles. 
We modify $h$ so that $\calI$ to be in the zero level of $h$ and  
$h< 0$ in the interior of $\calS$. 
%This process only introduces a saddle type singularity in 
%$N_\eps(\alpha_\calS)$.

We show that there are no critical points of positive type for $h$: 
If there is a critical point of $h$ with locally 
positive type where $h < 0$, then we see that in fact 
there exists a crescent of higher level near the critical 
point. The critical point is actually below the original
$I$-part hypersurface since only critical points above the $I$-part 
hypersurface are of saddle type. 
The crescent is obtained by a totally geodesic hyperplane 
containing the level set slightly above the critical point.
The level of the new crescent is one more than that of $\calS$,
which is greater than the highest level, which is absurd. 

Since there are no critical point of positive type, 
$\pi_1(N_{\eps}(\alpha_{\calS})) \ra \pi_1(N_\eps({\calS}))$
is surjective as shown by Freedman-McMullen \cite{FM}. 
There exists a compression body in $N_{\eps}(\calS)$ 
with a boundary $N_{\eps}(\alpha_{\calS}) \cup S'$ 
for an incompressible surface $S'$ in the interior of $\calS$.
Since every closed path in $N_{\eps}(\calS)$ is 
homotopic to one in $N_{\eps}(\alpha_\calS)$, it follows that 
$S'$ is parallel to $I$. Hence, $N_{\eps}(\calS)$ is 
a compression body homeomorphic to $\calI$ times an interval 
and $1$-handles attached at disks disjoint from $\calI$. 
($\calS$ is essentially obtained by pinching some points of $\calI$
together and pushing down a bit.) 

{\em Next, we reduce the number of components of $\calI$}\/:

Suppose now that $\calI$ is not connected. 
This means that there are $1$-handles attached to $\calI$ times an interval 
joining the components below $\calI$. 
Then $N_{\eps}(\calS)$ has a compressing disk $D$ 
for $N_{\eps}(\alpha_\calS)$ dual to the $1$-handles.
Since $\partial D$ bounds a disk $D'$ in $\tilde \Sigma$ by the 
incompressibility of $\tilde \Sigma$, the closed curve
$\partial D$ is separating in $\tilde \Sigma$.
Consequently also, $N_\eps(\alpha_{\calS})$ is a planar surface. 

The irreducibility of $\tilde M$ tells us that 
$D$ and $D'$ bound a $3$-ball $B$ in the closure of a component of 
$\tilde M - \tilde \Sigma$. Then $B$ contains at least one 
component of $I$. 

\begin{figure}[ht]
\centerline{\epsfxsize=3.8in \epsfbox{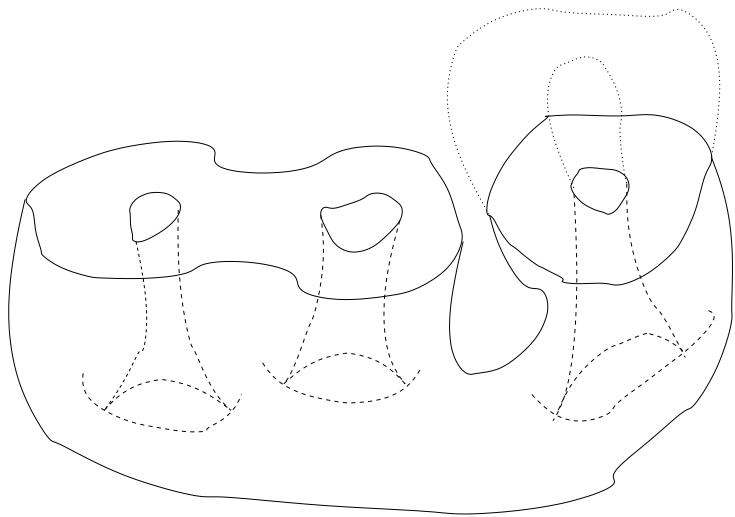}}
\caption{\label{fig:hlcrescent} The dashed arc indicates 
the tube from the bottom and the dotted arcs indicated 
the disks to be attached to the $\alpha$-parts}
\end{figure}
\typeout{<>}

By taking a maximal family of compressing disks dual to the 
$1$-handles and regarding the components of the complements as 
vertices, we see that the $1$-handles do not form a cycle. 
Therefore, we choose the compressing disk $D$ of $N_\eps(\alpha_\calS)$ 
to be the one such that $D$ and corresponding disk $D'$ in $\tilde \Sigma$ 
bounds a $3$-ball $B$ containing a unique component of $I$. 

We take a union of $B$ with $N_\eps(\calS)$. Then 
it is an $N_\eps$-neighborhood of a crescent $\calS'$,
which contains $\calS$ and $\alpha_{\calS'}$ 
containing $\alpha_\calS$ and $I_{\calS'}$, a subset of $I_\calS$. 

In fact, $\calS'$ is a union of $\calS$ and $B'$ with 
some parts in $\calI - \calI'$ removed.
Also, clearly, $N_\eps(\calS')$ is a compression body 
since it is obtained by taking a union of a cell with 
a compression body obtained from $N_\eps(\calS)$ 
by splitting along $D$.

By induction, we obtain a crescent $\calR''$ with 
$I_{\calR''}$ in $I_\calS$ and the surface $I'$ corresponding to 
$I$ connected. We say that $\calR''$ is {\em derived from $\calS$}.  
Since $\calR''$ is homeomorphic to a compression body,
$\calR''$ is homeomorphic to 
$I$ times an interval since there are no compressing disks. 

If there are any pinched points in $I_{\calR''}$ or 
disconnecting outer-contact points, then $I'$ would be disconnected.
$\alpha_{\calR''}$ has a closure that is a surface since 
there are no pinching points.

Since $\calR''$ is an $I$-bundle, it follows that 
the closure of $\alpha_{\calR''}$ and $I_{\calR''}$ are 
homeomorphic surfaces. 

Also, $\calR''$ is innermost: 
suppose not. Then there exists a component $C$ of $\calR'' - \tilde \Sigma$
so that a generic path in $\calR''$ from $\alpha_{\calR''}$ may meet 
$\tilde \Sigma$ more than once. Then $C \cap \calR$ is again 
a component of $\calR - \tilde \Sigma$, which is a contradiction.

%%%% April 4 11:30 some thinking? 

%%%% March 31 3:45 

{\em Now, we go to the final step}.
Recall that $\calS$ was a cut-off crescent from the original $\calR$. 
If there were more than one cut-off crescents $\calS$, then 
we obtain $\calR''$ for each $\calS$. 
Suppose that two cut-off crescents 
$\calS$ and $\calS'$ adjacent 
from opposite sides of some of the components of $I_\calS$.
Since the corresponding $\calR''$ and $\calR'''$ containing 
$\calS$ and $\calS'$ respectively  
does not have any pinched points or separating outer-contact edges, 
the unique components of $I_{\calR''}$ and $I_{\calR'''}$ either agree or 
are disjoint from each other. 
$\calR''$ and $\calR'''$ cannot be adjacent from opposite side
since we can then form a compact component of 
$\tilde \Sigma$ otherwise. It follows that one of 
$\calR''$ and $\calR'''$ is a subset of the other. 

Hence, choosing maximal ones among such derived crescents, 
we see that the conclusions of the proposition
hold if $\calR$ is compact. The final result 
is a product of its $I$-part since
the final compression body has no $1$-handles. 
This completes the proof in case $\calR$ is compact.

%Noncompact case
If $\calR$ is noncompact, we follow as before but
we choose $N_{\eps}(\calR)$ to be
tapered down near infinity. 

Recall that $\tilde M$ is $2$-convex, i.e., 
the boundary $\partial \tilde M$ has only convex or 
saddle vertices. Since $\tilde M$ can be considered a $2$-convex 
affine manifold, recall the main result of \cite{uaf} that any disk with 
a boundary in a totally geodesic hypersurface $P$ bounds a disk 
in $P$. 

%%% 0104 4:15  I need 2-convexity for M. Where to state? 
Only one 
component of $\calI$ maybe noncompact since the boundary of 
a compressing disk must bound a compact disk in $\tilde \Sigma$: 
Otherwise, we have a simple closed curve $c$ in $\alpha_\calR$ which 
separates the two noncompact components of $\calI$ and 
bounds a compact disk $D$ in $\tilde M$ in one side of $\tilde \Sigma$. 
$D$ can be pushed inside $\calR$ by the above paragraph. 
Thus $D$ cannot meet any of the lift of $\cl_1, \dots, \cl_n$. 
Since $\tilde \Sigma$ is incompressible in the complement of these lifts, 
it follows that $c$ also bounds a compact disk in $\tilde \Sigma$. 
This is a contradiction since $c$ separates the two noncompact 
components of $\calI$. 

Going a back to finding the ambient crescent with desired properties, 
components of $\calI$ involved in cell-attaching operation 
as above are compact. So the arguments for noncompact $\calR$ 
are the same as above. 
\end{proof}

%\begin{defn} 
%We say that $\calS$ as obtained from $\calR$ in 
%the above lemma is obtained by cutting along 
%the $I$-part hypersurface of $\calR$.  
%$\calS$ may not be a unique one so obtained.
%\end{defn}

%%%% Dec. 17 11:45

\subsection{Properties of secondary highest-level crescents}
%***** 2.1.3.0.4. Secondary highest level crescents

If a $0$-folded crescent $\calS$ contains a highest-level crescent 
$\calR$ so that $I^O_\calS$ is 
connected and is included in $I^O_\calR$, then 
we say that $\calS$ is a {\em highest-level} crescent as well. 
(Actually, it may not be highest-level since $\calS$ 
may not necessarily be contained in an $n$-folded crescent
but only a part of it.) 
More precisely, it is a {\em secondary highest-level crescent}.

\begin{cor}\label{cor:avoidg} 
Let $\calR$ be a secondary highest-level crescent.
Then the statements of Proposition \ref{prop:avoidg} 
hold for $\calR$ as well.
\end{cor}
\begin{proof}
The proof is exactly the same as that of 
Proposition \ref{prop:avoidg}.
\end{proof}

By taking a nearby crescent inside by changing the $I$-part hypersurface 
inwards, we see that a highest-level crescent could be {\em generically chosen}
so that the crescent is compact,
the $I$-part and the $\alpha$-part are surfaces,
and $I^O$-part is truly the interior of the $I$-part. 

\begin{cor}\label{cor:highestcrestop} 
Let $\calR$ be the compact secondary highest-level 
outer {\rm (} resp. inner {\rm )} crescent that is generically chosen.
Then $\calR$ is homeomorphic to the closure of $\alpha_{\calR}$ 
times $I$, and $I^O_{\calR}$ is isotopic to $\alpha_{\calR}$ 
by an isotopy inside $\calR$ fixing the boundary of $I_{\calR}$. 
\end{cor}

\subsection{Intersection properties of highest-level crescents}
%%% Oct 05 3:30

%**** 2.1.3.1. transversal intersections
%***** 2.1.3.1.1. No face off 

We say two crescents $\calR$ and $\calS$ {\em face each other} 
if $I_{\calR}$ and $I_{\calS}$ agree with each other in some 
$2$-dimensional part and have disjoint one-sided neighborhoods.

\begin{prop}\label{prop:nofacing} 
Assume $\tilde \Sigma$ is well-triangulated as above.
If two highest-level outer 
{\rm (}resp. inner{\rm )} crescents $\calR$ and $\calS$ face 
each other, then there exists a {\rm (}secondary{\rm )} highest-level
outer {\rm (}resp. inner{\rm )} crescent $\calT$ 
with connected $I^O_{\calT}$ containing both. 
\end{prop} 
\begin{proof} 
We replace $\calR$ and $\calS$ by secondary highest-level crescents
with connected $I^O$-parts. The replacements
still face each other or one becomes a subset of 
another since $I$-parts are unique boundary sets. 
In the second case, we are done. 

Since $I^O_\calR$ and $I^O_\calS$ meet in open subsets, 
either they are identical or we may assume without 
loss of generality that the boundary $\partial L$ of 
their intersection $L$ in $I^O_\calR$ is not empty.
$\partial L$ is a subset of $\tilde \Sigma$ and 
is a $1$-complex consisting of pinched arcs.
$\partial L$ is a set of outer-contact points of $\calS$ 
since a neighborhood of $\partial L$ in $\tilde \Sigma$ must 
be above $\calS$. However, then $\partial L$ must 
be disjoint from $I^O_\calS$ since $I^O_\calS$ is 
disjoint from outer-contact set by definition.
Therefore, $I^O_\calS = I^O_\calR$, and  
$\calS \cup \calR$ is bounded by 
a component subsurface of $\tilde \Sigma$,
which is absurd.

\end{proof}

%***** 2.1.3.2.2. Transversal intersection Defn and Proposition

%%%% December 17 4:00

\begin{defn}\label{defn:transvint} 
Two secondary highest-level outer (resp. inner) crescents 
$\calR$ and $\calS$ are said to meet 
{\em transversally} if $I_{\calR}$ and $I_{\calS}$ meet in 
a union of disjoint geodesic segment $J$, $J \ne \emp$, mapping into 
a common geodesic in $\hyps$, in a transversal manner such that 
\begin{itemize}
\item The the closure $\nu_{\calR}$ of the union of 
the components of $I_{\calR} - J$ 
in one-side is a subset of $\calS$ 
and the closure $\nu_{\calS}$ of the union of those of 
$I_{\calS} - J$ is a subset of $\calR$. 
\item The intersection ${\calR} \cap {\calS}$ is the closure of 
$\calS - \nu_{\calR}$ and conversely the closure of 
$\calR - \nu_{\calS}$. 
\item The intersection $\alpha_{\calR} \cap \alpha_{\calS}$ is 
a union of components of $\alpha_{\calR} - \nu_{\calS}$ in 
one-side of $\nu_\calS$
and, conversely, is a union of components of 
$\alpha_{\calS} - \nu_{\calR}$ in one side of $\nu_\calR$.
\item $\alpha_{\calR} \cup \alpha_{\calS}$ is an open surface in 
$\tilde \Sigma$. 
\end{itemize}
\end{defn}

%%% 0301 1:30

\begin{prop}\label{prop:transvint} 
Given two secondary highest-level 
outer (resp. inner) crescents $\calR$ and $\calS$, 
there are the following mutually exclusive possibilities: 
\begin{itemize}
\item $\calR$ and $\calS$ do not meet in $\tilde M -\tilde \Sigma$. 
\item $\calR \subset \calS$ or $\calS \subset \calR$. 
\item $\calR$ and $\calS$ meet transversally. 
\end{itemize}
\end{prop}
\begin{proof} 
The reasoning is exactly the same as \cite{cdcr1} and \cite{psconv} 
in dimension two or three.
\end{proof}

\section{The crescent-isotopy}\label{sec:crescent-iso}

The purpose of this section is to prove Theorem C and Corollary D:
Assume that $\Sigma$ is a closed well-triangulated surface in $M$ which is incompressible 
in $M$ with a number of closed geodesics removed.
In this section, 
we will describe our crescent-isotopy steps of $\tilde \Sigma$.
Let $\tilde \Sigma$ have a folding number $n$ achieved by 
outer and/or inner crescents. 
We may assume without loss of generality that there is an outer 
highest level crescent coming from an outer or inner crescent.
Using such outer crescents, we move first to reduce 
the outer level by $1$. 
\begin{description}
\item[Subsection \ref{subsec:smtr}] The first step is to truncate our surface
$\tilde \Sigma$ along vertices, edges or triangles 
in order to make highest-level crescents not 
have outer-contact points.
This may make $\tilde \Sigma$ only triangulated; however, 
we make small perturbations of vertices to make it well-triangulated. 
\item[Subsection \ref{subsec:criso}] We use the secondary 
highest-level crescents to move isotopy $\tilde \Sigma$ by isotopying 
the closure of the $\alpha$-parts to the $I$-parts.
One of the outer or inner levels strictly decreases.
\item[Subsection \ref{subsec:convpert}]   
The result may have some parts which are pleated with infinitely 
many and/or infinitely long pleating geodesics. We perturb these parts so
that we end up with a triangulated surface but with levels not increasing
\item[Subsection \ref{subsec:thmC}] We do the above for the level $n$ for inner highest-level 
crescents. This will decrease the inner level. 
(The steps are just the same if we reverse the orientation of 
$\tilde \Sigma$.) We keep doing this until our outer and inner level become 
$-1$ and we have obtained an saddle-imbedded surface 
isotopic to $\tilde \Sigma$, which proves Theorem C. 
Finally, we will prove Corollary D. 
\end{description}

\subsection{Small truncation moves}\label{subsec:smtr}
%*** 2.1.4 crescents and isotopy 

%**** 2.1.4.1. A simple lemma on crescent isotopy
\subsubsection{Isotopies}
First, we need:
\begin{lem}\label{lem:crescentiso}
Suppose that $\tilde \Sigma$ has been isotopied in 
the outward direction by a sufficiently small amount
and $\calR$ is an outer crescent. Then there exists 
a crescent $\calR'$ sharing the $I$-part hypersurface with 
$\calR$ and differs from $\calR$ by isotopying the 
$\alpha$-part only. 
Conversely, if $\tilde \Sigma$ has been isotopied in the inward direction 
and $\calR$ is an inner crescent, the same can be said. 
\end{lem}
\begin{proof} 
Straightforward.
\end{proof}

We say that $\calR'$ is {\em isotopied from $\calR$ with the 
$I$-part preserved}. (Of course, this is not literally so.)

\subsubsection{Small truncations}\label{subsubsec:smalltrunc}

We may ``truncate" $\Sigma$ at convex vertices and 
$\tilde \Sigma$ correspondingly and perturb: 
Let $v$ be a convex or concave vertex and $H$ a local half-open ball 
at $v$ containing $\tilde \Sigma$ locally 
with the side $F$ passing through $v$. 
\begin{itemize}
\item[(a)] We may move $F$ inside by a very small amount 
and then truncate $\tilde \Sigma$ using the displaced $F$ 
and add the trace disk $T$ 
of the truncation to the surface $\tilde \Sigma$.
\item[(b)] Then we introduce some equivariant triangulation 
of $T$ of the truncation and the truncated $\tilde \Sigma$
without introducing vertices in the interior of $T$. 
\item[(c)] We will have to do this for each vertex which is in the orbit
of $v$ so that resulting $\tilde \Sigma$ is still equivariant. 
\item[(d)] Finally, we perturb all the vertices of $\tilde \Sigma$
by a sufficiently small amount. Here, the perturbations 
must be so that the normal vectors to the totally geodesic 
triangles also move by small amounts, i.e., the normal 
vectors move continuously as well as the vertices themselves. 
Moreover, no triangle or edge degenerates to 
a lower-dimensional object.
\end{itemize}
The three steps (a)-(c) together are called the {\em small truncation move}.
Together with the final step (d), 
the move is called the {\em perturbed small-truncation move}.

For an edge or a triangle $e$, let $F$ be a neighborhood of $e$ in 
totally geodesic plane containing $e$ where $F - e$ lies 
outside $\tilde \Sigma$. We may move $F$ inside by a sufficiently 
small amount and truncate $\tilde \Sigma$. 
The rest is similar to the vertex case. 
They are also called {\em small truncation moves}
along edges or triangles.
After the perturbation, we call the move 
{\em perturbed small-truncation move}.

We denote by $\Sigma^\eps$ the perturbed $\Sigma$ where 
the trace disks are less than an $\eps$-distance away from the 
respective convex vertices and the normal vectors to
the triangles are also less than $\eps$-distances from the original 
vertices. 
Here, we assume that during the perturbations $\Sigma^\eps$ is 
isotopied from $\Sigma$ and the convexity of the dihedral
angles do not change under the isotopy. 
Thus, if an edge or a vertex is convex after being born, it will
continue to be so as $t \ra 0$ and as $t$ grows from $0$.

We may also assume that the convex vertex move is 
equivariant on $\tilde \Sigma$, i.e., the isotopy 
is equivariant. 

\subsubsection{Small truncation moves and crescents}
\label{subsubsec:smtrunc}

An {\em isotopy} of a crescent as we deform $\Sigma$ 
is a one-parameter family of crescents $\calR_t$ with 
$\alpha$-parts in $\Sigma$.
The above small truncation moves are isotopies.

We say that a crescent {\em bursts} if fixing the totally geodesic 
hypersurface containing the $I$-part of it and isotopying
the $\alpha$-parts in the isotopied $\Sigma$ 
cannot produce a crescent isotopied from the original one. 

Such an event happens when 
a parameter of vertices, edges, or triangles of $\tilde \Sigma$ 
go below the fixed totally geodesic hypersurface from the 
point of view of the crescent. 
Of course, a vertex could be a multivertex and 
all of the new vertices go down. 
The edge should be on the face that meets the $I$-parts of
the crescents and the vertex on the edge that meets the 
$I$-part of the crescent.
The event could happen simultaneously but 
the generic nature of the move shows that 
at most four vertex submersions, at most three edge submersions, 
at most two vertices and one edge submersions, or 
triangle-, edge- or vertex-submersions can 
happen simultaneously. (Basically, at most four vertices can 
lie on a totally geodesic plane while deforming.)
Moreover, at the event, 
the vertex and the edge must be in the $I$-part of 
the crescent and the triangles of $\Sigma$ must be placed 
in certain way in order that the bursting to take place.

\begin{prop}\label{prop:perturb}
Suppose that $\tilde \Sigma$ is well-triangulated. 
Under a small truncation move in the outer direction, 
we can isotopy $\tilde \Sigma$ to triangulated $\tilde \Sigma^\eps$
{\rm (}equivariantly{\rm )} so that
\begin{itemize}
\item[(i)] each outer crescent moves into itself by moving 
the $\alpha$-part in the outer direction
and preserving the $I$-part hypersurface. 
\item[(ii)] each inner crescent moves into itself union 
the $\eps$-neighborhood of $\tilde \Sigma$ by moving 
the $I$-part hypersurface in the outer direction or 
preserving the $I$-part hypersurface.
\end{itemize} 
Under a small truncation move in the inner direction, 
we can deform 
\begin{itemize}
\item[(iii)] each inner crescent into itself by moving 
the $\alpha$-part in the inner direction and 
preserving the $I$-part hypersurface.
\item[(iv)] each outer crescent into itself union 
the $\eps$-neighborhood of $\tilde \Sigma$ by moving 
the $I$-part hypersurface in the inner direction
or preserving the $I$-part hypersurface. 
\end{itemize} 
All crescents of $\tilde \Sigma^\eps$
can be obtained in this way.
The highest folding number may decrease only
under a convex vertex move, and the union of 
crescents of all levels strictly decreases under the moves. 
\end{prop}
\begin{proof} 
Essentially, the idea is that the move can only ``decrease" 
the associated crescents. 

Let $\calR$ be an outer crescent and $\tilde \Sigma$ moved 
in the outer direction. Lemma \ref{lem:crescentiso} implies (i).

Let $\calR$ be an inner crescent and $\tilde \Sigma$ be moved 
in the outer direction. Then again an isolated submerging 
vertex is a convex vertex.  In this case, we move the 
$I$-part inward so that the submerging vertex stay on
the boundary of the $I$-part. Other cases are treated similarly.
This proves (ii).

(iii) and (iv) correspond to (i) and (ii) respectively
if we change the orientation of $\Sigma$.

To show that all crescents of $\tilde \Sigma^\eps$
can be obtained in this way:
Given an outer crescent for $\tilde \Sigma^\eps$, 
we reverse the truncation move. 
If the $I$-part of a crescent avoids the 
trace disks of the truncation moves, then we simply
isotopy the $\alpha$-parts only. 

%%% December 29 5:10. Eliminate more choppy sentences...

The trace surface has only concave vertices and 
saddle-vertices. 

Let us start by reversing a vertex truncation move: 
Let $P'$ be a local totally geodesic hypersurface 
truncating the stellar neighborhood of a convex vertex 
$v$ of $\tilde \Sigma$ at some small distance from $v$ 
but large compare to our isotopy move distance.
Suppose that $v$ were involved in the convex truncation move.  
We may assume that $P'$ is parallel to the truncating totally
geodesic hyperplane near $v$ used to obtain $\tilde \Sigma^\eps$.

Clearly, $P'$ and a small stellar neighborhood of $v$ in 
$\tilde \Sigma^\eps$ bounds a small polyhedron $R^\eps$. 
Let $R$ be the small polyhedron bounded by $P'$ and 
$\tilde \Sigma$. 

Suppose that the $I$-part of a crescent $\calR$ for 
$\tilde \Sigma^\eps$ are contained in $R^\eps$. 
Then it is contained in a crescent whose $I$-part 
meets what are outside the part truncated by $P'$. 
It is sufficient to show that the ambient crescent 
is obtained by the above methods.

We thus assume without loss 
of generality that 
the $I$-part of a crescent $\calR$ for 
$\tilde \Sigma^\eps$ 
meets what are outside the part truncated by $P'$.
If the $I$-part does not meet the trace surface, 
then we only change the $\alpha$-parts to obtain 
a crescent for $\tilde \Sigma$ as above. 
We may assume that the $I$-part of $\calR$ meets 
the trace surface without loss of generality.
 
Assuming that our isotopy was very small, 
since the $I$-part meets one of the trace surface and $P'$ is separating,
the $I$-part meets $P'$.
$P'$ intersected with the closure of the exterior of
$\tilde \Sigma^\eps$ is a polygonal disk $D^\eps$.
Then $D^\eps$ intersected with the $I$-part is a disjoint union
of segments. 

Extending the $I$-part of $\calR$ in $R$ 
until they meet unperturbed $\tilde \Sigma$, 
the set of points in $\calR$ extends in $R$ into the polyhedrons 
bounded by $P'$ and the stellar neighborhood of $\tilde \Sigma$.

Since all vertex submersions of $\tilde \Sigma^\eps$ 
can happen by vertices near the convex vertices of
$\tilde \Sigma$ masked off by 
totally geodesic hypersurfaces such as $P'$,
we obtain a crescent $\calR'$ for $\tilde \Sigma$
preserving the $I$-part hypersurface of $\calR$. 

Therefore, $\calR$ were obtained from $\calR'$
by the convex truncation isotopy preserving the $I$-part.

For $\tilde \Sigma^\eps$ obtained from $\tilde \Sigma$ by
small truncations along edges and triangles, very similar 
arguments using totally geodesic planes as $P'$ parallel to those
used in the truncation process will show the desired results.

Therefore a crescent for $\tilde \Sigma^\eps$ is one 
we obtained by the process in (i). 

Let $\calR$ be an inner crescent for $\tilde \Sigma^\eps$.
Then since the vertices moved outward with respect to 
$\tilde \Sigma$, they move inward when we reverse the process and 
we see that $\calR$ is isotopied to a crescent for
$\tilde \Sigma$ by Lemma \ref{lem:crescentiso}
by preserving the $I$-part hypersurface.

To show that the highest folding number can only decrease:
For a crescent $\calR$ to increase the folding number, 
a vertex must move into $I_{\calR}$ during the isotopy. 
We see that such a vertex must be a convex one. 
However, the convex vertex can only move in the 
direction away from the interior of $\calR$. (Even ones 
after the births obey this rule.)

Also, we can do this construction for $\tilde \Sigma$ 
simplex by simplex so that the final result is an equivariant 
isotopy.
\end{proof}

\subsubsection{Perturbations}\label{subsubsec:perturb}

\begin{defn}\label{defn:perturb} 
By a perturbation of a triangulated surface, we mean 
the perturbations of vertices and corresponding edges and faces 
accordingly. 
\end{defn}
By definition, there cannot be generations of edges and 
triangles to lower-dimensional objects under perturbations.
This applies the the proof of 
Proposition \ref{prop:perturbation}:

\begin{prop}\label{prop:perturbation} 
Suppose that $\tilde \Sigma$ does not have 
outer-contact points for its outer highest-level crescents. 
Then for sufficiently small perturbation, the 
level of the outer highest-level crescents of $\tilde \Sigma$ 
does not increase. 
Moreover, the outer highest-level crescents of the perturbed $\tilde \Sigma$ do not 
have outer-contact points. 
\end{prop} 
\begin{proof}
The $I$-part of the crescent set has an image with 
a compact closure in the quotient manifold $M$.
Let $\tilde \Sigma^\eps$ be an equivariantly perturbed surface
parameterized by $\eps > 0$ where $\tilde \Sigma = \tilde \Sigma^\eps$.

Given $\delta > 0$, we can find $\eps >0$ so that
the union of all crescents of $\tilde \Sigma^\eps$ 
is in a $\delta$-neighborhood of that of $\tilde \Sigma$: 
If not, there exists a sequence of $R_i$ for 
$\tilde \Sigma^{\eps_i}$ 
converging geometrically 
to a crescent $R$ for $\tilde \Sigma$ where 
$R$ is not in the $\delta$-neighborhood of the union of 
all crescents. This is a contradiction.

Suppose that a level increased, say to higher than or equal to $n+1$, 
for a highest level outer crescent
at $\tilde \Sigma^\eps$ for some $\eps > 0$ if the level of $\tilde \Sigma$ 
were $n$ for an integer $n \geq -1$. 
Then by acting by deck transformation to put the level $(n+1)$-crescents to intersect 
the nearby fundamental domains of $\tilde \Sigma^{\eps}$s,
we can find a sequence of level $n+1$ inner-most crescent $R_i$ 
for $\tilde \Sigma^{\eps_i}$ converging to a some subset in $\tilde M$ 
as $i \ra 0$ where $\eps_i \ra 0$. 

We may assume without loss of generality that 
\begin{itemize}
\item There exists a sequence of points $p_i \in R_i$ converging to 
a point $p$ in $\tilde M$ where $R_i$s have level $n+1$. 
\item Let $P_i$ be the totally geodesic hyperplane containing $I_{R_i}$. 
$P_i$ has a geometric limit in a totally geodesic hyperplane $P$ and the normal
vectors to $P_i$ converges to that of $P$. 
\item $R_i$ converges geometrically to a closed subset in 
$\tilde M$ with no interior since otherwise there is a generalized limit 
crescent of level $n+1$. 
\end{itemize}

A path $\gamma_i$ from a point of $\alpha_{S_i}$ for the largest 
ambient crescent $S_i$ to $p_i$ must pass at least $n+1$ components of 
$\tilde \Sigma^{\eps_i}$. A subsequence of the closure 
of the component of $S_i - \tilde \Sigma^{\eps_i}$ 
must converge to a closed subset with no interior and 
and hence a subset of $\tilde \Sigma$. 
There is a point $x_i \in \tilde \Sigma^{\eps_i}$ on the path 
$\gamma_i$ very close to $p_i$ and the path in $\gamma_i$ 
between $x_i$ and $p_i$ does not meet $\tilde \Sigma^{\eps_i}$. 
A neighborhood of $x_i$ in $\tilde \Sigma^{\eps^i}$ is 
a union of triangles and they must all converge 
to triangles outside the interior of $R$ at the end. 
Since our perturbation is assumed to be very small, 
the distance on $\tilde \Sigma^\eps$ from $x$ and points of 
$\clo(\alpha_{R_i})$ does not change much as $\eps \ra 0$, 
the distance on $\tilde M$ from $x$ to $\clo(\alpha_{R_i})$ is bounded 
above as well.
Since $\alpha_R$ is a subset of the limit of a subsequence of $\clo(\alpha_{R_i})$
by Proposition \ref{prop:alphalimit},
$x_i$ converges to a point $x$ in $I^o_R$, away from the boundary of $\partial I_{R_i}$.
The triangles of $\tilde \Sigma$ containing $x$ are outside 
the interior of $R$ since otherwise we must have level increased. 
Therefore, there is an outer-contact point $x$, a contradiction. 
(See figure \ref{fig:ocp}.)

\begin{figure}[h]
\centerline{\epsfxsize=3.5in \epsfbox{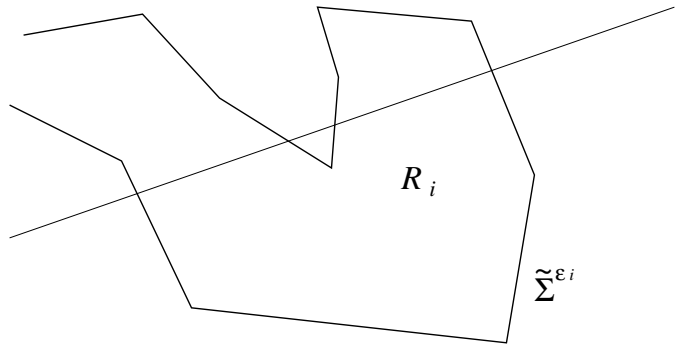}}
\caption{\label{fig:ocp} Movements of crescents.}
\end{figure}
\typeout{<>}

Finally, if $\tilde \Sigma^\eps$ has outer-contact points 
for all small $\eps >0$, we have a sequence of outer-contact
points $p_i$ for highest-level, i.e., level $n$, 
crescents $R_i$ for $\tilde \Sigma^{\eps_i}$ where $\eps_i \ra 0$. 
If $R_i$ does not degenerate, then, 
similarly to above, we obtain an outer-contact point for 
a crescent $R$ of $\tilde \Sigma$.

If the sequence of $R_i$ degenerates, then a sequence of triangles in the closure 
of their $\alpha$-part converges to a triangle which contains a limit point of $p_i$
as well. Since the distance of $p_i$ to the triangle on $\tilde \Sigma^{\eps_i}$ 
are bounded below as shown above, 
this contradicts the imbedded property of $\tilde \Sigma$.
\end{proof}

The final result of Subsection \ref{subsec:smtr} is:
\begin{app}\label{app:smtr}
As a consequence, we now perturb $\tilde \Sigma$ after
the small truncation moves for outer-contact points 
for highest-level outer crescents, and the result do not have outer-contact points and 
have the outer-level kept same or decreased.  
\end{app}

\subsection{The crescent isotopy}\label{subsec:criso}
\subsubsection{The crescent-sets}
\label{subsubsec:crset}
%**** 2.1.2.3 Crescent moves

First, we suppose that there are highest-level crescents 
whose innermost crescents are outer-direction ones. 
We will ``move" $\tilde \Sigma$ in the outer direction first to 
eliminate the outer highest-level crescents.

If the innermost ones are inner-direction, then 
we can simply change the orientation of $\tilde \Sigma$
and the discussions become the same.
Thus we assume the above without loss of generality.

%***** 2.1.2.3.2. highest level move itself

%****** 2.1.2.3.2.1 the replacing surface from the I-parts
As we did in \cite{cdcr1} and \cite{psconv}, 
we say that two highest-level crescents $\calR$ and $\calS$ 
are equivalent if there exists a sequence of 
transversally intersecting crescents from 
$\calR$ to $\calS$; that is, 
\[\calR = \calR_0, \calR_i \cap \calR_{i+1}^o \ne \emp, 
\calS = \calR_n \hbox{ for } i = 1, 2, \dots, n. \]

\begin{itemize}
\item We define $\Lambda({\calR})$ to be the union of 
all highest-level crescents equivalent to 
the highest-level crescent $\calR$. 
As before, $\Lambda({\calR})$ and $\Lambda({\calS})$ 
do not meet in the interior or they are the same. 
\item We define $\partial_I \Lambda({\calR})$
to be the boundary of $\Lambda({\calR})$ removed 
with the closure of the union of the $\alpha$-parts of 
the crescents in it. 
Then $\partial_I \Lambda({\calR})$ is 
a convex surface. 
\item We define $\partial_\alpha \Lambda({\calR})$ 
as the union of the $\alpha$-parts of the crescents 
equivalent to $\calR$. 
\end{itemize}

Recall that a pleated surface is a surface where 
through each point passes a geodesic. 

\begin{lem}\label{lem:pIlambda} 
The set 
\[\partial_I \Lambda({\calR}) \cap \tilde M - \tilde \Sigma\]
is a properly imbedded pleated surface.
\end{lem}
\begin{proof}
For each point of $x$ belonging to the above set, 
$x$ is an element of the interior of $\tilde M$ 
by Proposition \ref{prop:disjcrescent}. 
Let $B(x)$ be a small convex open 
ball with center at $x$. Then the 
crescents equivalent to $R$ meet $B(x)$ 
in half-spaces. Therefore the complement of 
their union is a convex subset of $B(x)$
and $x$ is a boundary point. There is a supporting 
half-space $H$ in $x$ and $H$ belongs to 
$\Lambda(\calR)$. 

If there were no straight geodesic 
passing through $x$ in the boundary 
set $\partial_I \Lambda(\calR)$, then
there exists a totally geodesic disk $D$ in $B(x)$  
with $\partial D$ in $\Lambda(\calR)$ 
but interior points are not in it. 

Since $\partial D$ is in $\Lambda(\calR)$, 
each point of $\partial D$ is in some crescent. 
We can extend $D$ to a maximal totally geodesic hypersurface 
and we see that a portion of the hypersurface
bounds a crescent $\calT$ containing $D$ in its $I$-part 
and overlapping with the other crescents.
Thus $\calT$ is a subset of $\Lambda(\calR)$
and so is $D$. 

Therefore, $\partial_I \Lambda(\calR)- \tilde \Sigma$ 
is a pleated surface. 
\end{proof}

\subsubsection{Outer-contact points of the crescent-sets}
\label{subsubsec:outerc}
%%%% June 24, 2004 11:40

We may have some so-called outer-contact points of $\tilde \Sigma$ at
$\partial_I \Lambda({\calR})$, i.e, those points with neighborhoods 
in $\tilde \Sigma$ outside $\Lambda(\calR)$. 
We can classify outer-contact points.

\begin{prop}\label{prop:outerpts}
Assume $\tilde \Sigma$ is well-triangulated.
The set of outer-contact points on $\partial_I \Lambda(\calR)$ 
is a union of the following components: 
\begin{itemize}
\item isolated points, 
\item an arc passing through the pleating locus
with at least one vertex. 
\item isolated triangles, 
\item union of triangles meeting each other 
at vertices or edges. 
\end{itemize}
\end{prop}
\begin{proof} 
This essentially follows by Proposition \ref{prop:sigmaI}.
\end{proof}

\subsubsection{Crescent isotopy itself}
\label{subsec:crescentiso}

Recall that the final resulting $\tilde \Sigma$ in State \ref{app:smtr} 
is a well-triangulated surface 
whose highest level outer crescents do not have any outer contact points.

First,we show that
\begin{equation}\label{eqn:replace}
\partial_I \Lambda({\calR}) 
\cup (\tilde \Sigma - \partial_\alpha \Lambda({\calR})) 
\end{equation} 
is a properly imbedded pleated surface. 

We do this for $\Lambda({\calR})$ for each highest-level crescent 
$\calR$ obtaining as the end result a properly imbedded surface 
$\tilde \Sigma'$. The deck transformation group 
acts on $\tilde \Sigma'$ since it acts on the union of 
$\Lambda({\calR})$. Thus, we obtain a new closed surface 
$\Sigma'$. 

Since the union of $\Lambda({\calR})$ for every highest-level crescent 
$\calR$ is of bounded distance from $\tilde \Sigma$ 
by the boundedness and the fact that $M$ is locally-compact, 
$\Sigma'$ is a compact surface.

\subsubsection{The isotopy}\label{subsubsection:isotopy}
%****** 2.1.2.3.2.2. Prove isotopy
We show that $\Sigma$ and $\Sigma'$ are isotopic. 

Let $N$ be the $\eps$-neighborhood of $\tilde \Sigma'$ in 
the closure of the outer component of $\tilde M - \tilde \Sigma$. 
There exists a boundary component $\partial_1 N$ nearer 
to $\tilde \Sigma$ than the other boundary component. 
The closure of a component $K$ of 
$\tilde M - \tilde \Sigma - \tilde \Sigma'$
contains $\partial_1 N$. 
Then $K$ projects to a compact subset of $M$. 
We can find a finite collection of generic secondary highest-level 
compact crescents 
$\calR_1, \dots, \calR_n$ and whose images under 
$\Gamma$ form a locally finite cover of $K$. 

We label the crescents by $\calS_1, \calS_2, \dots$. 
We know that replacing the closure of the $\alpha$-part of 
$\calS_1$ by the $I$-part is an isotopy. 
After this move, $\calS_2, \calS_3, \dots$
become new generic highest-level crescents 
by Proposition \ref{prop:transvint} and appropriate 
truncations. 

We define $\partial_I (\calS_1 \cup \calS_2 \cup \dots)$ 
as the boundary of $\calS_1 \cup \calS_2 \cup \dots$ 
removed with the union of the $\alpha$-parts of $\calS_1, \calS_2, \dots$. 
Again, this is a convex imbedded surface. 
Therefore, replacing the union of 
the $\alpha$-parts of $\calS_1, \calS_2, \dots$ 
by $\partial_I (\calS_1 \cup \calS_2 \cup \dots)$ is 
an isotopy as above. 

We obtain $\Sigma_{\calR_1, \dots, \calR_n}$ as 
the image in $M$, which is isotopic to $\Sigma$.
If $\eps$ is sufficiently small, then 
we see easily that 
$\partial_I (\calS_1 \cup \calS_2 \cup \dots)$ in $N$ 
can be isotopied to $\tilde \Sigma'$ 
using rays perpendicular to $\tilde \Sigma'$. 
Thus, $\Sigma_{\calR_1, \dots, \calR_n}$
is isotopic to $\Sigma'$.
We showed that $\Sigma'$ is isotopic to $\Sigma$. 
\qed

\subsubsection{The strict decreasing of the levels}

We show that our isotopy move decreases the level strictly. 

\begin{prop}\label{prop:nocrescent} 
If $\Sigma'$ is obtained from $\Sigma$ by 
the highest-level outer crescent move for level $n$ for an integer $n \geq 0$, 
then the union of the collection of crescents of $\tilde \Sigma$ 
contains the union of those of $\tilde \Sigma'$, and 
$\tilde \Sigma'$ has no outer highest-level crescent of level $n$ or higher
for an integer $n \geq 0$.
The level of the highest level outer crescents of 
$\tilde \Sigma'$ is now $\leq n-1$.
The level for inner highest-level crescents 
stays the same or may drop.
The same statements hold if we replace the word ``outer'' by 
``inner''.
\end{prop}
\begin{proof}
The outer highest-level crescents for $\tilde \Sigma'$ 
can be extended to ones for $\tilde \Sigma$ since 
their $I$-part can be extended across. 

The inner highest level 
crescents for $\tilde \Sigma'$ can be 
truncated to ones for $\tilde \Sigma$ by Lemma \ref{lem:crescentiso}
since the move from $\tilde \Sigma'$ to $\tilde \Sigma$ 
is inward and supported by the outer crescents of $\tilde \Sigma$.
Thus the first statement holds. 

If there were outer highest level 
crescent $\calR$ of level $n$ or higher, 
then we can extend $I_\calR$
across the $I_\partial$-parts so that we can obtain 
a level-$n$ or higher-level crescent. Such a crescent would have 
been included in the crescent set and should have been 
isotopied away.

If $\calR$ were inner highest level one, 
Lemma \ref{lem:crescentiso} implies the result.
\end{proof}

We apply our methods of Subsection \ref{subsec:crescentiso} to obtain $\Sigma'$ 
which we now let it be $\Sigma$. 
\begin{app}\label{app:criso}
$\tilde \Sigma$ now has no outer crescents of level $n$.
However, it may not be triangulated. 
\end{app} 

%* 3. Convex perturbation theory
\subsection{Convex perturbations} \label{subsec:convpert}

In this subsection, 
we modify $\tilde \Sigma$, which now has partial pleating, 
obtained above further. We discuss how a surface with a portion of 
itself concavely pleated by infinitely long geodesics and the remainder 
triangulated can be perturbed to a triangulated surface 
without introducing higher-level crescents. 
This is done by approximating the union of pleating geodesics
by train tracks and choosing normals in the concave direction 
and finitely many vertices at the normal and pushing 
the pleating geodesics to become geodesics broken at 
the vertices.

%%%% December 18 5:54
%** 3.1. Definition
\subsubsection{Pleated-triangulated surfaces}

We will describe $\tilde \Sigma$ as a ``pleated-triangulated" surface. 
We will then perturb the crescent-isotopied $\Sigma$ into 
a triangulated surface not 
introducing any higher-level crescents. 

Suppose that $\Sigma'$ is a closed imbedded surface in $M$.
$\Sigma'$ is a {\em pleated-triangulated surface} 
if 
\begin{itemize}
\item $\tilde \Sigma'$ contains a closed $2$-dimensional domain 
divided into locally finite collection
of closed totally geodesic convex
domains meeting each other in geodesic segments,
\item through each point of the complement in $\tilde \Sigma'$ passes 
a {\em straight} geodesic in the complement.
\end{itemize} 
We may also add finitely many straight geodesic segments 
in the surface ending at the domain.
\begin{itemize}
\item The domain union with the segments 
is said to be the {\em triangulated part} of $\Sigma'$. 
The boundary of the domain is a union of finitely pinched 
simple closed curves. 
\item The complement of the domain is an open surface, which is
said to be the {\em pleated part} where through each point 
passes a straight geodesic.
\item The pleated part has a locus where a unique 
straight geodesic passes through. This part is said to be the 
{\em pleating locus}. It is a closed subset of the complement and 
forms a lamination. 
\item If we remove the closure of the pleated part from 
$\tilde \Sigma'$, we obtain a locally finite collection of 
totally geodesic convex domains meeting each other in edges and vertices. 
The convex domains, edges, and vertices are in general position.
\end{itemize}

For later purposes, we say that $\Sigma'$ is 
{\em truly pleated-triangulated} if 
the triangulated part is a union of totally
geodesic domains that are convex polygons (i.e., finite-sided) and geodesic 
segments ending in the domains.

While the triangulated parts and pleating parts are 
not uniquely determined, we simply choose what seems 
natural. We also assume that the pleated part is locally convex or 
locally concave. We usually choose a normal direction 
so that the surface is locally concave at the pleated part. 

If $\Sigma'$ satisfies all of the above conditions, 
we say that $\Sigma'$ is a {\em concave pleated-triangulated} 
surface. If we choose the opposite normal-direction, 
then $\Sigma'$ is a {\em convex pleated, triangulated} 
surface.

%We show that the surface obtained in the previous 
%section is concave pleated-triangulated. 

\begin{prop}\label{prop:ptsurface}
Suppose that $\Sigma$ is as in State \ref{app:criso}
where $n$, $n \geq 0$, is the highest-level for the surface before the 
applications of methods in Subsection \ref{subsec:crescentiso} 
which is realized by outer-crescents. 
Then $\Sigma$ is a concave pleated-triangulated 
surface in the outer direction.
The level of highest-level outer crescents of 
$\tilde \Sigma$ is $\leq n-1$ and the level of 
highest-level inner crescent of $\tilde \Sigma$ is 
less than or equal to $n$.
Moreover, the statements are true if all ``outer" were 
replaced by ``inner" and vise-versa 
and the word ``concave'' by ``convex''.
\end{prop}
\begin{proof} 
The part $\partial_I \Lambda(\calR) - \tilde \Sigma$ 
for crescents $\calR$s are pleated by Lemma \ref{lem:pIlambda}. 
These sets for crescents $\calR$ are either identical or 
disjoint from each other as the sets of form 
$\Lambda(\calR)$ satisfy this property. 
The union of sets of form $\partial_I \Lambda(\calR)$ 
comprise the pleated part and the complement in 
$\Sigma'$ were in $\Sigma$ originally and 
they are the union of totally geodesic 
$2$-dimensional convex domains. 

The rest is proved already in Proposition \ref{prop:nocrescent}. 
\end{proof}

\subsubsection{The geometry and topology of pleating loci}

Two leaves in a pleating locus are {\em converging} if one 
is asymptotic to the other one 
(see Section \ref{sec:metricspaces}
for definitions);
i.e., the distance function from one leaf to the other 
converges to zero and conversely.
By an end of a leaf of a lamination {\em wrapping around} 
a closed set, we mean that the leaf converges to a subset of 
the closed set in the direction of the end. 

The main idea of classifying the pleating locus 
are from those of Thurston as 
written up in Casson-Bleiler \cite{CB}.

%%%% Oct 5 11:50

\begin{lem}\label{lem:nowrapping} 
Suppose that $\Sigma'$ is a closed concave pleated-triangulated 
surface with a triangulated part and pleated part assigned. 
Suppose that $l$ is a leaf.
Then 
\begin{itemize}
\item $l$ is either a leaf of a minimal geodesic lamination or 
a closed geodesic, or each end of $l$ wraps around  
a minimal geodesic lamination or a closed geodesic
or ends in the boundary of the triangulated part.
\item if $l$ is isolated from both sides, then 
$l$ must be a closed geodesic, and 
\item pleating leaves in a neighborhood of $l$ 
must diverge from $l$ eventually.
\end{itemize}
\end{lem} 
\begin{proof}
Since the pleated open surface carries an intrinsic metric 
which identifies it to a quotient of 
an open subset of the hyperbolic space, 
each geodesic in the pleating lamination will 
satisfy the above properties like the geodesic laminations 
on the closed hyperbolic surfaces. The first item is done. 

For second item, 
if $l$ is isolated from both sides, then there is a definite positive 
angle between two totally geodesic hypersurfaces ending at $l$.
Suppose that $l$ is not a simple closed geodesic. 
Then this angled pair of the hypersurfaces continues to wrap around
infinitely in $M$ accumulating at a point of $M$ and 
the sum of the angles violates the imbeddedness of $\Sigma'$.
(This is an argument in Thurston \cite{Thnote} to 
show the similar argument for the boundary of 
the convex hulls.)

For the third item,
if $l$ is not isolated but has converging 
nearby pleating leaves, the same reasoning will hold as in 
the second item.
\end{proof}

\begin{prop}\label{prop:pleatinglocus} 
Suppose that $\Sigma'$ is a closed concave pleated-triangulated 
surface. Let $\Lambda$ be the set of pleating locus of the pleated
part in $\Sigma'$. Then $\Lambda$ decomposes into finitely many 
components $\Lambda_1, \dots, \Lambda_n$ so 
that each $\Lambda_i$ is one of the following:
\begin{itemize}
\item a finite union of finite-length pleating leaves homeomorphic to 
a compact set times a line with endpoints in the triangulated 
part. {\rm (} A discrete set times a line 
if $\Sigma'$ is truly pleated-triangulated. {\rm )}
\item a simple closed geodesic.
\item a minimal geodesic lamination, which is a closed 
subset of the pleated part isolated away from the triangulated 
part. 
\end{itemize} 
Here each leaf is either bi-infinite or finite.
The union of bi-infinite leaves is a finite union of 
minimal geodesic laminations and is isolated away 
from the triangulated part and the union of 
finite-length leaves. 
\end{prop}
\begin{proof}
Let $l$ be an infinite leaf in the pleating locus. 
By Lemma \ref{lem:nowrapping}, $l$ is not 
isolated from both sides and the leaves in its neighborhood 
is diverging from $l$. 
If $l$ is not itself a leaf of a minimal lamination, 
then an end of $l$ must converge to a minimal lamination or 
a simple closed geodesic. This means that leaves in 
a neighborhood also converges to the same lamination
in one of the directions. However, this means that they 
also converges to $l$, a contradiction. 
Therefore, each leaf is a leaf of a minimal lamination or 
a closed geodesic or a finite length line. 

The union of all finite length lines in the pleating locus 
is a closed subset: Its complement in $\Lambda$ is 
a compact geodesic lamination in $\Sigma'$.
If a sequence of a finite length leaves $l_i$ converges to 
an infinite length geodesic $l$, then $l_i$ gets 
arbitrarily close to a minimal lamination or
a closed geodesic. If $l_i$ gets into a sufficiently
thin neighborhood of one of these, then
a neighborhood of an end of $l_i$ must be in a sufficiently 
thin neighborhood of one of these by
the imbeddedness property of $l_i$, i.e., 
cannot turn sharply away and go out of the neighborhood.
As $l_i$ ends in the triangulated part, 
the distance from the triangulated part to 
one of these goes to zero. 
Since the domains in the triangulated part 
are in general position, the boundary of
the triangulated part cannot contain 
a closed geodesic or the straight geodesic lamination. 
This is a contradiction. 

Looking at an $\eps$-thin neighborhood of $\Lambda$, 
we see that $\Lambda$ decomposes as described.
(See Casson-Bleiler \cite{CB} for background informations).
\end{proof}

\begin{rem}\label{rem:pleatinglocus} If $\Sigma'$ is truely pleated-triangulated, 
then there are only finitely many finite length 
pleating leaves since their endpoints are on 
the vertices of the polygons in the triangulated part.
\end{rem}

%%%% Jan 11 5:10 
\begin{lem}\label{lem:pleatinglocus} 
A line $l$ in the pleating locus of a pleated-triangulated surface cannot 
end in an interior of a segment $s$ in the boundary of the pleated part.
\end{lem} 
\begin{proof} 
First suppose $l$ is isolated. Then the boundaries of totally geodesic planes in its side must 
contain open segments in $s$, and the planes have to be identical, contradicting that $l$ is 
in the pleating locus. 

Suppose that $l$ is not isolated. Then the nearby leaves of the foliation must end at the same place as 
$l$; otherwise, we get that the nearby totally geodesic planes are identical. If they all 
end at the same place, again a similar reasoning shows that the nearby geodesic planes 
are identical. These contradicts the fact that $l$ is in the pleating locus.
\end{proof}

%*** 3.1.2. Convex vertex move first level n

%** 3.2. Convex perturbations

%*** 3.2.1 Define straight train tracks, laminations, in a limited sense
%%tie, switch, branch, 
%%N \tau, 

\subsubsection{The perturbation moves.}
Recall that our surface $\tilde \Sigma$ is pleated triangulated and 
has level $n$, achieved by innermost inner-crescents, for an integer $n \geq 0$.
The outer level is less than or equal to $n-1$.

A train track is obtained by taking a thin neighborhood of 
the lamination. We can think of the train track as a union of 
segment times an interval, so-called branches, joined up at the end of 
each segment times the intervals so that the intervals 
stacks up and matches. A point times the interval 
is said to be a {\em tie} and a tie where 
more than one branches meet a {\em switch}.
One can collapse the interval direction to obtain 
a union of graphs and circles. For more details, 
consult Casson-Bleiler \cite{CB}. 

We will perturb the pleating locus to obtain our 
well-triangulated 
surface. We describe our results divided 
in Theorem \ref{thm:convexperturb}.

\begin{description}
\item[The first step (I) of the perturbation move]
Let $l_1, \dots, l_k$ be the thin strips 
containing all the finite length open leaves ending at 
the triangulated parts. 
We find a thin totally geodesic hypersurface 
$P_i$ near $l_i$s nearly parallel to $l_is$. 
Then we cut off the neighborhood of $l_i$ in $\Sigma$ 
by $P_i$ and replace the lost part with the portion in
$P_i$. This introduces squares which are triangulated into
pairs of triangles. 

We now remove the union of the squares 
from the pleated part and add the union to 
the triangulated part. Now we retriangulate the 
triangulated part. They will be immobile during the pertubations now.

This forms a generalization of a small truncation move. 
We still call it a {\em small truncation move}. --(*)

As in Proposition \ref{prop:perturb}, we don't increase the level
and we remove some of the outer-contact points. 

Finally, we do remove all of the outer-contact points in the triangulated 
part, i.e., the closed set consisting of the compact totally geodesic domains. 

\item[The step (II) of the perturbation move] 
Now, every component of the pleating locus is 
contained in a minimal lamination. The set of pleated 
locus is a union of finitely many minimal laminations. 

\begin{defn}\label{den:minpleat}
In the pleated part, we define the {\em minimal pleated part} 
as the intrinsic convex hull in the closure of the 
pleated part of the union of bi-infinite pleating 
geodesics with respect to the intrinsic metric  
obtained by piecing the pleated parts together. 
The minimal pleated part is a subsurface of 
the pleated part which is the closure of the union of 
totally geodesic subsurfaces bounded by 
the bi-infinite pleating geodesics. 
\end{defn}

The boundary of the minimal pleated parts is 
a union of simple closed curve which might be  
a broken geodesic or just a geodesic.
We triangulate the closure of the complement of
the minimal pleated part, which will be added shortly 
to the triangulated part. This may introduce 
vertices at the boundary of
the minimal pleated part. 

We now add finite leaves of infinite length in 
the minimally pleated part so that the components of the 
complement of the union of the pleating locus and 
these leaves are all open triangles. These can start 
from the vertices at the boundary of the minimal pleated part.
This can be done even though the boundary of
the pleated part is not geodesic. 

By choosing sufficiently small $\eps$-thin-neighborhoods
of the union for $\eps>0$,  we obtain a train track. 
We obtain a maximal train track. 

We first choose switches for the endpoint of the finite 
length geodesics in the squares and added infinite length
finite leaves. (The switches are transverse arcs.)
We choose switches for the rest. 
We label them $I_1, \dots, I_m$. We may have a uniform bound on 
$m$ depending only on the Euler characteristic of 
the open pleated part surface and 
a uniform lower bound to the distances between any two of $I_i$. 
(Note here $m$ is bounded above by a constant depending only 
on $\Sigma$ as $\Sigma$ has \cat$(-1)$ metric and Thurston's theory 
of geodesic laminations hold for such surfaces)

By choosing $\eps >0$ sufficiently small, we can assume that 
the outer-normal vectors to totally geodesic hypersurfaces meeting
$I_i$ are $\delta$-close for a small $\delta > 0$
except the outer-normal vectors to the totally geodesic hypersurfaces 
corresponding to the complementary regions of the train track.
(Here the outer-normal is in the concave directions.)

The union of $I_i$s with the leaves of the train tracks 
divides the pleated part into infinitely squares with two sides 
within $I_i$s and polygons with some edges in $I_i$s. 
We regard the endpoints of $I_i$ in edges of some polygons 
to be a vertex and we retriangulated these accordingly.
We can assume that there is a lower bound to the length of 
each edges which are not in $I_i$s. In fact, we can 
assume that the ratios of the edges not in $I_i$ to those in 
$I_i$ are greater than $1000$. 

Topologically, the train track collapses to a union of 
graphs with vertices corresponding to $I_1, \dots, I_n$ 
and closed geodesics, and the complement 
becomes a union of totally geodesic triangles.

We now do the collapsing geometrically:
We choose one of the normal vectors and 
a generically chosen point $x_i$ on the normal vector
$\gamma$-close to $I_i$, where $\gamma$ is a small positive number.
We push all the points of $I_i$ to $x_i$ 
to obtain a train track $\tau_{\eps, \delta, \gamma}$ 
and the complementary regions 
move accordingly to disks divided into compact 
totally geodesic geodesic triangles with vertices in 
the train track $\tau_{\eps, \delta, \gamma}$ 
and in triangulated part.
($x_i$s are said to be {\em pleated part vertices}.) 

We claim that then the triangles are very close to the original 
triangles in their normal directions as well as in the
Hausdorff distance since the edge lengths of the triangles 
are bounded below.
Since there are no rapid turning 
of the complementary geodesic triangles, 
we can be assured that the new surface is imbedded
by integrations. 
We can see this by looking at a cross-section, 
which is a function of bounded variations.

The leaves of the laminations are moved to become a train track in 
the normal direction which is a concave direction.
Thus the leaf is bent against the concave direction.
The triangles meeting the endpoints of $I_i$ 
of almost the same direction as before and hence have angles $< \pi$
by concavity.
Lemma \ref{lem:cross} shows that the pleated part vertices 
are strict saddle vertices. 

\begin{figure}[h]
\centerline{\epsfxsize=3.5in \epsfbox{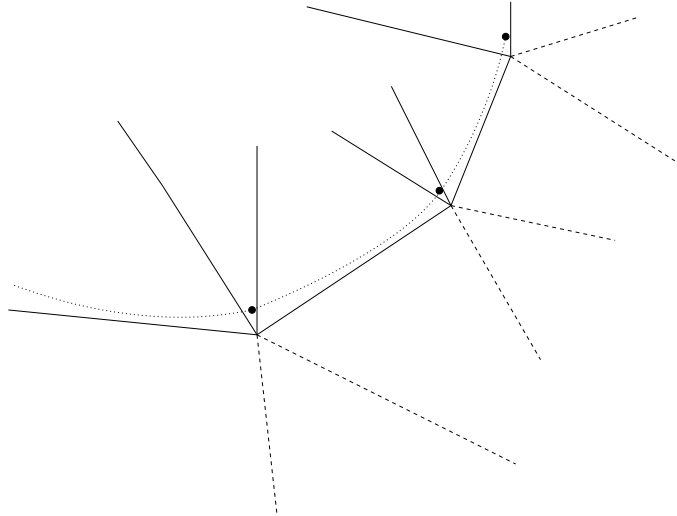}}
\caption{\label{fig:svertices} Making the boundary 
vertices of the pleated part into saddle-vertices.}
\end{figure}
\typeout{<>}

Of course, here there are many choices and the combinatorial 
structure of the triangulations will change according to the choices.

There is a slight modification when $n=0$: 
Since we started from well-triagulated $\tilde \Sigma$. 
All vertices in the interior of the triangulated parts are strict saddle vertices 
or convex vertices and the vertices in the boundary of the triangulated part 
are saddle vertices as there are no outer crescents by Corollary \ref{cor:cressv}.
By Lemma \ref{lem:deforming}, 
we perturb the vertices in the boundary of the triangulated part to be 
strict saddle vertices. The adjacent component of the pleated part 
with the pleating locus removed is totally geodesic. 
By introducing finite geodesic leaves, we see that 
the pleating locus will become larger including these leaves. 
By taking a sufficiently small perturbation at the boundary vertices, 
this will not change the vertex type of all vertices in 
the triangulated part. Also, we choose a train track for 
the pleating minimal lamination. The angles between the adjacent totally 
geodesic planes across the switches square are less than $\pi$. 
We further require that 
we perturb the boundary vertices by sufficiently small amount 
so that the angles are still less than $\pi$.
Finally, we perturb the pleated parts as 
above introducing pleated part vertices.

(Clearly, we can do the above collapsing 
and the perturbations equivariantly).
\end{description}

\begin{thm}\label{thm:convexperturb} 
Let $\Sigma'$ be a closed concave pleated, triangulated surface
where the outer highest level crescents have level $n-1$ for 
an integer $n$, $n \geq 0$.
Then one can find an imbedded isotopic 
well-triangulated surface $\Sigma''$ in 
any $\eps$-neighborhood of $\Sigma'$ by the above methods of 
perturbations so that
the following hold. 
\begin{itemize}
\item[(i)] The union of the set of crescents for $\tilde \Sigma''$ 
is in the $\eps$-neighborhood of that of $\tilde \Sigma'$
and vice-versa for a small $\eps$ 
if we choose $\tilde \Sigma''$ sufficiently close to $\tilde \Sigma'$.  
\item[(ii)] The pleated part vertices are strict saddle-vertices.
\item[(iii)] The level of the outer highest-level crescents
for the resulting surface 
$\tilde \Sigma''$ is less than or equal to $n-1$ 
and the level of inner highest level crescents for
$\tilde \Sigma''$ is less than or 
equal to that of $\tilde \Sigma'$.  
\item[(vi)] In particular, if $\tilde \Sigma'$ has 
no outer {\rm (}resp. inner{\rm )} crescent, 
then $\tilde \Sigma''$ contains no outer {\rm (}resp. inner{\rm )} crescent.
If there were no outer and inner crescents, $\tilde \Sigma''$ is 
saddle-shaped. 
\end{itemize}
\end{thm}
\begin{proof}
(i) This matters for crescents that are inner if 
the perturbations are inner and ones that are outer if 
the perturbations are outer: In other cases, 
Lemma \ref{lem:crescentiso} shows that reversing 
the perturbation process gives us back all of our old
crescents preserving the $I$-part hypersurface.

Suppose that we have the perturbed sequence 
$\Sigma'_i$ closer and closer to $\Sigma'$, 
and there exists a sequence of crescent $\calR_i$ 
for $\Sigma'_i$ not contained in
a certain neighborhood of the union of 
crescents for $\Sigma'$. Then the limit $\calR$ of $\calR_i$
is still a crescent for $\Sigma'$ and is not in 
the neighborhood. This is absurd. 

(ii) This is proved above.

(iii) Let $n-1$ be the level of outer highest-level
crescents of $\Sigma'$. 

As stated earlier, the perturbation step (I)
does not increase the level.

Let us discuss what are the outer-contact points of 
outer highest-level crescents, where we are no longer assuming 
the general position property of $\tilde \Sigma'$:  

Let $\calR$ be a highest-level crescent with outer-contact 
points. There could be points of the minimal sets pleated parts 
meeting tangentially the $I$-part of $\calR$.
If not, then they are on the triangulated part. 
In these cases, we do the small truncation move.
This will not increase the folding number. 

The minimal sets in pleated part may meet the crescents 
only tangentially at the $I$-part or meet the closures of 
the $\alpha$-parts of the crescents: If not, then 
the minimal pleated part must pass through the interior of 
a crescent and this implies that a bi-infinite pleating 
geodesic pass through the interior.
By Proposition \ref{prop:avoidg}, 
this is a contradiction. 

If a minimal set in the pleated part meets the $I$-part of a crescent, 
then the closure of a certain number of pleated leaves 
are in the $I$-part. Thus, it is clear that the complete 
geodesic leaf or the closure of a complement of the pleating locus 
in the minimally pleated part equals 
the set of outer-contact points of the $I$-part. 

We now move vertices of the train tracks of 
the pleating laminations by a very small amount 
according to the perturbation step (II). 
(The crescents do move in its $\eps$-neighborhood.) 

We choose the $\eps$-neighborhood sufficiently small 
so that any new component of $\tilde \Sigma'$ intersected 
with crescents may not arise as $\tilde \Sigma'$ is deformed: 
The small truncated places may be avoided 
by taking $\eps$-sufficiently small. 

Since any crescent during the perturbation cannot
meet the minimal pleated part in its interior and its $\alpha$-part 
by Proposition \ref{prop:avoidg}, 
the $I$-part of crescents and crescents themselves
close to the minimal pleated part
lie below the minimal pleated part or may meet the minimal
pleated part but cannot pass through it.

Let us choose a sequence $\eps_i \ra 0$ of positive real numbers 
$\eps_i$ and a sequence of approximating isotopic surfaces 
$\Sigma^{\eps_i} \ra \Sigma'$.  We assume that they all have 
the same number of vertices by introducing finer triangulations if necessary.
We can complete the sequence $\Sigma^{\eps_i}$ 
to a one-parameter family since any two triangulation with same number of
vertices on a closed surfaces are related by elementary moves. 
(Of course, there is a skipping around a tetrahedron during 
the elementary move. See Figure \ref{fig:tet})

\begin{figure}[h]
\centerline{\epsfxsize=2.5in \epsfbox{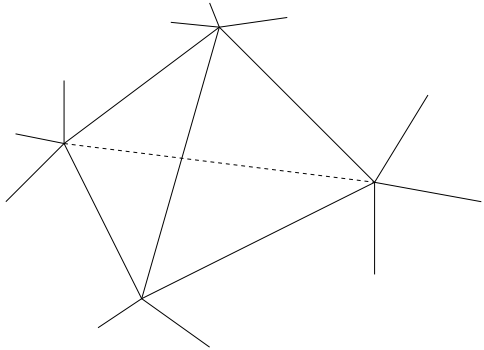}}
\caption{\label{fig:tet} 
An elementary move and a skipping a tetrahedron.
}
\end{figure}
\typeout{<>}

Suppose that the outer-level became $n$ or higher for $\Sigma^{\eps_i}$
for every $i$. Choose an index $i$.
Then there is a vertex $v$ of $\tilde \Sigma^{\eps_i}$ in the interior of 
a crescent $R$ of $\tilde \Sigma^{\eps_i}$ of level $n$ so that a component $C$
of $\tilde \Sigma^{\eps_i} \cap R^o$ containing $v$ is in the $\alpha$-part of
an innermost crescent realizing the level $n$. 

If $v$ is from the triangulated part and $n \geq 1$, then by letting $i \ra \infty$, 
we obtain an outer-contact point as in Proposition \ref{prop:perturbation}. 
This was ruled out by the perturbation step (I). 

Suppose that $n=0$.
Then all vertices in the interior of the triangulated parts are strict saddle vertices 
or convex vertices and the vertices in the boundary of the triangulated part 
are saddle vertices after the perturbation step (II) by construction. 

The vertex $v$ cannot come from the triangulated part. 

Since the vertex $v$ is from the pleated part of $\tilde \Sigma'$, 
it follows that $v$ is a saddle-type vertex. 
Moreover each vertex of $C$ are saddle-type by the same reason. 
Corollary \ref{cor:cressv} gives us a contradiction.

Since our move is toward outside only, we see that 
the set of inner crescents decrease only and hence 
the inner level do not increase. 

%If $\tilde \Sigma'$ has no inner crescents, 
%then we construct $\tilde \Sigma''$ with no inner crescent
%as above. Also, the outer-level do not increase as above. 

(iv) The first part follows from (iii).

If $\tilde \Sigma'$ has no inner or outer crescents, 
the first part of (iv) shows that 
we can construct $\tilde \Sigma''$ with no crescents.
This implies that $\tilde \Sigma''$ is saddle-imbedded.

\end{proof}

We now do small truncation moves to $\tilde \Sigma$ so that outer 
highest-level crescent of level now $n-1$ has no outer-contact 
points. 
We now perturb the triangles so that resulting $\tilde \Sigma''$ 
is well-triangulated. 
As in Application \ref{app:smtr}, we see that 
the outer level of $\tilde \Sigma''$ is still $n-1$.

\begin{app}\label{app:convpert} 
Applying the method of this section, we obtain a well-triangulated $\tilde \Sigma$ 
with outer level $n-1$ and the inner level $\leq n$. 
\end{app}

%%% Feb 18 11:30

\subsection{The proof of Theorem C and Corollary D}\label{subsec:thmC}

\subsubsection{Isotopy sequences}\label{subsec:isotopy}
%***** Convex perturbation theory....
We review the outer level $n$ crescent moves. We temporarily denote 
the result of the move after each move ---(*):
\begin{itemize}
\item[(i)] the small truncation moves for outer highest-level crescents: 
$\tilde \Sigma^{(i)}$.
\item[(ii)] the crescent-isotopy for outer highest-level crescents: 
$\tilde \Sigma^{(ii)}$.
\item[(iii)] the convex perturbations which involves 
small truncation moves: $\tilde \Sigma^{(iii)}$.
\end{itemize}

Let $n$ be the highest level for $\tilde \Sigma$.
First, we do the above outer highest-level crescent moves of level $n$.

We should worry about one issue here
with regards to $\cl_1, \dots \cl_n$:
After the step (ii), we may find $\tilde \Sigma^{(ii)}$ meeting a lift of 
one of these. In which case, a lift $l_i$ for some $i$ is contained in the $I$-part of a secondary 
highest level crescent. 
In this case, $\cl_{i_1}, \dots, \cl_{i_k}$ 
is a subset of the image surface $\Sigma^{(ii)}$ for some integer $k \geq 1$. 

Choose an arbitrary $j$. 
The deck transformation acting on $l_{i_j}$ preserves the $I$-part of the secondary highest level crescent
since otherwise we get more secondary highest level crescent around $l_{i_j}$ completing 
the $2\pi$-angle and hence $l_{i_j}$ can't be in $\tilde \Sigma^{(ii)}$. 
Thus the holonomy of $\cl_{i_j}$ is hyperbolic, i.e., loxodromic without rotational part.

Let $P_1, \dots, P_m$ be the $1$-complexes which are components of the union of 
$\cl_{i_1}, \dots, \cl_{i_k}$, which could just be the disjoint union of the closed 
curves themselves.

%%% 11:58 Feb 20

After the convex perturbation (iii), we find the $1$-complex $P'_i$
in $\Sigma^{(iii)}$ approximating $P_i$:
What happens actually is that in case $P_i$ is a simple closed geodesic, then 
$P'_i$ is perturbed to have saddle vertices only, and in case $P_i$ is a true 
$1$-complex with singularities, then there is a totally geodesic subsurface $S_i$ with boundary, 
and the boundary vertices are perturbed to become saddle vertices and 
$S_i$ is isotoped accordingly to $S'_i$ without introducing any interior vertices and
$P_i$ is isotopied to $P'_i$ on $S'_i$ accordingly.
Since the holonomies of closed curves on $S_i$ are without rotational parts, 
we assume without loss of generality that we moved $P_i$ to $P'_i$ in 
a fixed angle to the $\Sigma^{(ii)}$ (perpendicularly in the second case). 
Let $\tilde P_i$ be the inverse image of $P_i$ in $\tilde \Sigma^{(ii)}$, 
and $\tilde P'_i$ be that of $P'_i$. 
Thus $\tilde P_i$ and $\tilde P'_i$ lie on a union of 
totally geodesic planes $L_i$ in a fixed angle to the $I$-part of a crescent for 
$\tilde \Sigma^{(ii)}$, where a $2$-complex $D_i$ is bounded by $\tilde P_i$ and $\tilde P'_i$. 
($D_i$ maps to a compact $2$-complex in $M$.)

Recall that our isotopies are in the outer-direction and the convex perturbation
was also in the outer-direction.
\begin{lem}\label{lem:pertgeo} 
\begin{itemize}
\item If $P'_i$ is a closed geodesic, then a highest level 
crescent $\calR$ of $\tilde \Sigma^{(iii)}$ 
may meet a lift of $P'_i$ in $\tilde M$ only in an $\alpha$-part by a connected subarc. 
\item If $P'_i$ is a true $1$-complex with singularities, then a highest level crescent $\calR$ 
of $\tilde \Sigma^{(iii)}$ may meet 
a component of the inverse image of $S'_i$ in $\tilde M$ 
in a connected surface. 
\item In both case, when meeting, the $I$-part of a crescent for the isotopied 
$\tilde \Sigma$ passes $D_i$ and hence the $I$-part is in the inner direction of 
the isotopied $\tilde \Sigma$.
\end{itemize}
\end{lem}
\begin{proof} 
If $P'_i$ is a simple closed curve, then $P'_i$ is an arc bent in the inner-direction. 
Also, $S'_i$ is bent in the inner-direction since the deck transformations 
corresponding to $S'_i$ acts on a totally geodesic hypersurface $H_i$ containing $l_{i_j}$s
and hence is a Fuchsian subgroup.
The result follows from this. 
\end{proof}

%%%% JAN 07 2:30
Finally, we let $\Sigma = \Sigma^{(iii)}$.
During the isotopy, $\Sigma$ may pass through these $\cl_{i_j}$s. 
Thus, $\Sigma$ may no longer be incompressible. 

The next step is that we do the inner highest-level crescent moves of level $n$. 

Here, the incompressibility may not hold: however, 
the only place we need incompressibility is when we show that 
there is a secondary highest-level crescent containing a given highest-level crescent:
We need to  show that each boundary curve $c$ of the perturbed $I$-part of a highest-level crescent
which is an innermost component in the perturbed $I$-part bounds a disk in $\tilde \Sigma$. 
If $c$ does not meet $\tilde P'_i$, then the argument is the same since 
$\tilde \Sigma - \tilde P'_i$ is incompressible. If $c$ does meet $\tilde P'_i$, then 
we reverse the convex perturbation move and the isotopied $c'$ from $c$ 
is a closed curve in $\tilde M$ and hence bounds a disk $D$ in 
$\tilde \Sigma^{(ii)}$ before the convex perturbations (regarding the 
$\cl_{i_j}$ to be pushed outward from $\tilde \Sigma^{(ii)}$ by a small amount).
$D$ meets the the corresponding inverse images of $P_i$ and $H_i$
in regions. Since $P_i$ and $H_i$ are isotopied by a very small amount, 
we can isotopy $D$ to $D'$. Thus, $c$ bounds a disk in the
convex perturbed $\tilde \Sigma^{(iii)}$. 

If $\tilde \Sigma^{(iii)}$ passed over $P_i$, then Lemma \ref{lem:pertgeo} shows that 
the inner level $n$ crescent move can only move $P'_i$ closer to $P_i$ since the 
direction of the crescent is now reversed but cannot pass $P_i$. Since the perturbation
can be controlled, we see that the result $\tilde \Sigma'$ of the inner level $n$ crescent 
move, there is $P''_i$ very close to $P_i$ also.

We go to the level $n-1$ and so on. 
We see that the inner and outer level strictly decreases
until there are no more crescents. 
Therefore, the final result is saddle-imbedded. 

Moreover the union of the set of crescents 
is contained in the $\eps$-neighborhood of 
the union of the set of crescents in the previous step. 

This completes the proof of Theorem C. \qed

We now prove Corollary D:
Assume as in the premise of the corollary that $M$ is a codimension $0$ 
compact submanifold of a general hyperbolic manifold $N$.

We modify the boundary surface $\partial M$ according to 
Theorem C. The final result is an imbedded closed surface $\Sigma$
bounding a compact region by our construction. 
We let this region to be our isotopied manifold $M'$. \qed

\begin{rem}\label{rem:svertices} 
The vertices of the boundary of $M'$ are 
strict saddle-vertices: We start from general position 
$\partial M$ so that each saddle-vertex is a strict saddle-vertex
(see Theorem \ref{thm:convexperturb}).
Newly created vertices in the interior of the pleated 
part after the crescent isotopies are all strict saddle-vertices.
The boundary vertices of the pleated part at level $0$ 
are perturbed to be strict saddle-vertices.
Finally, we no longer 
have any convex or concave vertex left at $M'$.
\end{rem}

%\begin{rem}\label{rem:lamination} 
%In fact, we can even assume that $\cl_i$ are 
%geodesic laminations since we can modify 
%Proposition \ref{prop:geodesicavoid}
%using Proposition \ref{prop:avoidg}.
%\end{rem}

%% Jan 08 12:20

\part{General hyperbolic $3$-manifolds and convex hulls of 
their cores}

\section{Introduction to Part 2}

%** Definitions and Main result

A {\em hyperbolic triangle} is a subset of a metric space isometric 
with a geodesic triangle in the hyperbolic plane $\hypp$.
If the ambient space is a $3$-dimensional metric space, 
then we require it to be totally geodesic as well 
and develop into a totally geodesic plane in the hyperbolic 
space $\hyps$. 

A {\em general hyperbolic manifold} $M$ is a metric space 
with a locally isometric immersion $\dev: \tilde M \ra \hyps$
from its universal cover $\tilde M$. 
$\dev$ has an associated homomorphism 
$h: \pi_1(M) \ra \Isom(\hyps)$ given 
by $\dev \circ \gamma = h(\gamma)\circ \dev$ 
for each deck transformation $\gamma \in \pi_1(M)$. 
We require that the boundary of $M$ 
triangulated by totally geodesic hyperbolic triangles. 

A {\em $2$-convex} general hyperbolic manifold is 
a general hyperbolic manifold
so that a isometric imbedding from $T^o \cup F_1 \cup F_2 \cup F_3$ 
where $T^o$ is a hyperbolic tetrahedron in the standard 
hyperbolic space $\hyps$  and $F_1, F_2, F_3$ its faces extend to 
an imbedding from $T$.

The main result of Part 2 is:
\begin{thm}\label{thm:Gromovh}
The universal cover of 
a $2$-convex general hyperbolic manifold is \cat$(-1)$. 
\end{thm}

A {\em hyperbolic-map} from a triangulated hyperbolic surface is 
a map sending hyperbolic triangles to hyperbolic triangles 
and the sum of angles of image triangles around a vertex is 
greater than or equal to $2\pi$.

\begin{thm}\label{thm:GaussB2}
Let $\Sigma$ be a compact hyperbolically-mapped surface relative to 
$v_1, \dots, v_n$ into a general hyperbolic manifold
and suppose that each arc in $\partial \Sigma - \{v_1, \dots, v_n\}$ 
is geodesic.
Let $\theta_i$ be the exterior angle of $v_i$ 
with respect to geodesics in the boundary of $\Sigma$.
Then 
\begin{equation}\label{eqn:GB2}
\Area(\Sigma) \leq \sum_i \theta_i - 2\pi \chi(\Sigma). 
\end{equation} 
\end{thm}

Note that in this part, by a geodesic we mean the geodesic 
with respect to the ambient manifold $M$ unless we state otherwise 
that it is a geodesic in an submanifold, say of codimension-one.

The convex hull of a homotopy-equivalent 
closed subset of a general hyperbolic manifold $M$
is the image in $M$ of the smallest closed convex subset of
containing the inverse image of the closed subset 
in the universal cover $\tilde M$ of $M$.
We can always choose the core 
$\core$ to be a subset of $M^o$.

\begin{thm}\label{thm:convh}
Let $\convh(\core)$ be the convex hull of the core $\core$ of 
a $2$-convex general hyperbolic manifold.
We assume that $\core$ is chosen to be a subset of 
$M^o$ and $\partial \core$ is saddle-imbedded. 
Suppose that $\convh(\core)$ is compact.
Then $\convh(\core)$ is homotopy equivalent to the core and 
the boundary is a truly pleated-triangulated hyperbolic-surface.
\end{thm}

%** outline

In the preliminary, Section \ref{sec:hyp}, 
we recall the definition of 
\cat($\kappa$)-spaces for $\kappa \in \bR$
using geodesics and triangles. 
We also define $M_\kappa$-spaces, the simplicial metric 
spaces needed in this paper. We discuss the link conditions
to check when $M_\kappa$-space is \cat($\kappa)$-space, 
the Cartan-Hadamard theorem, and Gromov boundaries of these spaces.
Next, we discuss the $2$-dimensional versions
of these spaces. Define the interior angles, and prove the Gauss-Bonnet 
theorem. Finally, we discuss general hyperbolic manifolds.

In Section \ref{sec:hmaps}, we show that the universal cover of 
a $2$-convex general hyperbolic 
manifold, which we used a lot in Part 1, is 
a $M_{-1}$-simplicial metric space and a \cat($-1$)-space
and a visibility manifold. 
Next, we define hyperbolic-maps of surfaces. These are similar 
to hyperbolic surfaces as defined by Bonahon, Canary, and Minsky. 
We define Alexandrov nets and A-nets as generalized triangles. 
We show that maps from surfaces can be homotopied to
hyperbolic-maps relative to a collection of boundary points
in $2$-convex general hyperbolic manifolds.
We prove the Gauss-Bonnet theorem for such surfaces 
and find area bounds for polygons.
%We finally discuss $2$- and $3$-simplicies. 

In Section \ref{sec:convhgen}, we discuss the convex hull of 
the core $\core$ in a general hyperbolic manifold.
First, we show that the convex hull and $\core$ is homotopy 
equivalent. We finally show that 
the boundary of the convex hull can be deformed to 
a nearby hyperbolic-imbedded 
surface, which is truly pleated-triangulated.
We show this by finding a geodesic in the boundary of the convex hull
through each point of the boundary.

%%% Oct. 07 10:00

%%% Oct 07 4:00

\section{Hyperbolic metric spaces} \label{sec:hyp}

\subsection{Metric spaces, geodesic spaces, and cat(-1)-spaces.}
\label{sec:metricspaces}
%N (X, d), l, t, t' c, c', p, q, r \epsI

We will follow Bridson-Haefliger \cite{BH}.

%*** 0.3.1. geodesics
Let $(X, d)$ be a metric space. A {\em geodesic path} from 
a point $x$ to $y$, $x, y \in X$ is a map 
$c: [0, l] \ra X$ such that $c(0)=x$ and $c(l) = y$ 
and $d(c(t), c(t')) = |t - t'|$ for all $t, t' \in [0, l]$. 

A {\em local geodesic} is a map $c: I \ra X$ from an interval 
$I$ with the property that for every $t \in I$ there exists $\eps > 0$ 
such that $d(c(t'), c(t'')) = |t'- t''|$ for 
$t', t''$ in the $\eps$-neighborhood of $t$ in $I$. 

$(X,d)$ is a {\em geodesic metric space} if every pair of points of 
$X$ is joined by a geodesic. 

%*** 0.3.2. angles 
We denote by $\bE^2$ the plan $\bR^2$ with the standard Euclidean metric. 
A {\em comparison triangle} in $\bE^2$ 
of a triple of points $(p, q, r)$ in $X$
is a triangle in $\bE^2$ with vertices $\bar p, \bar q, \bar r$ 
such that $d(p,q) = d(\bar p, \bar q)$, 
$d(q, r) = d(\bar q, \bar r)$, and 
$d(r, p) = d(\bar r, \bar p)$, 
which is unique up to isometries of $\bE^2$. 

The interior angle of the comparison triangle at $\bar p$ 
is called the comparison angle between $q$ and $r$ at $p$ 
and is denoted $\bar \angle_p(q, r)$. 

Let $c:[0, a] \ra X$ and $c': [0, b] \ra X$ be two 
geodesics with $c(0) = c'(0)$. 
We define the {\em upper angle} $\angle_{c, c'} \in [0, \pi]$ 
between $c$ and $c'$ to be 
\begin{equation}\label{eqn:angle} 
\angle(c, c') := \limsup_{t, t' \ra 0} \bar \angle_{c(0)}(c(t), c'(t'))
\end{equation}
The angle exists in strict sense if the limsup equals the limit. 

The angles are always less than or equal to $\pi$ by our construction. 
We define angles greater than $\pi$ in 
two-dimensional spaces by specifying sides 
and dividing the side into many parts
(see Subsection \ref{subsec:singhyp}).

%**** 0.3.2.1. Geometric convergence 
We say that a sequence of closed subsets $\{K_i\}$ converge to 
a closed subset $K$ if for any compact subset $A$ of $X$, 
$\{A \cap K_i\}$ converges to $A \cap K$ in Hausdorff sense.

%*** 0.3.3. length metric 
Let $(X, d)$ be a metric. We can define a {\em length-metric} $\bar d$ so that 
$\bar d(x, y)$ for $x, y \in X$ is defined as the infimum of the lengths of 
all rectifiable curves joining $x$ and $y$. 
We note that $d \leq \bar d$ and 
$(X, d)$ is said to be a length space if $\bar d = d$. 

%%% 2:50 02/09/2004
%**** 0.3.3.1. length space, geodesic space 
\begin{prop}[Hopf-Rinow Theorem] 
Let $(X, d)$ be a length space. If $X$ is complete and locally 
compact, then every closed bounded subset of $X$ is compact 
and $X$ is a geodesic space. 
\end{prop} 

%**** 0.3.3.2. Riemannian space  
As an example, a Riemannian space with path-metric is a geodesic 
metric space. A covering space of a length space has an obvious induced 
length metric

%***  quasi-isometry 
%****  quasi-geodesics 

%*** 0.3.4. cat(\kappa)-space 
%**** 0.3.4.0. defn $M_\kappa$.
%NN M_\kappa, \kappa, X, d, x, y, \bar x, \bar y , D_\kappa

We define $M_\kappa$ to be the $3$-sphere of constant curvature $k$, 
Euclidean space, or the hyperbolic $3$-space of constant curvature $k$
depending on whether $\kappa > 0,=0,< 0$ respectively. 

Let $D_\kappa$ denote the diameter of $M_\kappa$ if $\kappa > 0$ 
and let $D_\kappa = \infty$ otherwise. 

Let $(X, d)$ be a metric space. Let 
$\tri$ be a geodesic triangle in $X$ with parameter 
less than $2D_\kappa$  and $\bar \tri$ the comparison
triangle in $M_\kappa$. 
Then $\tri$ is said to satisfy {\em \cat{\rm (}$\kappa${\rm )}-inequality} if 
$d(x, y) \leq d(\bar x, \bar y)$ for 
all $x, y$ in the edges of $\tri$ and their comparison 
points $\bar x, \bar y$, i.e., of same distance from 
the vertices, in $\bar \tri$.
If $\kappa < 0$, a {\em \cat{\rm (}$\kappa${\rm )}-space} is a geodesic space all of whose 
triangles bounded by geodesics satisfy \cat($\kappa$)-inequality. 
If $\kappa > 0$, then $X$ is called a {\em \cat($\kappa$)-space} if 
$X$ is $D_\kappa$ geodesic and all geodesic triangles in 
$X$ of perimeter less than $2D_\kappa$ satisfy 
the \cat($\kappa$)-inequality. 

Angles exists in the strict sense for 
\cat($\kappa$)-spaces if $\kappa \leq 0$. 

A \cat($\kappa$)-space is a \cat($\kappa'$)-space if $\kappa \leq \kappa'$.

A \cat($0$)-space $X$ has a metric $d:X \times X \ra \bR$
that is convex; i.e., given any two geodesics $c:[0,1] \ra X$ 
and $c':[0,1] \ra X$ parameterized proportional to length, 
we have 
\begin{equation}\label{eqn:geoconv} 
d(c(t), c'(t)) \leq (1-t) d(c(0), c'(0)) + td(c(1), c'(1)).
\end{equation} 

%**** 0.3.4.1. characterization of cat(\kappa)-spaces 
%N K, X, \kappa, \lambda, \lambda', L, q, p_{\lambda}, h_{\lambda, \lambda'} 

%***** 0.3.4.1.1. Def of M_\kappa simplicial complex 
A geodesic $n$-simplex in $M_\kappa$ is the convex hull of 
$n+1$ points in general position. 

An {\em $M_\kappa$-simplicial complex} $K$ is 
defined to be the quotient space of the 
disjoint union $X$ of a family of geodesic $n$-simplicies 
so that the projection $q: X \ra K$ induces the injective projection
$p_\lambda$ for each simplex $\lambda$
and if $p_\lambda(\lambda) \cap p_{\lambda'}(\lambda')\ne \emp$, 
there exists an isometry $h_{\lambda, \lambda'}$ from a face of 
$\lambda$ to $\lambda'$ such that 
$p_\lambda(x) = p_{\lambda'}(x')$ if and only if 
$x' = h_{\lambda, \lambda'}(x)$. 

In this paper, we will restrict to the case when locally there 
are only finitely many simplicies, i.e., $X$ is locally convex. 
We do not assume that we have a finite isometry types of 
simplicies as Bridson does in \cite{Bridson}.

A {\em geodesic link} of $x$ in $K$, denoted by $L(x, K)$  is the set of 
directions into the union of simplicies containing $x$. 
The metric on it is defined in terms of angles. 
(For details, see Chapter I.7 of \cite{BH}.)

%***** 0.3.4.1.2 link condition 
\begin{defn}\label{defn:linkcond} 
An $M_\kappa$-simplicial complex satisfies the link condition if 
for every vertex $v$ in $K$, the link complex 
$L(v, K)$ is a \cat($1$)-space. 
\end{defn}

%***** 0.3.4.1.3. link condition prop.
The following theorem can be found in Bridson-Haefliger \cite{BH}:
\begin{thm}[Ballman]\label{thm:linkcond} 
Let $K$ be a locally compact $M_\kappa$-simplicial complex.
If $\kappa \leq 0$, then the following conditions are equivalent\rmc
\begin{itemize} 
\item[(i)] $K$ is a \cat($\kappa$)-space. 
\item[(ii)] $K$ is uniquely geodesic. 
\item[(iii)] $K$ satisfies the link condition and contains no isometrically
imbedded circle. 
\item[(iv)] $K$ is simply connected and satisfies the link condition. 
\end{itemize} 
If $\kappa > 0$, then the following conditions are equivalent \rmc
\begin{itemize}
\item[(v)] $K$ is a \cat($\kappa$)-space. 
\item[(vi)] $K$ is $\pi/\sqrt{\kappa}$-uniquely geodesic. 
\item[(vii)] $K$ satisfies the link condition and contains no
isometrically embedded circles of length less than $2\pi/\sqrt{\kappa}$. 
\end{itemize} 
\end{thm} 
\begin{proof} 
See Ballmann \cite{Ballmann} or Bridson-Haefliger \cite{BH}.
\end{proof}

\begin{lem}\label{lem:2Dlinkcond} 
A $2$-dimensional $M_\kappa$-complex $K$ satisfies the link condition 
if and only if for each vertex $v \in K$, every injective loop 
in $Lk(v, K)$ has length at least $2\pi$. 
\end{lem} 

%**** 0.3.4.2. curvature

\begin{defn} 
A metric space $X$ is said to be of {\em curvature $\leq \kappa$} 
if it is locally isometric to a $\cat(\kappa)$-space, 
i.e., for each point $x$ of $X$, 
there exists a ball which is a $\cat(\kappa)$-space. 
\end{defn}

%***** 0.3.4.2.1. Cartan-Hadamard theorem 
\begin{thm}[Cartan-Hadamard]\label{thm:CaHa}
Let $X$ be a complete metric space. 
\begin{itemize} 
\item[(i)] If the metric on $X$ is locally convex, then 
the induced length metric on the universal cover $\tilde X$
is globally convex. 
{\rm (}In particular, there is a unique geodesic connecting 
two points of $\tilde X$, and geodesic segments in $\tilde X$ vary 
continuously with respect to their endpoints.{\rm )} 
\item[(ii)] If $X$ is of curvature $\leq \kappa$ where $\kappa \leq 0$, 
then $\tilde X$ is a \cat$(\kappa)$-space. 
\end{itemize} 
\end{thm}

%***  0.3.5. $\delta$-hyperbolicity 
Let $\delta$ be a positive real number. 
A geodesic triangle in a metric space 
$X$ is said to be {\em $\delta$-slim} if each of its sides is 
contained in the $\delta$-neighborhood of the union of 
the other two sides. 

For $\kappa < 0$, \cat($\kappa$)-space is $\delta$-hyperbolic.

%*** 0.3.6. Gromov boundary 
For positive real numbers $\lambda, \eps$, 
{\em $(\lambda, \eps)$-quasi-geodesic} in $X$ is a map 
$c: I \ra X$ such that 
\begin{equation}
1/\lambda |t - t'| - \eps 
\leq d(c(t), c(t')) \leq \lambda |t - t'| + \eps  \hbox{ for any pair} t, t' \in \bR
\end{equation} 

Let $X$ be a $\delta$-hyperbolic space. 
Two quasi-geodesic rays $c, c'$ are 
{\em equivalent} or {\em asymptotic} if 
their Hausdorff distance is finite, or, equivalently
$\sup_{t >0} d(c(t), c'(t))$ is finite.
We define the {\em Gromov boundary} $\partial X$  
as the space of equivalence classes of quasi-geodesic
rays in $X$. One can show that $\partial X$ is 
the space of equivalence classes of geodesics rays 
as well. 

%**** 0.3.6.1. topology of the boundary 
If $X$ is a proper metric space, then $X$ is a visibility 
space: For each pair of points $x$ and $y$ in $\partial X$, 
there exists a geodesic limiting to $x$ and $y$. 
Topology and metrics are given on $\partial X$ to 
compactify $X \cup \partial X$. 
The group of isometry acts as homeomorphisms on $\partial X$.

%** 0.4. Singular hyperbolic surface
\subsection{Sigular hyperbolic surfaces}\label{subsec:singhyp} 

%% \theta( ) 
%*** 0.4.1. singular hyperbolic space
%N \Sigma \sng(\Sigma), v, X, x, 

%***** 0.4.1.0.1
A {\em hyperbolic triangle} in a metric space is 
a subset isometric to a triangle in $\hypp$ bounded 
by geodesics. 
Sometimes, we need a {\em degenerate hyperbolic triangle}. 
It is defined to be a straight geodesic segment or a point
where the vertices are defined to be the two endpoints and 
a point, which may coincide. 

A {\em hyperbolic tetrahedron} in a metric space is a subset 
isometric to a tetrahedron in $\hyps$ bounded by 
four totally geodesic planes with six edges geodesic 
segments and four vertices. 
Again degenerate ones can obviously be defined on a hyperbolic triangle, 
a segment, and a point with various vertex and edge structures.

A {\em hyperbolic cone-neighborhood} of a point $x$ in a surface 
$\Sigma$ with a metric is a neighborhood of $x$ which divides into 
hyperbolic triangles with vertices at $x$.
The {\em cone-angle} is the sum of angles of the triangles 
at $x$.  
The set of singular points is denoted by 
$\sing(\Sigma)$ and the cone-angle at $x$ by 
$\theta(x)$. 

By a {\em singular hyperbolic surface}, we mean 
a complete metric space $X$ locally isomorphic to 
a hyperbolic plane or a hyperbolic cone-neighborhood 
with cone-angle $\geq 2\pi$ so that the set of 
singular points are discrete.  
We will also require that $X$ is triangulated by 
hyperbolic triangles in this paper (i.e., is 
a metric simplicial complex 
in the terminology of \cite{BH}). 

%%% Is this the right place? 
By Lemma \ref{lem:2Dlinkcond}, the universal cover $\tilde X$ of $X$ 
is a \cat($-1$)-space. 

\begin{defn}
Let $X$ be a singular hyperbolic surface.
Clearly, $X$ has an induced length metric
and is a geodesic space. 
\begin{itemize}
\item We say that a geodesic in $X$ is {\em straight} if it is 
a continuation of geodesics in hyperbolic triangles 
meeting each other at $\pi$-angles in the intrinsic sense.
\item We can measure angles greater than $\pi$ in 
singular hyperbolic surface by dividing the 
angle into smaller ones. In this case, we need 
to specify which side you are working on. 
In general a path is {\em geodesic} if it is a continuation  
of straight geodesic meeting each other at 
greater than or equal to $\pi$-angles from both sides. 
\item We also say that a boundary point $x$ is {\em bent} if 
the two straight geodesics end at the point
not at $\pi$-angle in the interior. 
\item For $x \in \partial X$, we define 
the {\em interior angle} to be the sum of angles of triangles with 
vertices at $x$ and 
the {\em exterior angle} $\theta(x)$ to be $\pi$ minus the interior angle. 
It could be negative. 
\item We will denote by $\sing(X)$ the set of singular points in 
the interior of $X$ and $\sing(\partial X)$ the 
set of bending points of $\partial X$.
\end{itemize}
\end{defn} 

%*** 0.4.2. Gauss-Bonnet theorem 
%N \theta 

\begin{prop}[Gauss-Bonnet Theorem] \label{prop:GaussBonnet}
Let $\Sigma$ be a compact singular hyperbolic surface
with piecewise straight geodesic boundary. 
Then 
\begin{equation}\label{eqn:GaussBonnet} 
-\Area(\Sigma) + \sum_{v \in \sing(\Sigma)} (2\pi-\theta(v)) 
+ \sum_{v \in \sing(\partial \Sigma)} \theta(v) 
= 2\pi \chi(\Sigma).
\end{equation}
\end{prop}

%**** 0.4.2.1. bigons, closed geodesics
From the Gauss-Bonnet theorem, we can show that
there exists no disk bounded by two geodesics. 
This follows since if such a disk exists, 
then $\theta(v) \geq 2\pi$ for all singular, 
the exterior angles at virtual vertices $\leq 0$,
the exterior angles at common end points $< \pi$,
and the area is less than $0$. 

This implies: 
Given a compact singular hyperbolic space
and a closed curve, we can homotopy the curve into
a closed geodesic, and the closed geodesic is 
unique in its homotopy class. 

%%% corollary?

Moreover, two closed geodesics meet in a minimal number of 
times up to arbitrarily small perturbations: 
that is, the minimum of geometric intersection number 
under small perturbation is the true minimum under 
all perturbations. 
(We may have two geodesics agreeing on an interval 
and diverging afterward unlike the hyperbolic plane.)

%** 0.5. general hyperbolic 3-manifolds 
%***0.5.1. \pi_1(M) and (dev, h), 
%%N M, \dev, h, \tilde M , \vth
\subsection{ General hyperbolic $3$-manifolds} 

By a {\em general hyperbolic manifold}, 
we mean a manifold $M$ with an atlas of charts to 
$\hyps$ with transition maps in $\Isom(\hyps)$. 
The metric on it will be the length metric given 
by the induced Riemannian metric. We require the metric 
to be complete. As a consequence, this is a geodesic space
by local compactness \cite{BH}. 
In general we assume that $\partial M$ is not empty. 
If it is not geodesically complete, 
$M$ need not be a quotient of $\hyps$ which 
are the usual subject of the study in $3$-manifold theory.

Also, we will require general hyperbolic manifolds to have 
hyperbolic triangulations, i.e., it has a triangulation 
so that each tetrahedron is isometric with 
a hyperbolic tetrahedron in $\hyps$. 
Moreover, we assume that the vertices of 
the triangulations are discrete and
the induced triangulation on the universal cover 
map under $\dev$ to a collection of tetrahedra in 
general position in $\hyps$. 
We also require the following mild condition: 
Every boundary point of a general hyperbolic 
manifold has a neighborhood isometric with a subspace of 
a metric-ball in $\hyps$. By subdivisions and small modifications, we can 
always achieve this condition. 

We will say that $M$ is {\em locally convex} if there is
an atlas of charts where chart images are convex subsets of 
$\hyps$. Thus, $M$ is locally convex if $\partial M$ is empty. 
(In this paper, we will be interested in the 
non-locally-convex manifolds.)

Given a general hyperbolic manifold $M$, 
its universal cover $\tilde M$ has an immersion
$\dev: \tilde M \ra \hyps$, which is not in general 
an imbedding or a covering map, and a homomorphism $h$ 
from the deck transformation group $\pi_1(M)$ to 
$\Isom(\hyps)$ satisfying 
\begin{equation}\label{equ:dev} 
\dev\circ \vth = h(\vth) \circ \dev, \vth \in \pi_1(M). 
\end{equation}
$\dev$ is said to be a {\em developing map} and 
$h$ a holonomy homomorphism. 

\begin{thm}[Thurston]\label{thm:localconv} 
Let $M$ be a metrically complete general hyperbolic 
$3$-manifold and is 
locally convex. Then its developing map $\dev$ is 
an imbedding onto a convex domain, and $M$ is 
isometric with a quotient of a convex domain in $\hyps$
by an action of a Kleinian group.
\end{thm}
\begin{proof} 
See \cite{Thnote}.  
\end{proof} 

%***0.5.2. define drilled hyperbolic manifolds 
In this paper, we will often meet
{\em drilled hyperbolic manifolds} obtained 
by removing the interior of 
a codimenion-zero submanifold of a general hyperbolic manifolds. 
They are of course general hyperbolic manifolds.

Of course, special hyperbolic manifolds are general hyperbolic 
manifolds and drilled hyperbolic manifolds.

%***0.5.3. define straight, bent, external angle, 
%length, curvature,

%N S, l_1, l_2, x, l , l'_1

Since a general hyperbolic manifold has a geodesic metric, 
we can define geodesics. 
A {\em straight} geodesic is a geodesic which 
maps to geodesic in $\hyps$ under the charts. 
Geodesics are in general a union of straight geodesics. 
Thus, it has many bent points in general. 
The bent points in the interior of the geodesic segments 
are said to be {\em virtual vertices}.

We define angles as above for metric spaces. 
Then at a virtual vertex, the angle is
equal to $\pi$ since if not, then we can 
shorten the geodesics.

\begin{prop}\label{prop:2Dangle} 
Let $l$ be a geodesic with a bent point $x$ in 
its interior. Let $S$ be a simplicially immersed 
surface containing $l$
in its boundary. Then for $S$ with an induced length metric,
the interior angle at $x$ in $S$ is always greater than 
or equal to $\pi$. 
\end{prop}
\begin{proof} 
If the angle is less than $\pi$, 
we can shorten the geodesic. 
\end{proof} 

Given an oriented geodesic $l_1$ ending at a point $x$ 
and an oriented geodesic $l_2$ starting from $x$, 
we define an {\em exterior angle} between $l_1$ and $l_2$ to be 
$\pi$ minus the angle between the geodesic $l'_1$ 
with reversed orientation and the other one $l_2$. 

%% 02/10/2004 11:30

%*1. 2-convex ghm theory 

\section{ $2$-convex general hyperbolic manifolds 
and {h}-maps of surfaces}\label{sec:hmaps}

\subsection{$2$-convexity and general hyperbolic manifolds}

%%Setting: X domain or general hyperbolic manifold. 
%%N X, \partial X , X^\zerosk, X^\onesk, X^\twosk, \Gamma 

In Part 1, we showed that a general hyperbolic 
manifold was $2$-convex if the vertices of the boundary 
were either saddle-vertices or convex vertices.

%** 1.1. 2-convex domains, manifolds 

We recall the definition of $2$-convexity: 
\begin{defn}\label{defn:2conv}
A general hyperbolic manifold is 
{\em $2$-convex} if given a compact subset $K$ mapping to 
a union of three sides and the interior $T^o$ 
of a tetrahedron $T$ in $\hyps$ 
under a chart $\phi$ of the atlas, there exists a subset $T'$ 
mapping to $T$ by a chart extending $\phi$.
\end{defn} 

\begin{prop}\label{prop:uni2conv} 
If $M$ is a $2$-convex general hyperbolic manifold, 
then $M$ is a $K(\pi_1(M))$, i.e., 
its universal cover is contractible.  
\end{prop} 
\begin{proof}
Since the universal cover $\tilde M$ has an affine 
structure with trivial holonomy induced from 
the affine space containing $\hyps$ from the Klein model, 
this follows from \cite{uaf}. Also,
this follows from Theorem \ref{thm:2convcat}. 
\end{proof}

%*** 1.1.1. Property of 2-convex ghm.
%% c, x, P, L(x, \tilde M), D^c_1, D^c_2, H, H', y

%%Dec 22 4:58
\begin{thm}\label{thm:2convcat} 
Let $M$ be a $2$-convex general hyperbolic manifold. 
Then its universal cover $\tilde M$ is 
a $M_{-1}$-simplicial complex and a \cat$(-1)$-space. 
{\rm (}$M$ has a curvature $\leq -1$.{\rm )}
\end{thm} 
\begin{proof} 
Using Theorem \ref{thm:linkcond} (iv), we need to show that 
for each vertex $x$ of $\tilde M$, the link $P=L(x, \tilde M)$ is 
a \cat($1$)-space. 

To show $P$ is a \cat($1$)-space, we use (vii) of Theorem \ref{thm:linkcond};
i.e., we show that $P$ satisfies the link condition and contains 
no isometric circle of length $< 2\pi$. By the boundary condition on $M$, 
$P$ is isometric to the unit sphere if $x$ is the interior point
or is isometric to a subspace of the unit sphere if $x$ is 
the boundary point. Clearly the former satisfy the 2-dimensional link condition. 

Let $P$ be a proper subspace of a unit sphere $\SI^2$ and $c$ an isometrically 
imbedded circle of length $< 2\pi$. By Lemma \ref{lem:lengthhemisphere}, 
$c$ is disjoint from a closed hemisphere $H$ in $\SI^2$. 

The closed curve $c$ meets $\partial P$ since otherwise $c$ has 
to be a great circle of length $2\pi$ being a geodesic. 
As $c$ may never cross-over the circle
$\partial P$, let $D^1_c$ and $D^2_c$ denote the disks 
in $\SI^2$ bounded by $c$. Then $\partial P$ is a subset of 
$D^1_c$ or $D^2_c$. Assume without loss of generality that 
the former is true. 

Since $c$ is disjoint from the hemisphere $H$, it follows that 
$H$ is a subset of $D^1_c$ or $D^2_c$. 
In the second case, $H \subset P$. Looking at this situation, from the vertex $x \in M$ again, we see that 
$2$-convexity is violated since we can find a triangle in $M$ containing $x$ in its interior
whose one-sided neighborhood with $x$ removed is in the interior of $M$.

Suppose that $H$ is a subset of $D^1_c$. Let $H'$ be the complementary 
open hemisphere. Then $c \subset H'$ and $\partial P$ is 
outside the disk $D^2_c$ in $H'$ bounded by $c$. 
Since $H'$ has a natural affine structure, and $c$ is compact, it follows that 
the convex hull $K'$ of $c$ in $H'$ is compact. Let $y$ be an extreme point of 
$K$, where $y \in c$ as well. In the one-sided neighborhood of $y$ inside $c$, there are no 
points of $\partial P$ implying that we can shorten $c$ in $P$
contrary to the fact that $c$ is isometrically imbedded.
\end{proof} 

\begin{lem}\label{lem:lengthhemisphere} 
Let $\gamma$ be a broken geodesic loop in the sphere $\SI^2$ of 
radius $1$. If the length of $\gamma$ is less than $2\pi$, 
then there exists an open hemisphere containing it 
{\rm (}and hence a disjoint closed hemisphere{\rm ).}
\end{lem}
\begin{proof} 
We can shorten the loop without increasing the number of 
broken points to a loop as short as we want. 
A sufficiently short loop is contained in an open hemisphere. 

Let $l_t, t\in [0,1]$ be a homotopy so that 
$l_1$ is the original loop and $l_0$ is a constant loop. 
Then let $A$ be the maximal connected set containing $0$ 
so that $l_t$ for $t \in A$ is contained in an open hemisphere, say $H_t$. 

The set $A$ is an open set since the small change in $l_t$ does not 
violate the condition. Suppose that the complement of $A$ is not empty, and 
let $t_0$ be the greatest lower bound of 
the complement of $A$. 
Then $l_{t_0}$ is contained in a closed hemisphere, say $H'$, since 
we can find a geometric limit of the closure of $H_t$s. 

Suppose that $\partial H' \cap l_{t_0}$ is contained in 
a subset of length strictly less than $\pi$. 
Then we can rotate $H'$ along a pivoting antipodal 
pair of points on $\partial H'$ outside the subset. 
Then the new hemisphere contains $l_{t_0}$ in its interior, 
a contradiction. 

Suppose that $\partial H' \cap l_{t_0}$ contains a pair of antipodal points. 
Let $s_1$ and $s_2$ be the corresponding points of $[0,1]$
and suppose $0 < s_1< s_2 < 1$ without loss of generality. 
Then two arcs $l_{t_0}|[s_1, s_2]$ and $l_{t_0}| [s_2, 1] \cup [0, s_1]$,
must have length greater than or equal to $\pi$, a contradiction. 
Therefore, no subsegment in $\delta H'$ of length $\leq \pi$ contains 
$\partial H' \cap l_{t_0}$. 

Suppose now that there are three points $p_1, p_2, p_3$
in $\delta H' \cap l_{t_0}$ are not contained in 
a subsegment in $\delta H'$ of length $\leq \pi$ 
and no pair of them are antipodal. 

The sum of lengths of segments $\ovl{p_1p_2}, \ovl{p_2p_3}, \ovl{p_3p_1}$ 
equals $2\pi$. 
This is clearly less than or equal to that of $l_{t_0}$ since 
the shortest arcs connecting the pairs 
$(p_1, p_2)$, $(p_2, p_3)$, $(p_3, p_1)$ are these segments respectively. 
This is again a contradiction. 

Thus $A$ must be all of $[0,1]$. 
\end{proof}

The following proves Theorem \ref{thm:Gromovh} in detail.
%**** 1.1.1.1. Geodesic curvature property of M
\begin{prop}\label{prop:Mgeodcond} 
Let $\tilde M$ be a universal cover of 
a compact $2$-convex general hyperbolic 
manifold $M$. Then the following hold{\rm :} 
\begin{itemize}
\item $\tilde M$ is uniquely geodesic.
\item Geodesic segments of $\tilde M$ depend continuously on their endpoints. 
\item The metric is locally convex. 
\item $\tilde M$ is $\delta$-hyperbolic and hence it is a visibility manifold. 
\item $M$ has curvature $\leq -1$. 
\item Given any path class on $M$, there exists a unique geodesic 
segment, which depends continuously on endpoints.
\end{itemize}
\end{prop}
\begin{proof} 
These are direct consequences of the fact that 
$\tilde M$ is a \cat$(-1)$-space. 
\end{proof}

%% June 30 9:30
%*** 1.1.2. hyperbolic-maps into 2-convex manifolds Defs
\subsection{Hyperbolic-maps of surfaces into $2$-convex general
hyperbolic manifolds}

A {\em triangulated hyperbolic surface} is a metric surface with 
or without boundary triangulated and each triangle is 
isometric with a hyperbolic triangle or a degenerate 
hyperbolic triangle in $\hypp$. 
A {\em half-space} of $\hyps$ is a subset bounded by 
a totally geodesic plane.  

\begin{defn}\label{defn:h-map} 
Let $\Sigma$ be a compact triangulated hyperbolic surface, $M$ 
a general hyperbolic $3$-manifold, and $f:\Sigma \ra M$
a map which sends each triangle to a hyperbolic triangle in $M$. 
Let $\partial \Sigma$ have distinguished vertices 
$v_1, \dots, v_n$. 
Then $f$ is a {\em hyperbolic-map relative to} $\{v_1, \dots, v_n\}$ 
if the sum of the angles of the image triangles of 
the stellar neighborhood of each interior vertex 
$v$ is greater than or equal to $2\pi$ and 
the sum of angles of the image triangles of 
the stellar neighborhood of the boundary vertex $v$, $v \ne v_i$, 
is greater than or equal to $\pi$. 
\end{defn}  

A hyperbolic-map is a completely analogous concept to a hyperbolic-map by
Bonahon \cite{Bonahon}, Canary and Minsky and so on.
Note that if the boundary portion between 
$v_i$ and $v_{i+1}$ is geodesic for each $i$, 
then the boundary angle conditions are satisfied also.

%**** 1.1.2.1 Define A-net, a-net

%\begin{defn}\label{defn:arcvertex}
%A point on an arc is a {\em vertex} if every 
%open interval containing it is not geodesic. 
%A point on an arc is a {\em virtual vertex} if every
%open interval containing it is not straight geodesic.
%\end{defn}

\begin{defn}\label{defn:A-net}
Given an arc or a point $\alpha$ and an arc $\beta$ in $M$, 
an {\em Alexandrov net with ends $\alpha$ and $\beta$}
is a map $f: I \times I \ra M$
so that $s \in I \mapsto f(t, s)$ is geodesic for 
each $s$ and $t \mapsto f(t,0)$ is $\alpha$ and 
$t \mapsto f(t, 1) \in \beta$.
\end{defn}

\begin{lem}\label{lem:closegeo} 
Let $M$ be a $2$-convex general hyperbolic manifold.  
Let $\gamma$ be a geodesic in $M$. Then for any geodesic $\gamma'$
sufficiently close to $\gamma$,
there exists a homotopy $H: I \times I \ra M$
so that the following hold {\rm :} 
\begin{itemize}
\item $s \mapsto H(0,s)$ is $\gamma$ 
and $s \mapsto H(1, s)$ is $\gamma'$. 
\item $H$ is a simplicial map with 
a triangulation of $I\times I$ with 
all vertices at $\{0,1\} \times I$. 
\end{itemize}
\end{lem} 
\begin{proof} 
For each virtual vertex of $\gamma$, we choose 
a real number $\eps> 0$ such that the $\eps$-neighborhood of 
the vertex is a stellar neighborhood. 
If $\gamma'$ is in an $\eps$-neighborhood of $\gamma$ 
for $\eps>0$ for any such $\eps$s, 
then we can find the desired $H$. 
\end{proof}

\begin{def}\label{defn:a-net} 
Given a point or an arc $\alpha$ and another arc $\beta$,
an {\em A-net} $f:I \times I \ra M$ with ends $\alpha$ and 
$\beta$ is a map such that 
\begin{itemize}
\item $s \mapsto f(t_i, s)$ for a finite subset 
$\{t_1=0, t_2, \dots, t_n=1\}$ of $I$ is a geodesic 
for each $i$,
\item $t \mapsto f(t, 0)$ is $\alpha$ 
and $t \mapsto f(t, 1)$ is $\beta$. 
\item $f$ is a hyperbolic-map relative 
to vertices of the arcs $\alpha$ and $\beta$
with a triangulation of $I\times I$ with all the vertices in 
$\{t_1, \dots, t_n\} \times I$.   
\end{itemize}
\end{def}

\begin{prop}\label{prop:a-netexist} 
Given a point or an arc $\alpha$ and another arc $\beta$,
there exists an A-net with ends $\alpha$ and $\beta$.
\end{prop}
\begin{proof}
We find an Alexandrov net $f: I \times I \ra M$ 
with ends $\alpha$ and $\beta$. 
We take sufficiently many $t_i$'s so that 
geodesics $s \mapsto f(t_i, s)$ are very close. 
By Lemma \ref{lem:closegeo}, we can 
find a simplicial map $F: I \times I \ra M$. 
Since $s \mapsto F(t_i, s) = f(t_i, s)$ are 
geodesics, the sum of angles at each of 
the sides of a vertex on this geodesic 
is greater than $\pi$. Hence, the sum of 
angles at an interior vertex is greater than or 
equal to $2\pi$. At the vertices of $s \mapsto F(0, s)$ 
or $s \mapsto F(1, s)$, the sum of angles are 
greater than $\pi$. Therefore, $F$ is a hyperbolic-map.
\end{proof} 

%**** 1.1.2.2 Existence of hyperbolic-maps rel boundary rel closed geodesics 
%%N, v, v_n, v_m, v_l

\begin{prop}\label{prop:himmexists} 
Let $\Sigma$ be a compact triangulated surface, $M$ a general 
hyperbolic $3$-manifold, and let $f:\Sigma \ra M$ be a map
with an injective induced homomorphism $f_*: \pi_1(\Sigma) \ra \pi_1(M)$. 
\begin{itemize}
\item Let $v_1, \dots, v_n$ be the distinguished vertices in $\partial \Sigma$
and $l$ be a union of disjoint simple closed curve in $\Sigma$
which is disjoint from $\{v_1, \dots, v_n\}$ 
and is a component of $\partial \Sigma$ or is disjoint from $\partial \Sigma$.
\item We suppose that $\{v_1,\dots, v_n\} \cup l$ is not empty.  
Suppose that $f$ maps each arc in $\partial \Sigma$ 
connecting two distinguished vertices to a geodesic
and each component of $l$ or $\partial \Sigma$ without 
any of $v_1, \dots, v_n$ to a closed geodesic. 
\end{itemize} 
Then in the relative homotopy class of $f$ with 
$f|\partial \Sigma$ fixed, 
there exists a hyperbolic-map $f': \Sigma' \ra M$
relative to $v_1, \dots, v_n$ where $\Sigma'$ is 
$\Sigma$ with a different triangulation in general and 
$f'$ agrees with $f$ on $\partial \Sigma \cup l$. 
\end{prop}

From now on, we will just use $v_i$ for $f(v_i)$ and so on
since the reader can easily recognize the difference. 
By the angle of a triangle,
we mean the corresponding angle measured in the image triangle of $f$. 

\begin{proof} 
First, we find a topological triangulation $\Sigma$ 
so that all the vertices are in the union of 
$\{v_1, \dots, v_n\} \cup l \cup \partial \Sigma$. 
We find a geodesic in the right path-class 
for each of the edges of the triangulations.
For each triangle, we extend by choosing 
a vertex and the opposite geodesic edges 
and finding A-nets with these ends. 

At each interior point of an edge, 
we see that the sum of angles of any of its side is greater 
than or equal to $\pi$ since the edge is geodesic. 
Since A-nets are hyperbolic-maps, we see that 
the whole map is a hyperbolic-map. 
\end{proof}

%*** 1.1.3. Gauss-Bonnet theorem for hyperbolic-maps
%N \area, \chi
%**** 1.1.3.1. Prop. \area(S) <= sum ext - 2\pi \chi(S)
\subsection{Gauss-Bonnet theorem for hyperbolic-maps}

\begin{prop}\label{prop:GB} 
Let $\Sigma$ be a compact hyperbolically-mapped surface relative to $v_1, \dots, v_n$. 
Let $\theta_i$ be the exterior angle of $v_i$ 
with respect to geodesics in the boundary of $\Sigma$.
Then 
\begin{equation}\label{eqn:GB} 
\Area(\Sigma) \leq \sum_i \theta_i - 2\pi \chi(\Sigma). 
\end{equation} 
\end{prop}
\begin{proof} 
The interior angle with respect to $\Sigma$ is larger than 
the angle in $M$ itself. 
Thus the exterior angle with respect to $\Sigma$ is 
smaller than the exterior angle in $M$. 

Since the interior vertices have the angle sums greater than or 
equal to $2\pi$ and the boundary virtual vertices have 
the angle sum greater than or equal to $\pi$, 
the proposition follows from the Gauss-Bonnet theorem.
\end{proof}

%**** 1.1.3.2. Corollary. \area(triangle) <= \pi \area(square) <= 2\pi 
An $n$-gon is a disk with boundary a union of geodesic segments
between $n$ vertices. 

\begin{cor}\label{cor:ngons} 
Let $S$ be a hyperbolically-mapped $n$-gon. Then 
$\Area(S) \leq (n-2)\pi$. 
\end{cor}
\begin{proof} 
The exterior angle of a bent virtual vertex on a geodesic
is always less than $\pi$. 
\end{proof}

\section{Convex hulls in $2$-convex general hyperbolic manifolds}
\label{sec:convhgen}

%** 2.0. Setting
Let $M$ be a $2$-convex general hyperbolic manifold 
with finitely generated fundamental group, and
$\core$ denote a core of $M$.

%** 2.1. Definition of the convex hull 

Let $\tilde M$ be the universal cover of $M$. 
Since $\core \ra M$ is a homotopy equivalence 
the subset $\tilde \core$ in $\tilde M$ which is
the inverse image of $\core$ is connected and 
is a universal cover of $\core$. 
A subset of $\tilde M$ is {\em convex} if any two 
points can be connected by a geodesic in the subset.

The {\em convex hull} $\convh(K)$ of a subset $K$ of $\tilde M$ 
is the smallest closed convex subset containing $K$. 
%We define the {\em convex hull} $\convh(\tilde \core)$ of the core 
%$\core$ in $M$ as the smallest closed convex 
%set containing $\tilde \core$ in $\tilde M$:
Since $\tilde \core$ is deck-transformation group invariant,
and the convex hull is the smallest convex subset,  
$\convh(\tilde \core)$ is deck-transformation group invariant.
Therefore, $\convh(\tilde \core)$
covers its image. We define the image as 
$\convh(\core)$, i.e., $\convh(\tilde \core)$ quotient by 
the deck-transformation group action.

Since $\core$ is a $3$-dimensional domain, 
$\convh(\core)$ is a $3$-dimensional closed set. 

%%% Dec 23, 12:10
%*** 2.1.1. The homotopy property of the convex hull. 

\begin{prop}\label{prop:homprop} 
The convex hull $\convh(\core)$ of the compact core $\core$ of 
$M$ is homotopy equivalent to $\core$. 
\end{prop}
\begin{proof}
Let $\tilde \core$ be the inverse image of $\core$
in the universal cover $\tilde M$ of $M$. 
Then $\tilde \core$ and $\tilde M$ are both contractible 
as $M$ and $\core$ are irreducible $3$-manifolds. 

A closed curve in the convex hull $\convh(\tilde \core)$ of 
$\tilde \core$ bounds a disk in $\convh(\tilde \core)$ since a distinguished point on
the curve can be connected by a geodesic in any other point 
of the curve. Similarly, a sphere always bounds a $3$-ball. 
Therefore, $\convh(\tilde \core)$ is contractible. 
\end{proof}

%The following definition shows how to define the convex hull 
%interms of the domain covers. 

%\begin{prop}\label{prop:hatconv}
%Suppose that $M$ is covered by a domain $\hat M$ in
%$\hyps$. Let $\hat \core$ be the inverse image of $\core$ in
%the domain cover of $\hat M$ of $M$.
%Define $\convh(\hat \core)$ as the smallest closed 
%convex subset of $\hat M$ containing $\hat \core$. 
%Then under the cover $\tilde M \ra \hat M$, 
%$\convh(\tildec \core)$ covers $\convh(\core)$. 
%\end{prop}
%\begin{proof}
%The cover $\hat M$ is of form $\tilde M/\Gamma$ 
%for a deck transformation group $\Gamma$.
%Under the covering map $\tilde M \ra \hat M$, 
%$\convh(\tilde \core)$ covers $\convh(\hat \core)$: 
%Since the inverse image of $\convh(\hat \core)$ in $\tilde M$ 
%is convex, $\convh(\tilde \core)$ maps into $\convh(\hat \core)$ 
%under the covering map. 
%Since $\convh(\tilde \core)$ is $\Gamma$-invariant, 
%it follows that $\convh(\tilde \core)$ covers its image, say $B$.
%Since any two points $p, q$ of $B$ can be connected by 
%a geodesic segment by taking inverse points of $p$ and $q$, 
%it follows that $B$ is convex. 
%Therefore, $B = \convh(\hat \core)$. 
%\end{proof}

%** 2.2. The structure of the boundary of the convex hull

%*** 2.2.1. hyperbolicity of the boundary of the convex hulls 

A surface is {\em pleated} if through each point of 
it passes a straight geodesic. 

Recall that the pleated-triangulated surface 
is an imbedded surface where a closed subdomain is 
a union of a locally finite collection of 
totally geodesic convex domains meeting 
each other in geodesic segments and the complementary open 
surface is pleated. 

A pleated-triangulated surface is {\em truly pleated-triangulated} 
if the triangulated part are union of totally geodesic polygons
in general position.

A {\em truly pleated-triangulated hyperbolic-surface}
is a truly pleated-triangulated surface where 
each vertex of the triangles is a hyperbolic-vertex.

\begin{lem}\label{lem:geohvertex}
If a geodesic in $M$ contained in $S$ 
passes through a vertex in 
the triangulated part of $S$, then the vertex is a hyperbolic-vertex.
\end{lem}
\begin{proof} 
A neighborhood of a point of the triangulated part is 
stellar. If a geodesic passes through, 
the angles in both sides are greater than or equal to $\pi$:
otherwise, we can shorten the geodesic. 
Hence, the sum of the angles is greater than or 
equal to $2\pi$. 
\end{proof} 

%\begin{lem}\label{lem:pleatinglines2} 
%A straight pleating line $l$ in the boundary of a convex hull in $M^o$
%is extendable to a pleating line but does not go inside the interior of 
%the convex hull. 
%\end{lem} 
%\begin{proof}
%Using convex hull of an open ball neighborhood of $x$, 
%we find that $l$ is in the interior of the convex hull.
%\end{proof}

The following proves Theorem \ref{thm:convh}:

\begin{prop}\label{prop:convbdhsurf} 
Let $K$ be a deck-transformation-group invariant codimension $0$
submanifold of $\tilde M$ with $\partial K$ saddle-imbedded. 
Also, suppose $K$ is a subset of $\tilde M^o$.
The boundary of $\convh(K)$ can be given the structure of
a convex truly pleated-triangulated hyperbolic-surface.
\end{prop} 
\begin{proof}
We will show that 
\begin{itemize} 
\item through each point of $\partial \convh(K)$
a geodesic in $\partial \convh(K)$ passes or 
\item the point is in the triangulated part and is a saddle-vertex 
or a hyperbolic-vertex. 
\end{itemize}

Let $x$ be a boundary point of $\convh(K)$:
Suppose that $x$ is a point of the manifold-interior of $\tilde M$.
Take a ball $B_\eps(x)$ in the interior for a sufficiently 
small $\eps$. Then $\convh(K) \cap B_{\eps}$ is 
the convex hull of itself. Since $B_{\eps}$ is 
isometric with a small open subset of $\hyps$, 
the ordinary convex hull theory shows that 
there exists a geodesic in the boundary of
the convex hull through $x$: If not, 
we can find a small half-open ball to decrease the 
convex hull as the side of the half-open ball cannot
meet $\partial K$ by the saddle-imbeddedness of $\partial K$.

Suppose that $x$ is in the topological interior of 
$\convh(K)$ but in $\partial \tilde M$. 
There exists a neighborhood of $x$ in $\convh(K)$ with manifold-boundary 
in $\partial \tilde M$. If $x$ is in the interior of 
an edge or a face of $\partial \tilde M$, 
then there is a geodesic through $x$ obviously. 
Suppose that $x$ is a vertex of $\partial \tilde M$. 
$x$ can be a saddle-vertex or a convex vertex
(see Proposition \ref{prop:s-bd2-conv}).
\begin{itemize}
\item If $x$ is a convex vertex of $\partial \tilde M$, we can find a truncating 
totally geodesic hyperplane and a sufficiently 
small disk in it bounding 
a neighborhood of $x$ in $\tilde M$. 
Since $K$ is disjoint from the disk, 
we see that $x$ is not in the convex hull.
This is absurd. 
\item If $x$ is a saddle-vertex of $\partial \tilde M$, then $x$ is 
a saddle-vertex of $\partial \convh(K)$. 
\end{itemize} 

Assume from now on that $x$ is a point in the topological boundary of 
$\convh(K)$ and on $\partial \tilde M$. This means that $x$ is in the frontier
of the open surface $C = \partial \convh(K) - \partial \tilde M$. 

(a) Suppose $x$ is a point of the interior of a triangle $T$ in 
$\partial \tilde M$. The set $T \cap\convh(K)$ is a convex 
subset and $x$ lies in the boundary. The boundary must 
be a geodesic since we can use a small half-open ball 
to decrease the convex hull otherwise. 
Hence there is a geodesic through $x$. 

(b) Suppose that $x$ is a point of the interior of an edge
in $\partial \tilde M$. We take a small ball $B_\eps(x)$ around 
$x$, which is isometric with a ball in $\hyps$ of 
the same radius and a wedge removed. The line $l$ of the wedge 
passes through $x$. 

Let $P_1$ and $P_2$ be the totally geodesic plane extended 
in $B_\eps(x)$ from the sides of the wedge.
We denote by $P'_1$ the set $\partial B_\eps(x) \cap P_1$ 
and $P'_2$ the set $\partial B_\eps(x) \cap P_2$. 
We can form two convex subsets $L_1$ and $L_2$ in $B_\eps(x)$ 
that are the closures of the components of 
$B_\eps(x) - P_1 - P_2$ and adjacent to $P'_1$ and $P'_2$ 
respectively.  

The set $L_1 \cap \convh(K)$ is a convex subset of $L_1$
and $L_2 \cap \convh(K)$ one of $L_2$. 
The open surface $C$ may intersect $L_1$ or $L_2$ or both. 

If $C$ is disjoint from $L_1$, then 
it maybe that a one-sided neighborhood of $x$ in a triangle in 
$\partial \tilde M$ is a subset of $\convh(K)$ 
and an edge of the triangle is a geodesic through $x$.
(The side $P_1 - P_1^{\prime, o}$ of $B_\eps(x)$ 
is in $\convh(K)$.) 
Otherwise, $\convh(K)$ is contained in a convex subset 
of $B_\eps(x)$ bounded by $P_1$. In this case, 
the ordinary convexity in $\hyps$ holds and there is a geodesic 
in $\partial \convh(K)$ through $x$.

%%%% Alternative saddle-vertex or convex or concave vertex....

Since the same argument holds with $L_2$ as well,
we assume that $C\cap L_1$ and $C \cap L_2$ are both not empty. 

If $C\cap L_1$ or $C \cap L_2$ are totally geodesic surfaces, 
then $x$ is on a pleating locus that is the edge of the wedge. 

We may assume without loss of generality that $C \cap L_1$ is 
not totally geodesic in a ball of radius $\eps$ about $x$ for every 
sufficiently small $\eps > 0$.
Then there exists a sequence of points 
$\{x_i \in L_1 \cap \tilde M^o\}$ 
converging to $x$ and a sequence of
pleating lines 
\[\{l_i \subset \partial \convh(K) \cap \tilde M^o\}\]
so that $x_i \in l_i$. 
By Lemma \ref{lem:pleatinglines}, 
we only have to worry about the case when all $l_i$s 
end at $x$. In this case there exists a small 
neighborhood $B(x)$ of $x$ such that $C \cap L_1 \cap B(x)$
is a stellar set with vertex at $x$. 

By a same argument, $C\cap L_2 \cap B(x)$ is a stellar set also
with a vertex at $x$. 
Considering $C \cap L_1 \cap B(x)$ and $C\cap L_2 \cap B(x)$ at the same time,
in order that at $x$, the convexity 
to hold true and $x$ to be in $\partial \convh(K)$, 
we see that $C\cap L_1 \cap B(x)$ has to have a unique pleating geodesic with 
a convex dihedral angle as seen from $\convh(K)$ 
and so does $C \cap L_2 \cap B(x)$. 
Furthermore, their unique pleating geodesics must extend each other as geodesics
passing through $x$.

%%% Dec 23 1:30

(c) Now assume that $x$ in $\partial \convh(K)$ is 
a vertex of $\partial \tilde M$. 

Let $B_\eps(x)$ be a small neighborhood of $x$ so 
that $B_\eps(x) \cap \tilde M$ is a stellar set from $x$. 
As before $x$ is in the boundary of $C$. 

Suppose first that there are no pleating lines with a sequence of points 
on them converging to $x$. We can choose a small 
$\eps$ so that $B_\eps(x) \cap \convh(K)$ is a stellar set.

Let $M'$ be an ambient general hyperbolic manifold containing $M$
in its interior which is homeomorphic to the interior of $M$ 
as there are always such a manifold.
We claim that $x$ is a saddle-vertex of 
$B_\eps(x) \cap \partial \convh(K)$: 
If not, we can find a small half-open ball $B$ in $\tilde M'$ with 
a totally-geodesic side passing through $x$ with
$B^o$ disjoint from $\partial \convh(K)$. 
By stellarity, $B^o$ is disjoint from $\convh(K)$ and 
we can decrease $\convh(K)$ if $x \not\in K$, which is a contradiction. 
If $x \in K$, then there is no such $B$ as $x$ is a saddle vertex of $K$ itself. 
Therefore, $x$ is a saddle-vertex.

We assume that $B_\eps(x) \cap \convh(K)$ is not a stellar set.
Suppose now that there exists a sequence of points 
$\{x_i \in l_i\}$ converging to $x$ where $l_i$ is 
a distinct pleating line for each $i$ and does not end at $x$.
Here, $l_i$ are infinitely many.
Lemma \ref{lem:pleatinglines} shows that 
the endpoints of $l_i$ are bounded away from $x$. 
Therefore, a subsequence of $l_i$ converges to a geodesic $l$ passing
through $x$. 

We proved the two items above, and 
$\partial \convh(K)$ is a pleated-triangulated surface.

As a final step, we show that 
$\partial \convh(K)$ is a truly pleated-triangulated 
surface: Let $A$ be the closure of the 
set of all points in $\partial \convh(K)$ intersected with
the interiors of triangles in $\partial \tilde M$.
Then $A$ is a locally finite union of totally geodesic polygons and 
segments. The complement of $A$ in $\partial K$ is pleated 
since it lies in the interior of $\tilde M$. 
Any pleated lines in $\tilde M^o$ must end at $A$ or 
is infinite. By Remark \ref{rem:pleatinglocus}, 
the set of pleating lines ending at $A$ is isolated. 

%As before, 
%we see that a pleating line ending at a point of $A$ 
%must be a finite length segment. 
%Take a union of these finite length segments with $A$ 
%o form $A'$. The infinite length pleating lines 
%are in the finite union of minimal laminations 
%bounded away from $A'$ by Proposition \ref{prop:pleatinglocus}.
%Since the minimal laminations are bounded away from $A'$, 
%we obtain a truly pleated-triangulated surface.

By Lemma \ref{lem:geohvertex}, it is a hyperbolic-surface as well.
The convexity is obvious. 
\end{proof}

\begin{lem}\label{lem:pleatinglines} 
Let $l_i$, $i \in I$, be a collection of 
mutually distinct straight pleating lines 
$\partial \convh(K) - \partial \tilde M$ 
for a convex hull $\convh(K)$ of a closed subset $K$ of 
$\tilde M$ and a countably-infinite index set $I$. 
Suppose $x_i \in l_i$ form a sequence converging to $x$
but $x$ is not on $l_i$s. 
Then the endpoints of $l_i$s are bounded away from $x$
and a subsequence of $l_i$ converges to a line segment in the pleating 
locus containing $x$ in its interior. 
\end{lem}
\begin{proof}
Suppose that the endpoints of $l_i$ are bounded away from $x$. 
Then the second statement holds obviously.

Suppose that the endpoints $q_i$ of $l_i$
form a sequence converging to $x$. 
Then we may assume without loss of generality that 
$q_i$ lies in an arc or a point $\alpha$ in a triangle in $\partial \tilde M$.
If the arc $\alpha$ is a convex curve, we can decrease $\convh(K)$ 
further, a contradiction. Thus  $\alpha$ is a geodesic
or a point.

By Lemma \ref{lem:pleatinglocus}, $\alpha$ cannot be a line.
If $\alpha$ is a point, then $\alpha \ne x$, and the conclusion holds. 
\end{proof}

\part{The proof of the tameness of hyperbolic $3$-manifolds}

\section{Outline of the proofs}
%** outline of the proof

We will prove Theorems A and B in this part: 
%*** exhaustions
The strategy is as follows.
Suppose that the unique end $E$ of $M$ not associated 
with incompressible surface is not geometrically finite and 
is not tame. 
We find an exhausting sequence $M'_i$ in $M$ 
so that $M'_i$ contains neighborhoods of all tame and geometrically 
finite ends and meets the neighborhood of the infinite 
end in a compact subset.

%*** closed geodesics and compressions
\begin{description}
\item[Step 1] Using the work of Freedman-Freedman \cite{FrFr}, we can
modify $M'_i$s to be compression bodies:
Since $E$ is not geometrically finite, we can choose 
a sequence of closed geodesics $\cl_i \ra \infty$. 
We fix a sufficiently small Margulis constant $\eps$.
We assume without loss of generality that $\cl_i \subset M'_i$.
%and that the Margulis tubes that $M'_i$ meets are contained in 
%$M'_i$. 
Let $\mu_i$ be the union of closed geodesics that
are in the Margulis tubes in $M_i$.
We further modify $M'_i$ so that $\partial M'_i$ is incompressible
in $M - \cl_1 - \cdots \cl_i - \mu_i$ with 
the compact core $\calC$ removed. 

The manifold $N_i$ is obtained from compressing disks 
for $M'_i$ in $M -\cl_1 - \cdots -\cl_i -\mu_i - \calC$.  
Let $A_i$ be a homotopy in $M$ between $\cl_i$ and 
the closed curve in $\core$, which can be homotopied to be in $N_i$. 

Now we modify $N_i$ to $M''_i$ so that $\cl_i \subset M''_i$ 
and $M''_i$ are compression bodies and the boundary 
component of $M''_i$ corresponding to $E$ is incompressible 
in $M$ removed with the closed geodesics $\cl_1 \dots, \cl_i$
and $\calC$ and 
contains any Margulis tubes that $M''_i$ meets.
(See Subsection \ref{subsec:step1}.)

%*** 2-convex hulls
\item[Step 2] As the boundary is incompressible, 
we take a $2$-convex hull $M_i$ of $M''_i$ using crescents
(see Part 1).
This implies that $M_i$ is a polyhedral hyperbolic space 
and hence is \cat$(-1)$. 
Since $M_i$ is isotopic to $M''_i$, the homotopy 
$A_i$ between $\cl_i$ and a closed curve in $\core$
still exists. 
%Also, we will be able to 
%make $M_i$ homeomorphic to compression bodies.

%*** Margulis constant
We show that we can choose a Margulis constant 
independent of $i$ and the thin part of $M_i$ are 
contained in the original Margulis tubes of $M$
and are homeomorphic to solid tori with nontrivial 
homotopy class in the original tubes. 
(See Subsection \ref{subsec:step2}.)
%*** covers
\item[Step 3] 
We take the cover $L_i$ of $M_i$ corresponding to the fundamental group of 
the fixed compact core $\core$ of $M$. Since $M_i$ is tame, the cover $L_i$ is 
shown to be tame (this is from ideas of Agol).

%*** boundary of the convex hull
The core $\core$ lifts to the cover $L_i$ and can be 
considered a subset. We take a convex hull $K_i$ of 
$\core$ in $L_i$, which is shown to be compact.
Since $K_i$ is homotopy equivalent to $\core$ by Proposition 
\ref{prop:homprop} in Part 2, 
the boundary component $\Sigma_i$ of $K_i$ 
corresponding to $E$ has the same genus as that of 
a boundary component of $\core$. 
$\Sigma_i$ is a ``hyperbolic surface"
(see Part 2).
Since $\cl_i$ is an exiting sequence, and an $\eps$-neighborhood of $M_i$ contains 
$\cl_i$, it follows that 
$p_i|\Sigma_i$ is an exiting sequence of surfaces. 
This proves Theorem A. 
(See Subsection \ref{subsec:step3}.)

%*** drillout and convex hulls
\item[Step 4] We push $\core$ inside itself so that $\core$ does not meet 
$\partial K_i$. 
We now remove the core from $K_i$ to obtain $K_i - \core^o$.
We can find a simple closed curve $\alpha$ in 
$\Sigma_i$ compressible in $K_i$. 
We realize $\alpha$ by a closed geodesic 
$\alpha^*$ in $K_i - \core^o$. 
Using it, as Bonahon does \cite{Bonahon}, we obtain a simplicial 
hyperbolic surface $T_i$ meeting $\partial \core$.

Then by compactness of bounded simplicial
hyperbolic surfaces of Souto \cite{Souto}, infinitely many immersed $T_i$s 
are isotopic in $M - \core^o$ and hence infinitely many of $p_i|\Sigma_i$ are 
isotopic. Since $p_i|\Sigma_i$ are exiting
and are isotopic in $M - \core^o$,
this implies that the end $E$ is topologically tame, proving Theorem B.
(See Subsection \ref{subsec:step4}.)
%This is a contradiction. 

\end{description}

%But we should point out that we had similar 
%ideas as early as Agol but without much details. 

%** Proof itself

\section{The Proof of Theorem A} 

%*** Setting exhaustions other ends

%*** boundary of the core 
\subsection{Choosing the right exhausting sequence} \label{subsec:step1}

\subsubsection{Choosing the core}
As a preliminary step, we choose the compact core more carefully so that 
$\partial \core$ is saddle-imbedded:
We choose incompressible closed 
surfaces $F_i$ associated with incompressible 
ends $E_i$ to be strictly saddle-imbedded by Theorem C and disjoint 
from one another (see Remark \ref{rem:svertices}).
We choose a number of closed geodesics in 
$E_i$ and choose a mutually disjoint submanifold homeomorphic to 
$F_i \times I$ disjoint and between these curves for each $i$.
Then by Theorem C, we find a mutually disjoint collection of
manifolds in the respective neighborhoods of $E_i$
between these curves whose boundary components 
are strictly saddle-imbedded.

%%Oct 16 6:30

Essentially $\partial \core$ is considered as  
a regular neighborhood of the union of the essential surfaces 
$F_1, \dots, F_n$ and a number of arcs connecting them 
in some manner. 

We choose each of the arcs to be the shortest 
path in $M$ among the arcs connecting 
the surfaces $F_i$s with the respectively given homotopy classes.
Their endpoints must be in the interior of 
an edge of a triangle. 
By perturbing $F_i$ if necessary, we may 
assume that they are all disjoint geodesics. 
We first take thin regular neighborhoods of $F_i$s 
which are triangulated. 
We take thin regular neighborhoods of the geodesics
which are triangulated and all of whose vertices lie in
$F_i$s. 

We take the union of the regular neighborhood of these 
geodesics with those of $F_i$s to be our core $\core$. 
We may assume that $\partial \core$ is strictly saddle-imbedded as well.
(We may need to modify a bit where the neighborhoods meet.)
As stated above, we choose $\core$ to be in the interior 
$M^o$. Obviously, if necessary, we push $\core$ inward 
itself without violating strict saddle-imbeddedness of $\partial \core$. 

%%%%% Need that crescent move does not go over the core....

%**** other ends
\subsubsection{Choosing a compression body exhaustion}
Let $M$ be as in the introduction, and let $U_1, \dots, U_n$ 
be mutually disjoint neighborhoods of incompressible ends
$E_1, \dots, E_n$. Suppose that the end $E$ is a geometrically 
infinite but not geometrically tame. 

Let $\hat M$ be the $2$-convex hull of 
$M$ with $U_1, \dots, U_n$ removed.
The boundary components $F_1, \dots, F_n$ 
corresponding to $U_1, \dots, U_n$ 
of $\hat M$ are saddle-imbedded respectively.

%**** exhaustion by compression bodies
Let $M'_i$ be an exhaustion of $M$
by compact submanifolds containing $\hat M$. 
We extend $M'_i$ by taking a union with $U_1, \dots, U_n$
so that $M'_i$ meets neighborhoods of $E$ in 
compact subsets or in the empty set. 
We assume that $M'_i$ contains 
the boundary components $F_1, \dots, F_n$ and 
contains the core $\core$ of $M$ always. 

%***** Mtubes and disks
\begin{lem}\label{lem:mtubes}
A disk with boundary outside the union of 
Margulis tubes may be isotopied with the boundary 
of the disk fixed so that 
the intersection is the union of meridian disks. 
\end{lem}
\begin{proof}
First, we perturb the disk to obtain transversal intersection. 
A disk may meet the Margulis tubes in 
a union of planar surfaces. 
The boundary of the Margulis tube meets the disk 
in a union $C$ of circles. If an innermost component 
is outside the tube, then since the boundary tube is 
incompressible in $M$ with the interior of 
the Margulis tubes removed, it follows that we 
can isotopy it inside. This means that 
the innermost components are disks. 

If an innermost component of $C$ is a circle bounding 
a disk in the boundary of the Margulis tube, 
then we can isotopy the bounded disk away from the tube. 
Now, each component of $C$ is a meridian circle. 
\end{proof}

We fix a small Margulis constant $\eps_M >0$.

\begin{prop}[Freedman-Freedman, Ohshika]
\label{prop:compressing}
We obtain a new exhaustion $M'_i$ so that 
each $M'_i$ is homeomorphic 
to a compression-body with incompressible
boundary components removed. $M'_i$ contains 
the Margulis tubes that $M'_i$ meets
by taking union with these.
\end{prop}
\begin{proof} 
We essentially follow Theorem 2 of Freedman-Freedman \cite{FrFr}. 
We assume that each $M'_i$ includes any Margulis tubes it meets.
$\partial M$ has no incompressible closed surface other than 
ones parallel to $F_1, \dots, F_n$. 
Hence, we compress the boundary component 
$\partial M'_i$ until we obtain a union of 
balls and manifolds homeomorphic to 
$F_i$ times an interval. 
An exterior disk compression adds a disk times $I$ to the compressed 
manifold from $M'_i$ but an interior disk compression removes 
a disk times $I$ from the manifold.
The exterior disk may meet Margulis tubes outside 
$M'_i$. By Lemma \ref{lem:mtubes}, if the disk meets a Margulis tube 
at an essential disk, then we include the Margulis tube.
If not, we push the disk off the Margulis tube. 
This operation amounts to adding $1$-handles to 
the manifolds.

We recover our loss to $M'_i$ by interior disk compressions
by attaching $1$-handles each time we do 
interior compressions. The core arcs of $1$-handles 
may meet the exterior compression disk many times. 
We add a small neighborhood of the cores first. 
Then we isotopy to make it larger and larger to recover 
the loss due to interior disk compression
while fixing the Margulis tubes outside $M_i$. 
(This may push the exterior compression disks.)
We also recover all the Margulis tubes originally in $M'_i$.

From the surface times the interval components, we 
add all the $1$-handles to obtain the desired compression body.

\end{proof} 

\begin{rem}\label{rem:intext}
We do interior compressions first and 
then exterior compressions. This is sufficient
to obtain the union of $3$-balls and $F_i$ times 
the intervals. The reason is that we can isotopy 
any interior compressing curve away from the traces of 
exterior compression disks.
\end{rem}

%%% July 21, 2004 11:30

%***** closed geodesc s Margulis tubes
Since $E$ is geometrically infinite, there exists 
a sequence of closed geodesics $\cl_i$ tending to $E$
by Bonahon \cite{Bonahon}.
We assume that $\cl_i \subset M'_i$ for each $i$
since $M'_i$ is exhausting. Let $\calC_i$ denote the union of 
$\cl_1, \dots, \cl_i$.
We assume without loss of generality that 
$M'_i$ has a free-homotopy $A_i$ between $\cl_i$ and a 
closed curve in $\core$ since $M'_i$s form an exhausting sequence.
Let $\mu_i$ be the union of the simple closed geodesics 
in the Margulis tubes that $M'_i$ contains. 

%***** Recall McCollough
\begin{rem}\label{rem:diskslide}
We know that given a surface there exists 
a finite maximal collection of exterior essential compressing 
disk so that any other exterior essential compressing 
disks can be pushed inside the regular neighborhood of their union and 
some arcs on the surface connecting them. 
This is from the uniqueness of the compression body.  
(See Theorem 1 in Chapter 5 of McCullough \cite{McNote}.)
\end{rem}

%**** remove core and compress boundary N_i
\subsubsection{Compression bodies $M''_i$ with ``incompressible" boundaries}

\begin{description}
\item[$N_i$ with ``incompressible'' boundaries]
If we remove the interior $\core^o$ of the core from $M'_i$, 
and compress the boundary 
$\partial M'_i$ in $M - \core^o - \calC_i - \mu_i$ 
and remove resulting cells to obtain a manifold 
with incompressible boundary containing $\partial \core$. 
Then we join the result with $\core$, $\calC_i$ and $\mu_i$.
Let us call the resulting $3$-manifold $N_i$. 
Note that $N_i$ need not be a compression body 
and may not form an exhausting sequence. 
However, $N_i$ contains $\core, \calC_i, \mu_i$
for each $i$.

The exterior compressing disk of $M'_i$ 
may meet the Margulis tubes outside $M'_i$. 
We may assume that $N_i$ meets these Margulis 
tubes in meridian disks times intervals
by Lemma \ref{lem:mtubes}.

%**** homotopy A_i
The manifold $N_i$ is obtained from compressing disks 
for $M'_i$. Let $A_i$ be a homotopy in $M'_i$ between $\cl_i$ and 
a closed curve in $\core$. Then a compressing disk 
for the sequence of manifolds obtained from $M'_i$ 
by disk-compressions
in $M - \core^o - \calC_i - \mu_i$ may meet $A_i$ in 
immersed circles. Since the circles bound 
immersed disks in the compressing disk, and $\cl_i$ 
is not null-homotopic in $M$, we can modify $A_i$ so
that $A_i$ has one less number of components where 
$A_i$ meets the compressing disk. In this manner, 
we can find $A_i$ in $N_i$.

%**** Our operations M''_i
\item[$N_i^{III}$]
We do the following steps:
\begin{itemize} 
\item We find a maximal collection of 
essential interior compressing disks for $N_i$
and do disk compressions. By incompressibility,
the disk must meet one of $\core, \calC_i, \mu_i$ essentially.
We call $N^I_i$ the component of the result containing $\partial \core$. 
\item We find a maximal collection of 
essential exterior compressing disk for the result of the first step
and do disk compressions. We call $N^{II}_i$ the component 
containing $\partial \core$. 
\item We add $1$-handles lost in the first step.
We call the result $N^{III}_i$. 
\end{itemize} 

%%%¤ÔAugust 23 11:30 very careful here....

Clearly $N^{III}_i$ includes $N_i$.

\end{description}

\begin{prop}\label{prop:incompressible2} 
The submanifold $N^{III}_i$ is homeomorphic to a compression body 
and is contained in a compression body $M''_i$ whose 
boundary component is incompressible in 
$M - \core^o - \calC_i - \mu_i$. A Margulis tube is 
either contained in $M''_i$ or the tube meets $M''_i$ 
in meridian disk times intervals. 
\end{prop}
\begin{proof} 
The fact that $N^{III}_i$ is a compression body follows as before. 
Using the fact that 
$N^{III}_i$ is contained in some compression body $M'_j$ for a large $j$,
we will now show that $N^{III}_i$ is contained in
a compression body $M''_i$ with the above property.

By construction of $N^{III}_i$, it follows that 
an interior compression disk can be isotopied inside the 
regular neighborhood of the union
the disks of the $1$-handles and arcs in the boundary.
Therefore, each interior compressing disk must 
intersect at least one of $\core, \calC_i, \mu_i$ essentially.
Hence $\partial N^{III}_i$ has no interior compression 
disk in $M - \core^o - \calC_i - \mu_i$. 

There could be an exterior compression disk for $N^{III}_i$. 
We take a maximal mutually disjoint 
family $D_1, \dots, D_n$ of them where no two of 
$\partial D_i$ are parallel.
We choose $j$ sufficiently 
large so that a compression body $M'_j$ includes all of them
as $M'_i$s form an exhausting sequence.

We find a $3$-manifold $X$ isotopic to $N^{III}_i$ in $M'_j$:
$M'_j$ decomposes into a union of cells or submanifolds homeomorphic to
$F_l$ times intervals by a maximal family of interior compression 
disks. We suppose that no two of the disks are parallel
and call $D$ the union of these disks.  

We consider $X$ to be a thin regular neighborhood of 
a $1$-complex with the unique vertex in a fixed base cell $\calB$ of $M'_j$, 
fix a handle decomposition of $X$ corresponding to the $1$-complex 
structure, define the complexity of the imbedding of $X$ in $M'_j$ 
by the number of components of $X \cap D$,
choose $X$ with minimal complexity, and put all things in general positions.
For each disk $D_k$, we first get rid of any closed circles 
by the innermost circle arguments. 
We may find an edgemost arc if 
$\partial D_k$ meets $D$ bounding a component of 
$D_k -D$. Then there must be a handle of 
$X$ following the arc in $\partial D_k$. 
This handle can be isotopied away,
and then using the innermost circle argument again 
if necessary, we can reduce the complexity. 
Therefore, it follows that 
each $D_k$ is in the fixed base cell of $M'_j$. 
Also, a handle where $D_k$ passes essentially
must lie in the base cell also.

We look at the handles in the base cell $\calB$  and 
the disks. The union of the handles and the ball around 
the basepoint is a handle body.
Our disks $D_1, \dots, D_n$ are in the cell. 

From Corollary 3.5 or 3.6 of Scharlemann-Thomson \cite{ScTh}, 
we see that there exists an unknotted cycle in the $1$-complex 
or the $1$-complex has a separating sphere. 
In the first case, we cancel the corresponding cycle 
by an exterior compression. In the later case, 
a sphere bounds a ball which we add to $X$, i.e., we engulf it. 
We do the corresponding topological operations to 
$N^{III}_i$ while $X$ and $N^{III}_i$ are still compression bodies. 
In both cases, we can reduce the genus of 
the boundary of the handle body $X$ or $N^{III}_i$.
We do this operations until there are no more 
exterior compressing disks. (The genus complexity shows 
that the process terminates.)
We let the final result to be $M''_i$.

The imbedding $\partial M''_i \ra M - M^{\prime \prime, o}_i$ 
is incompressible since the boundary of any exterior compressing disk 
for $\partial M''_i$ can be made to avoid the traces of 
handle-canceling exterior compressions which are pairs of disks
or the disks from the $3$-ball engulfing.
Therefore these must be exterior compressing disks for 
$N^{III}_i$. 

The imbedding $\partial M''_i \ra M^{\prime \prime, o}_i 
- \core^o - \calC_i - \mu_i$ 
is incompressible since the boundary of any interior compressing disk
can be made to avoid the traces of handle-canceling 
exterior compressions.

The statement about Margulis tubes follows by Lemma \ref{lem:mtubes}.
\end{proof}

%**** The crescent isotopy to $M_i$

\subsection{Crescent-isotopy}\label{subsec:step2}
We will now modify $M''_i$ by crescent-isotopy.

\begin{lem}\label{lem:crescentcore} 
A secondary highest-level crescent of 
$\tilde \Sigma$ does not meet the interior of $\tilde \core$. 
\end{lem}
\begin{proof} 
If not, then $\tilde \core^o$ meets $I_\calR$ for a secondary 
highest-level crescent $\calR$. We may assume that $\calR$ is 
compact by an approximation inside. 
Again find a Morse function by totally geodesic planes 
parallel to $I_\calR$. $\core \cap \calR$ has a maximum 
inside as $\core \cap \calR$ is compact. 
But at the maximum point, a totally geodesic plane 
bounds a local half open ball disjoint from $\core^o$. 
This contradicts the saddle-imbeddedness of $\partial \tilde \core$. 
\end{proof}

We define $M_i$ to be the $2$-convex hull of $M''_i$.
The boundary components of $M_i$ are saddle-imbedded.
Let $\core'$ be the core obtained from $\core$ by pushing $\partial \core$
inside $\core$ by an $\eps$-amount.
Note that during the crescent-isotopy
$\core'$ is not touched by the interior of secondary highest 
level crescents since 
$\partial \core$ is saddle-imbedded by Lemma \ref{lem:crescentcore}.

Define $M_i^\eps$ be the regular $\eps$-neighborhood of $M_i$.
$\calC_i, \mu_i \subset M_i^\eps$ 
by Proposition \ref{prop:avoidg} in Part 1.
We may assume $A_i$ is in $M_i^\eps$ since we isotopied $M''_i$ 
to obtained $M_i$.

%***** thin parts
The universal cover $\tilde M_i$ of $M_i$ 
with the universal covering map $p_i$
is an $M_{-1}$-space and is $\delta$-hyperbolic. 

We define the {\em thin part} of $M_i$ as the subset of $M_i$ 
where the injectivity radius is $\leq \eps$. 
Since $\tilde M_i$ is a uniquely geodesic,  
through each point of the thin part of $M_i$
there exists a closed curve of length $\leq \eps$
which is not null-homotopic in $M_i$.

\begin{prop}\label{prop:Margulis}
The $\eps$-thin part of $M_i$ is homeomorphic to a disjoint union of 
solid tori in Margulis tubes in $M$ parallel to 
a multiples of the shortest geodesics in the respective tubes.
Furthermore, we can choose $\eps >0$ independent of $i$. 
\end{prop}
\begin{proof}
Let $\gamma$ be a closed curve of length $\leq \eps$ 
which is not null-homotopic in $M_i$. 
Then if $\gamma$ has nontrivial holonomy, 
then $\gamma$ lies in a Margulis tube of $M$ 
which either is in $M_i$ or disjoint from $M_i$. 

Suppose that 
$\gamma$ is null-homotopic in $M$. 
Then $\gamma$ bounds an immersed disk $D$ in $M$. 
The diameter of $D$ is $\leq \eps$.

Suppose that $D$ cannot be isotopied into $M_i$.
By incompressibility of $\partial M_i$, 
$D$ must meet $\cl_i$ or $\mu_i$ or $\core$. 
Also, $D/M_i$ must be nonempty.
There must be an innermost disk $D'$ such that  
$\partial D'$ maps to $\partial M_i$ and 
$D'$ maps into $M_i$ and 
meets $\mu_i$, $\core$ or $\calC_i$
and a component of $D/\partial M_i$ 
adjacent to $D'$ lies outside $M_i$.

If $D'$ meets $\core$, then the diameter of $\partial D'$ 
is not so small, and hence that of $D$. 
If $D'$ meets $\calC_i$, then since $A_i$ is in $M_i$ 
$D'$ cannot be bounded outside by a component 
in $M /M_i^o$. 

Suppose that an innermost disk $D'$ meet $\mu_i$. 
Since $\partial D'$ is very close to $\mu_i$
due to its size, and 
the distance from $\mu_i$ to $\partial M''_i$ 
is bounded below by a certain constant, 
it follows that the boundary of $D'$ lies 
in the union of $I$-parts of some crescents
or its perturbed images
obtained during the crescent-isotopies. 
Since the length of components of $\mu_i$ is short, 
we see that the the $I$-parts meeting $\partial D'$
would extend for long lengths along the geodesics
near a component of $\tilde \mu_i$. 
Therefore, we see that at the last stage of
the isotopies, we have the inverse image of 
torus bounding a component of $\mu_i$. 
Since our crescent moves are isotopies and $\partial M'_i$ is 
not homeomorphic to a torus, this is a contradiction. 

We conclude that $\gamma$ bounds a disk in $M_i$. 
Since $\gamma$ is not null-homotopic in $M_i$, 
this is a contradiction. 
Therefore, the thin part of $M_i$ is in the intersection of 
the Margulis tubes of $M$ with $M_i$. 

Note here that the Margulis constant $\eps >0$ could be 
chosen independent of $i$ as the above argument 
passes through once $\eps$ is sufficiently small 
regardless of $i$.
\end{proof}

%***** Characterization of the thin parts
%We will denote by $T_d$ the Margulis tube containing
%the simple closed geodesic, say $d$.
By above discussions, it follows that any $\eps$-short closed 
curve in $M_i$ is a multiple of the simple closed geodesic 
in a  Margulis tube. We may assume that given a component of 
the thin part of $M_i$,
an $\eps$-short closed curve of a fixed homotopy class passes  
through each point of the component.
Therefore, components of thin parts are solid tori in Margulis 
tubes. 

During the crescent move, the shortest closed 
geodesic in the Margulis tube may go outside 
particularly during the convex perturbations.
But there are short closed curves in the result homotopic 
to the closed geodesic.
Therefore the thin part are union of solid tori 
parallel to some multiples of the shortest geodesics.

%*** compactness of the convex hulls
\subsection{Covers $L_i$s}\label{subsec:step3}
Assume without loss of generality that 
we have an inclusion map $i:\core \ra M_i$ for each $i$. 
Let $L_i$ be $\tilde M_i/i_*\pi_1(\core)$
with the covering map $p_i:L_i \ra M_i$. 
$L_i$ has ends corresponding to $F_1, \dots F_n$, 
and another end $E$ corresponding to $E$. (We abused notation 
a little here.)

\begin{prop}\label{prop:compK_i}
The convex hull $K_i$ of $\core$ in $L_i$ is compact.
\end{prop} 
\begin{proof}
Since $M_i$ is tame, its cover $L_i$ is tame. 
For any compact set, we can find a compact core 
of $L_i$ containing it. By choosing a large 
compact subset of $M_i$, we obtain a compact 
core $\core'$ of $L_i$ containing it which is obtained 
as the closure of the appropriate component of $L_i$ with 
a finite number of disks removed. 

Certainly $\core$ is a subset of it. 
The disks lifts to disks in the universal cover $\tilde L_i$ 
of $L_i$. They bound the universal cover $\tilde \core'$ of $\core'$. 

The convex hull of a disk is a compact subset of 
$\tilde L_i$ since the convex hull of a compact 
subset is compact in the universal cover.
Since the convex hull of $\tilde \core'$ is in the union 
of $\tilde \core'$ and convex hulls of the boundary disks,
the convex hull of $\core$ itself is compact
being a subset of the union of a compact set $\core'$ and 
finitely many compact sets.  
\end{proof}

Since $\core$ is homotopy equivalent to $L_i$, 
there exists a convex hull $K_i$ of $\core$ in $L_i$ homotopy 
equivalent to $\core$ by Proposition \ref{prop:homprop} in Part 2.
Obviously, $K_i$ contains $F_1, \dots, F_n$. 
Let $\Sigma_i$ be the unique boundary component of 
$\partial K_i$ associated with $E$. 

%**** increasing nature of $K_i$.
Any closed geodesic homotopic to 
a closed curve in $\core$ in $L_i$ is contained in $K_i$: 
If not, we can find a hyperbolic-imbedded annulus $B_i$ 
with boundary consisting of  
the closed geodesic and a closed curve on $\partial K$, 
intrinsically geodesic,  
where the interior angles in $B_i$ are always greater than or equal to 
$\pi$. Such an annulus clearly cannot exist.  

We can find a quasi-geodesic $c'_i$ $\eps$-close to $\cl_i$ 
in $M_i$. Since $c'_i$ is homotopic to a closed curve in 
$\core$ by a homotopy $A'_i$ in $M_i$ modified from $A_i$,
we have that $c'_i \subset L_i$.
Choose a geodesic representative $c''_i$ in $L_i$
which is again arbitrarily close to $c'_i$ and hence to $\cl_i$.
Therefore $c''_i \subset K_i$ for each $i$.

%Since the distance from $c_i$ to $\core$ in $M$ is tending to $\infty$, 
%it follows that the maximal distance from $p|\Sigma_i$ to $\core$ 
%in $M_i$ is tending to $\infty$. 
Since $\Sigma_i$ is a truly pleated-triangulated convex hyperbolic-surface,
the intrinsic metric in the pleated part is a Riemannian 
hyperbolic ones. Thus, $\Sigma_i$ carries a triangulated 
hyperbolic-surface structure intrinsically. 
Since $\Sigma_i$ is intrinsically a hyperbolic-imbedded surface
and $c''_i$ forms an exiting sequence, 
$p|\Sigma_i$ is one also and hence form an exiting sequence in $E$. 

More precisely, the parts of boundary of the image of $K_i$ 
form an exiting sequence in $E$. Any part of the boundary of 
the image of $K_i$ is in the image of $\Sigma_i$. 
Hence, there exists an exiting sequence of parts 
of $\Sigma_i$. By the uniform boundedness of $\Sigma_i$, 
it follows that $\Sigma_i$ form an exiting sequence in $E$. 

%%Remark 03/25 2:30
\begin{rem}\label{rem:e-thinh-surf}
The uniform nature of the Margulis constant plays 
a role here. Any $\eps$-thin part of a hyperbolic-immersed surface 
must be inside a Margulis tube in $L_i$ and by 
incompressibility the thin parts are homeomorphic 
to essential annuli. Since $\Sigma_i$ is incompressible 
in $L_i - \core^o$, we see that the thin part of $\Sigma_i$ 
is a union of essential annuli which are not homotopic 
to each other. Thus, outside the Margulis tubes, 
$\Sigma_i$s have bounded diameter independent of $i$. 
\end{rem}

%%% March 29, 5:15 2005

%** Proof of Theorem B
\subsection{The Proof of Theorem B} \label{subsec:step4}

%*** remove core
We recall that $\core$ was pushed inside itself somewhat so that 
$\Sigma_i$ and $\partial K_i$ does not meet $\core$. 
In $K_i$, we may remove $\core^o$ and 
we obtain a compact submanifold $Q_i$ of 
codimension $0$ bounded by saddle-surfaces
including $\Sigma_i$ since $\partial \core$ is saddle-imbedded. 
$Q_i$ is a $2$-convex \cat$(-1)$-space. Finally,
$\Sigma_i$ is incompressible in $Q_i$
since any compressing disk of $\Sigma_i$ not meeting the core 
would reduce the genus of $\Sigma_i$ but the genus of 
$\Sigma_i$ is the same as that of the component of $\partial \core$
corresponding to the end $E$ since $K_i$ is homotopy equivalent 
to $\core$.

As $K_i$ is homeomorphic to a compression body,
we choose a compressing curve $\alpha$ in $\Sigma_i$.
Then $\alpha$ bounds a disk $D$ in $K_i$
and the core $\core$ must meet $D$ in its interior. 
Let $\hat \alpha$ be the geodesic realization of 
$\alpha$ in $K_i - \core^o$, which must be in $K_i$. 

If $\hat \alpha$ does not meet $\partial \core$, 
then it maps to a geodesic in $M$, which is absurd since 
the holonomy of $\alpha$ is the identity. 
Therefore, $\hat \alpha$ meets $\partial \core$. 

We form a triangulation of $\Sigma_i$ with the only vertex $p$
at a point of $\alpha$ and including $\alpha$ as an edge. 
Then choosing a vertex $\hat p$ in $\hat \alpha$ 
and a path from $p$ to $\hat p$, we isotopy each edge of 
the triangulation to a geodesic loop in $K_i - \core^o$ 
based at $\hat p$. Each triangle is isotopied to an A-net
spanned by new geodesic edges. The resulting 
surface $T_i$ is a hyperbolic-surface since each of the triangles 
is an A-net and a geodesic passes through each 
point of the $1$-complex.

%**** Finite isotopy classes.
Each surface $q_i: T_i \ra M - \core^o$
has the same genus and is homotopic to 
$p_i| \Sigma_i$ in $M - \core^o$. 
Since they are hyperbolic-imbedded, and meet $\partial \core$, 
they are in a bounded neighborhood of $\core$ by the boundedness of hyperbolic-imbedded 
surfaces. They form a pre-compact sequence. 
Thus infinitely many of $q_i|T_i$ are isotopic 
in $M - \core^o$. Therefore, infinitely many of $p_i|\Sigma_i$ 
are isotopic in $M-\core^o$. 
Since $\Sigma_i$ bounds larger and larger domains in a cover of $M$ 
and $\Sigma_i$ projects to a surface far from $\core$, 
the above fact shows that our end $E$ is 
tame as shown by Thurston \cite{Thnote}.
(This is essentially the argument of Souto \cite{Souto}
simplifiable in our setup.)

%*Bibliography

\bibliographystyle{plain}
%\bibliography{all}

\end{document}